\documentclass{article}

\usepackage{amsmath, amssymb, amsthm}
\usepackage{graphicx}
\usepackage{enumerate}
\usepackage{tikz}

\newcommand{\prf}[1][]{\noindent {\bf Proof{#1}:} }
\newcommand{\Ex}{\mathbb{E}}
\newcommand{\Var}{\mathsf{Var}}
\newcommand{\Cov}{\mathsf{Cov}}
\newcommand{\diam}{\mbox{diam}}
\newcommand{\cl}{\mbox{cl}}

\newtheorem{theorem}{Theorem}
\newtheorem{prop}[theorem]{Proposition}
\newtheorem{lemma}[theorem]{Lemma}
\newtheorem{cor}[theorem]{Corollary}
\theoremstyle{remark}
\newtheorem{remark}[theorem]{Remark}

\title{Testing for linearity in boundary regression models with application to maximal life expectancies}
\author{J\"urgen Kampf, Alexander Meister}

\begin{document}

\maketitle

\begin{abstract}
We consider a regression model with errors that are a.s.\ negative. Thus the regression function is not the expected value of the observations but the right endpoint of their support. 
We develop two goodness-of-fit tests for the hypotheses that the regression function is an affine function, study the asymptotic distributions of the test statistics in order
to approximately fix the sizes of the tests, derive their finite-sample properties based on simulations and apply them to life expectancy data. 
\end{abstract}

\section{Introduction}
Over the last years there has been an increasing interest in \emph{boundary regression models}, i.e.\ in models of the form
\[ Y_i=g(x_i)+\epsilon_i, \quad i=1,\dots, n,\]
where $Y_i$ is the observed data, $g:\mathbb{R}\to\mathbb{R}$ is the (unknown) regression function, $x_i\in[0,1]$ are the design points and the errors $\epsilon_i, \, i=1, \dots, n,$ are i.i.d.\ random variables with $\mathbb{P}(\epsilon_i<0)=1$, but $\mathbb{P}(\epsilon_i < -\delta) < 1$ for all $\delta>0$. The dependence of $x_i$ and $\epsilon_i$ on $n$ is suppressed in the notation. Nonparametric estimators for the function $g$ have been proposed and investigated by  \cite{HvK09}, \cite{JMR14} and \cite{DNS17} and a goodness-of-fit test for the hypotheses that $g$ equals a fixed function $g_0$ has been developed in \cite{ReWa17}. 

Tests for the hypotheses that the true regression function belongs to a given set of functions are a classical topic in nonparametric mean regression analysis, see e.g.\ \cite{HuKu15} for an overview. Such tests have also been considered in quantile regression, see e.g. \cite{Zh98}, \cite{DLF19} or \cite{MvKY19}. Except for the test from \cite{ReWa17} which works only if the null hypotheses is a single function no such test has been proposed in boundary regression. See \cite{GC13} for a survey of general goodness-of-fit tests in regression analysis. 

The aim of the present paper is to develop a test for
\[ H_0: \mbox{$g$ is an affine function}\quad \mbox{vs.}\quad H_1: \mbox{$g$ is not affine} \]
in boundary regression. 

This endeavor is essential motivated by an observation of Oeppen and Vaupel \cite{OV02} in the field of demography. They considered the maximal life expectancy in each year, where the maximum is taken over all states in the world. They observed an ``extraordinarily linear'' growth in time. Still the test we construct will be able to reject the null hypotheses with highest significance. 

As the test statistic we use the squared $\mathcal{L}^2$-distance between a non-parametric estimator $\hat g$ of $g$ and the space of affine functions, enriched with a bias reduction method and a split of the sample technique. It will be shown to be asymptotically Gaussian with mean zero under the null hypotheses and the limiting variance can either be bounded by elementary calculations or calculated exactly using a Poisson approximation. These two possibilities yield two tests based on the same test statistic.

This paper is organized as follows. In Section \ref{s:methodology} we define the tests precisely. The asymptotic distribution of the test statistic will be calculated in Section \ref{s:asymptotic} both under the null hypotheses and under the alternative hypotheses. Simulation results together with the details on the application to the life expectancy data will be presented in Section \ref{s:application}. In Section \ref{s:discuss} we discuss our results and some further questions. The proofs are deferred to the appendix.

\section{Methodology}\label{s:methodology}

Here we construct our test statistic and propose the critical values for the two tests. 

A first step in constructing the test statistic is estimating the regression function $g$. We use the non-parametric estimator $\hat g$ from \cite{HvK09}. This estimator works as follows. At first some bandwidth $h>0$ is chosen. Then for each $x\in[0,1]$ we define $\hat g(x)$ as the minimal $p_0$ such there is $p_1$ with 
\begin{equation} Y_{2i} \le p_0 + p_1(x_{2i}-x)\label{Y_above} \end{equation}
for all $i$ with $|x-x_{2i}| < h$. We will not be concerned with boundary effects in this paper, but we will assume that there are some observation points to the left and to the right of the interval $[0,1]$ so that this construction makes sense for all $x\in[0,1]$. Requiring \eqref{Y_above} only for data points with even indices is a preparation for the split-of-the-sample technique as will become clear later. 

As preliminary test statistic consider the squared discrete $\mathcal{L}^2$-distance from $\hat g$ to the set of affine functions,
\begin{equation} 
T_1= \min_{m,c}\sum_{i=1}^{n/2} \big( \hat g(x_{2i-1})-mx_{2i-1}-c\big)^2, 
\label{e:Tdef}\end{equation}
where we assume $n$ to be even. 

At first we derive an explicit form of $T_1$. We abbreviate
\[
R := \sum_{i=1}^{n/2} x_{2i-1}^2\quad \mbox{and}\quad 
S := \sum_{i=1}^{n/2} x_{2i-1},\] 

\begin{prop}\label{p:T1_explicit}
Let $f:\mathbb{R} \to \mathbb{R}$ be any affine function.

\begin{enumerate}[(i)]
\item  An explicit form of $T_1$ is given by 
\[ T_1= \begin{aligned}[t]
\sum_{i=1}^{n/2} \big(\hat g(x_{2i-1})-f(x_{2i-1})\big)^2  &- \frac{2}{n} \cdot \Big( \sum_{i=1}^{n/2} \big(\hat g(x_{2i-1}) - f(x_{2i-1}) \big) \Big)^2\\
& - \frac{(n/2) \cdot \Big(\sum_{i=1}^{n/2} \big(\hat g(x_{2i-1}) -f(x_{2i-1}) \big) \cdot \big(x_{2i-1}-\frac{S}{n/2}\big)  \Big)^2}{Rn/2-S^2}.\end{aligned} \]
\item For equidistant design points $x_i = \frac{i}{n}, i=1, \dots, n,$ we obtain
\begin{align}
T_1 = \sum_{i=1}^{n/2} \big(\hat g(x_{2i-1})-f(x_{2i-1})\big)^2 & - \frac{2}{n}\Big( \sum_{i=1}^{n/2} \big(\hat g(x_{2i-1}) - f(x_{2i-1}) \big) \Big)^2\nonumber\\
& - \frac{n \cdot \Big(\sum_{i=1}^{n/2} \big(\hat g(x_{2i-1}) -f(x_{2i-1}) \big) \cdot \big(x_{2i-1}-\tfrac{1}{2}\big) \Big)^2}{\frac{1}{24}n^2-\frac{1}{6}}. \label{e:T_explicit}\end{align}
\end{enumerate}
\end{prop}

In using $T_1$ as the test statistic there is the following problem. Both of our tests will be based on the fact that
\[ \frac{T-\Ex[T]}{\sqrt{\Var(T)}} \to \mathcal{N}(0,1), \quad n\to\infty, \]
in distribution. Now $\Ex[T]$ and $\Var(T)$ depend on the unknown distribution of the errors $\epsilon_i$ and hence they can neither be calculated nor simulated. So these numbers have to be bounded or approximated by other quantities $E$ and $V$ which can be calculated or simulated. If we want to set $T=T_1$, these quantities have to fulfill
\begin{equation} 
\limsup_{n\to\infty}\frac{|E-\Ex[T_1]|}{\sqrt{\Var(T_1)}} <\infty \qquad \mbox{and} \qquad \liminf_{n\to\infty} \frac{V}{\Var(T_1)} >0\label{e:asy_equiv} \end{equation}
in order to ensure that
$ (T_1-E)/\sqrt{V}$
is still bounded in probability. However, $\Ex\, T_1\in O(n^{-1}h^{-2})$, while $\Var\, T_1 \in O(n^{-2}h^{-3})$, as we shall see in Proposition \ref{p:T1prop} below and of course we will have to assume $h\to 0$.  Since we strongly conjecture that these rates are optimal, the first condition in \eqref{e:asy_equiv} is stronger than requiring that $E$ and $\Ex[T_1]$ are asymptotically equivalent -- there is no hope to achieve this. 

We solve this problem by using a biases-corrected version of $T_1$ following ideas of \cite{ReWa17}. We assume that there is a constant $\gamma>0$ with
\begin{equation}
\lim_{{t\to 0}\atop{t<0}} \frac{1-F(t)}{|t|} = \gamma. \label{e:lim=gamma}
\end{equation}
where $F$ is the distribution function of the $\epsilon_i$. We put
\begin{align*}
T = \sum_{i=1}^{n/2} \hat g(x_{2i-1})^2 + \frac{2}{\gamma'}\sum_{i=1}^{n/2} Y_{2i-1}\mathbf{1}_{\{Y_{2i-1}\ge \hat g(x_{2i-1})\}} &- \frac{2}{n} \Big( \sum_{i=1}^{n/2} \hat g(x_{2i-1})+ \frac{1}{\gamma'}\sum_{i=1}^{n/2} \mathbf{1}_{\{Y_{2i-1}\ge \hat g(x_{2i-1})\}}\Big)^2 \\
&- \frac{ (n/2)\Big(\sum_{i=1}^{n/2} \big(\hat g(x_{2i-1}) + \frac{1}{\gamma'} \mathbf{1}_{\{Y_{2i-1} \ge \hat g(x_{2i-1})\}} \big) \cdot \big(x_{2i-1}-\frac{S}{n/2}\big) \Big)^2 }{ Rn/2-S^2} 
\end{align*}
for some $\gamma'\approx\gamma$. We reduce the bias only sufficiently if $\gamma'$ is close to $\gamma$. However, $\gamma$ is unknown in many applications. Hence we will plug in an estimator $\hat\gamma$ fulfilling certain consistency assumptions for $\gamma'$.

Any estimator $\hat\gamma$ fulfilling the assumptions (G\ref{i:G1})-(G\ref{i:G_alternative}) in Section \ref{s:asymptotic} below will be eligible. However, the only estimator for which we checked these assumptions is 
\begin{equation} \hat \gamma := \frac{ \frac{2k}{n}}{ \hat\epsilon_{\frac{n}{2}:\frac{n}{2}} - \hat\epsilon_{\frac{n}{2}-k:\frac{n}{2}} }, \label{e:negHill} \end{equation}
where $\hat\epsilon_{1:\frac{n}{2}}, \dots, \hat\epsilon_{\frac{n}{2}:\frac{n}{2}}$ are the ascendingly sorted order statistics of the residuals
\[ \hat\epsilon_{2i} := Y_{2i} - \hat g(x_{2i}). \]
It will be necessary to apply the estimator $\hat g$ with a bandwidth $h_1$, which is much larger than $h$, here. 
This estimator is motivated by the negative Hill estimator from \cite{Fa95}.

The reason for the application of the split-of-the-sample technique is that we are not able to treat the variance of $T$ without this assumption, because then $Y_i$ and $\hat g(x_i)$ are no longer independent. Since we do not know what the variance of $T$ would be without the split-of-the-sample technique, we cannot tell whether this theoretical problem would be the only problem or whether there is a practical problem hidden behind it.

Let $z_q$ denote the $q$-quantile of the $\mathcal{N}(0,1)$-distribution. Fix a level $\lambda\in(0,1)$. Now the two tests we are considering are
\[\varphi_1 = \begin{cases} 1 & \mbox{if } T \ge  z_{1-\lambda} \cdot \sqrt{\frac{8}{(C_x)^3}  n^{-2}h^{-3}\hat\gamma^{-4}} \\
													0 & \mbox{else}, 
						\end{cases}
\]
where $C_x$ is the constant defined in equation \eqref{e:ass_scatter} below (in particular $C_x=\frac{1}{2}$ if the observation points are equidistant), and 
\[ \varphi_2 = \begin{cases} 1 & \mbox{if } T \ge z_{1-\lambda} \cdot \sqrt{n^{-2}h^{-3}\hat\gamma^{-4}A_1} \\ 0 & \mbox{else},  \end{cases} \]
where the constant $A_1$ is defined by the following construction.  For $\gamma>0$ and two point measures $\phi^o$ on $[0,1]\times (-\infty,0]$ and $\phi^e$ on $[-1,2]\times (-\infty,0]$ we define 
\[  G(\phi^o, \phi^e, \gamma) := \frac{1}{2}\int_0^1\tilde g(x)^2 \, dx + \frac{2}{\gamma} \sum_{i=1}^\infty \mathfrak{Y}_i  \mathbf{1}_{\{ \mathfrak{Y}_i \ge \tilde g(\mathfrak{X}_i)\}}, \]
where $\phi^o=\{ (\mathfrak{X}_i, \mathfrak{Y}_i) \mid i=1,\dots\}$ and $\tilde g$ is the estimator $\hat g$ at bandwidth $h=1$ applied to the data set $\phi^e$. We let $\Phi_l^o$ and $\Phi_l^e$, $l=1,\dots,5$, be Poisson processes on $[0,1]\times (-\infty,0]$ and $[-1,2]\times(-\infty,0]$ respectively with intensity $\gamma$ and put $A_\gamma:= \sum_{l=1}^5\Cov\big( G(\Phi_3^o,\Phi_3^e, \gamma), G(\Phi_{l}^o, \Phi_{l}^e, \gamma) \big)$.\label{A1def} Since in Subsection \ref{s:A1sim} we see that $A_1$ is approximately $13.7$, we replace it by $13.7$ in applications.

The motivation for the choices of the critical values comes from the asymptotic distribution of the test statistic that will be established in the appendix. It will be shown to be asymptotically Gaussian with mean zero and a variance that is asymptotically bounded by $\frac{8}{(C_x)^3}  n^{-2}h^{-3}\gamma^{-4}$ and that asymptotically equals $n^{-2}h^{-3}\gamma^{-4}A_1$. 

The advantages of the two tests will be discussed in Section \ref{s:asymptotic}.

\section{Asymptotic behavior of the tests}\label{s:asymptotic}

In this section we analyse the asymptotic behavior of $T$.  

The bandwidth $h$ has to be chosen depending on $n$ in such a way that $\lim_{n\to\infty} nh=\infty$, but $\lim_{n\to\infty} h=0$. Furthermore we have to assume $\lim_{n\to\infty} nh^2=\infty$ in order to ensure that the standard deviation of the test statistic $T$ asymptotically dominates its bias and that we hence do not run into the same problems we get when using $T_1$ as test statistic. \medskip

In order to ensure an appropriate asymptotic behavior of $T$, we need an assumption on the distribution $F$ of the errors $\epsilon_i$ which is slightly stronger than \eqref{e:lim=gamma}. We shall assume that 
\begin{equation}  \mbox{there is some constant $C_F$ with }  |F(z)-(1-\gamma |z|)| \le C_F|z|^2 \mbox{ for all $z<0$.} \label{e:ass_distr} \end{equation}
Furthermore, we assume that 
\begin{equation} \mbox{there are $\Gamma>0$ and $\tilde C_F>0$ such that } F(t) \le \tilde C_F \cdot |t|^{-\Gamma}. \label{e:ass_distr_tail} \end{equation}\medskip

We shall assume in the sequel that the data points $x_i, \, i=1,\dots, n,$ are scattered regularly enough such that
\begin{equation}
C_x:= \liminf_{n\to\infty} \frac{\inf\big\{ \#\{i\in\{1,\dots, n\} \mid x_i\in[x,x+h/2)\} \mid x\in(-h,1+h/2) \big\}}{nh} >0 \label{e:ass_scatter}
\end{equation}
and that
\begin{equation}
C'_x:= \limsup_{n\to\infty} \frac{\sup\big\{ \#\{i\in\{1,\dots, n\} \mid x_i\in[x,x+h/2)\} \mid x\in(-h,1+h/2) \big\}}{nh} <\infty. \label{e:ass_scatter_up}
\end{equation}
In order to ensure that $\varphi_2$ has asymptotically the correct size, we will need a much stronger assumption, namely: 
\begin{enumerate}[{(A}1{)}]
\item There is a bandwidth function $H$ with $\lim_{n\to\infty}\frac{H}{h} =0$ such that\label{i:ass_scatter_strong}
\begin{align*} 
\lim_{n\to\infty} &\frac{ \sup\big\{ \#\{i\in\{1,\dots,n\} \mid x_i\in[x,x+H) \} \mid x\in (-h/2,1+h/2-H] \big\}}{nH} \\
&=  \lim_{n\to\infty} \frac{ \inf\big\{ \#\{i\in\{1,\dots,n\} \mid x_i\in[x,x+H) \} \mid x\in (-h/2,1+h/2-H] \big\}}{nH} =1. 
\end{align*}
\end{enumerate}
For the investigation of the behavior of the tests under the alternative, we need a weaker version of \eqref{e:ass_scatter}, namely
\begin{equation}
 \lim_{n\to\infty} \frac{\inf\big\{ \#\{i\in\{1,\dots, n\} \mid x_i\in[x,x+h/2)\} \mid x\in(-h,1+h/2) \big\}}{h} =\infty, \label{e:ass_scatter_alt}
\end{equation}\medskip
and the following assumption:
\begin{enumerate}[{(A}1{)}]
\setcounter{enumi}{1}
\item There is a measure $\mu$ on $[0,1]$ such that $\frac{2}{n} \sum_{i=1}^{n/2} \delta_{x_{2i-1}} \to \mu$ weakly. \label{i:ass_scatter_converge}
\end{enumerate}
This assumption is implied by (A\ref{i:ass_scatter_strong}) and independent of \eqref{e:ass_scatter} and \eqref{e:ass_scatter_up}.\medskip

For the estimator $\hat\gamma$ we have to impose the following assumptions:
\begin{enumerate}[{(G}1{)}]
\item $\hat\gamma$ is independent of $Y_{2i-1}, i=1, \dots, n/2,$ under $H_0$ \label{i:G1}
\item $\Ex[(1/\hat\gamma-1/\gamma)^4] \in O((nh)^{-2})$ under $H_0$\label{i:G4}
\item $\mathbb{P}( \hat \gamma < t) \le \hat C t^{9/2}$ for some constant $\hat C$ independent of $n$ under $H_0$\label{e:gamma_large-dev}
\item $\sup_n \Ex\big[ (\frac{1}{\hat\gamma})^2 \big]<\infty$ both under $H_0$ and under $H_1$\label{i:G_alternative}  
\end{enumerate}
Some of our lemmata will hold under weaker versions of (G\ref{i:G4}), namely one of the following:
\begin{enumerate}[{(G}1{)}]
\setcounter{enumi}{4}
\item $\Ex[(1/\hat\gamma-1/\gamma)^2]\le \hat C/(nh)$ for some constant $\hat C$ \label{i:G2}
\item $\lim_{n\to\infty} \Ex[(1/\hat\gamma-1/\gamma)^4]=0$ \label{i:G3}
\end{enumerate}

If one uses specifically the estimator $\hat\gamma$ from \eqref{e:negHill}, then these assumptions are implied by the following ones, as we shall see in Subsection \ref{s:negHill}.  We assume that
\begin{align}
 \limsup_{n\to \infty} \frac{n^2h^2}{h_1^5k^4} &< \infty, \label{e:ass1}\\
 \limsup_{n\to \infty} n^\alpha\frac{nh}{k} &=0 \mbox{ for some $\alpha>0$,} \label{e:ass2}\\
\lim_{n\to\infty} \frac{k^2h}{n} &=0, \label{e:ass2a} \\
\lim_{n\to\infty} \frac{n}{k} &=\infty, \label{e:ass3} \\
\lim_{n\to\infty} \frac{n}{k^\beta} &= 0 \mbox{ for some $\beta>0$,} \label{e:ass4}\\
\lim_{n\to\infty} \frac{n}{k}\omega_g(h_1) &\to 0 \mbox{ for all $g$ from the alternative} \label{e:ass5}\\
\lim_{n\to\infty} kh_1^2&=\infty \label{e:ass6}\\
\mbox{and }\lim_{n\to\infty} \frac{nh}{n^\gamma} &= \infty \mbox{ for some $\gamma>0$.} \label{e:ass8}
\end{align}
\begin{remark} If all functions $g$ from the alternative are H\"older continuous of order $2/3+\epsilon$ for some $\epsilon\in(0, \frac{5}{24})$, then assumptions \eqref{e:ass1}--\eqref{e:ass8} together with $\lim_{n\to\infty} nh^2=\infty$, which is assumed in most of our results, can be fulfilled by choosing $h=n^{-1/2+\epsilon/5}$, $h_1=n^{-3/8+\epsilon/5}$ and $k=n^{3/4-\epsilon/5}$.
\end{remark}

\begin{remark} There is some redundancy in the assumptions \eqref{e:ass1}--\eqref{e:ass8}. For example, \eqref{e:ass2} together with the assumption $\lim_{n\to\infty} nh^2=\infty$ implies \eqref{e:ass4}. We accept that in order to state each result only under the assumptions used in its proof. \end{remark}

\begin{remark} Observe that the two bandwidths $h$ and $h_1$ have to be chosen differently. Indeed, \eqref{e:ass1} and \eqref{e:ass2a} imply
\[ \lim_{n\to\infty} \frac{h^4}{h_1^5}  =  \lim_{n\to\infty} \frac{n^2h^2}{h_1^5k^4} \cdot \big( \frac{k^2h}{n} \big)^2 =0. \]
 \end{remark}

We are now in the position to state the main results of this paper. 

\begin{theorem}\label{c:test}
Assume that the design points fulfill \eqref{e:ass_scatter} and \eqref{e:ass_scatter_up}, that the error distribution fulfills \eqref{e:ass_distr} and \eqref{e:ass_distr_tail}, that $\lim_{n\to\infty} nh^2=\infty$ and that the estimator $\hat\gamma$ for $\gamma$ fulfills (G\ref{i:G1}) and (G\ref{i:G4}). Then for any $\lambda \in (0,1)$ the test
\[\varphi_1 = \begin{cases} 1 & \mbox{if}\quad T \ge  z_{1-\lambda} \cdot \sqrt{\frac{8}{(C_x)^3}  n^{-2}h^{-3}\hat\gamma^{-4}} \\
													0 & \mbox{else} 
						\end{cases}
\]
has asymptotically at most size $\lambda$. 
\end{theorem}	

\begin{theorem}\label{c:tests}
Let the level $\lambda\in(0,1)$. 
Assume that the design points fulfill (A\ref{i:ass_scatter_strong}), the errors satisfy \eqref{e:ass_distr} and \eqref{e:ass_distr_tail}, that $\lim_{n\to\infty} nh^2=\infty$ and that $\hat\gamma$ is an estimator for $\gamma$ with (G\ref{i:G1}), (G\ref{i:G4}) and (G\ref{e:gamma_large-dev}). Then 
	\[ \varphi_2 = \begin{cases} 1 & \mbox{ if } T \ge z_{1-\lambda} \sqrt{n^{-2}h^{-3}\hat\gamma^{-4}A_1} \\ 0 & \mbox{ if } T<z_{1-\lambda} \sqrt{n^{-2}h^{-3}\hat\gamma^{-4}A_1} \end{cases} \]
	has asymptotically size $\lambda$.
\end{theorem}

The advantage of $\varphi_1$ is that it can cope with irregularly scattered data points. In order to apply this test, we only have to assume \eqref{e:ass_scatter} and \eqref{e:ass_scatter_up}, whereas we need (A\ref{i:ass_scatter_strong}) in order to apply $\varphi_2$. 

On the other hand several estimates in the construction of $\varphi_1$ make it very conservative. In particular, the level of $\varphi_1$ only bounds the limit of the size of $\varphi_1$, while the level of $\varphi_2$ equals the limit of the size of $\varphi_2$. Therefore it is not surprising that $\varphi_1$ is much weaker in detecting alternatives than $\varphi_2$ as demonstrated in the simulation study (Section \ref{s:application}).

\begin{theorem}\label{t:consistent}
Let $g:\mathbb{R}\to\mathbb{R}$ be a uniformly continuous function that is not affine on $[0,1]$. Assume that the errors fulfill \eqref{e:lim=gamma} and \eqref{e:ass_distr_tail}, that the estimator $\hat\gamma$ fulfills (G\ref{i:G_alternative}) and that the design points satisfy \eqref{e:ass_scatter_alt} and (A\ref{i:ass_scatter_converge}). Then 
\[ \lim_{n\to\infty} \mathbb{P}(\varphi=1) =1, \]
where $\varphi$ is either of the two tests $\varphi_1$ or $\varphi_2$.  
\end{theorem}

Notice that while we only have to assume that $g:\mathbb{R}\to\mathbb{R}$ is uniformly continuous in order to ensure that Theorem \ref{t:consistent} is mathematically correct, we do not know any estimator $\hat\gamma$ fulfilling (G\ref{i:G_alternative}) unless $g$ is H\"older continuous of order $\frac{2}{3}+\epsilon$ with known $\epsilon>0$.

\section{Simulations and real data application}\label{s:application}

\subsection{Simulation of $A_1$}\label{s:A1sim}

In order to apply the test $\varphi_2$ it is important to know the value of the constant $A_\gamma$ from page \pageref{A1def} at least in the case $\gamma=1$. There is no hope for an analytic expression for $A_1$ and also evaluating it numerically seems to be very difficult. But a simulation technique is obvious from its definition: Produce a large sample $\Phi_{l,i}^e, \Phi_{l,i}^o$, $l=1,\dots, 5$, $i=1, \dots, n$ of i.i.d.\ copies of $\Phi_{l}^e, \Phi_{l}^o$, $l=1,\dots, 5$ and then approximate $A_1$ by the empirical covariance of $G(\Phi_{3,i}^e, \Phi_{3,i}^o, 1)$, $i=1,\dots, n$ and $\sum_{l=1}^5 G(\Phi_{l,i}^e, \Phi_{l,i}^o, 1)$, $i=1,\dots, n$. We carried out this experiment with $n=10^5$ and got an empirical covariance of  $A_1\approx 13.7$. 

In order to evaluate the accuracy of this estimation, we calculate its standard deviation. The variance of the empirical covariance of two samples $X_i,\, i=1,\dots,n,$ and $Y_i,\, i=1, \dots, n,$ is
\[\frac{1}{n} \Big( \tilde \mu_4- \frac{n-2}{n-1} (\eta^2)^2 + \frac{1}{n-1}\sigma^2\tau^2 \Big),\]
where $\tilde\mu_4:= \Ex[(X_1-\mu)^2(Y_1-\nu)^2]$, $\mu:=\Ex[X_1]$, $\nu:=\Ex[Y_1]$, $\sigma^2:=\Var(X_1)$, $\tau^2:=\Var(Y_1)$ and $\eta^2:=\Cov(X_1,Y_1)$. We estimate this quantity by plugging in the mean values of the same sample that was used for the estimation of $A_1$ itself for the expected values, and the empirical (co-)variances for the (co-)variances. We got a variance of the empirical covariance of $0.06$, i.e.\ a standard deviation of $0.25$. So already the last digit we reported above is not really reliable.

\subsection{Simulation under the null hypotheses}\label{ss:sim_null}

In this subsection we determine the size of the tests $\varphi_2$ for finite sample size based on simulations.

For bandwidths $h=h_1=0.2$ and $n\in\{20,50,100,200\}$ we considered design points $-h, -h+\frac{1}{n}, -h+\frac{2}{n}, \dots$. As true regression function we considered $g\equiv 0$ in this subsection and we used two different error distributions, namely the uniform distribution on $[-1,0]$ and the negative of an exponentially distributed random variable with rate $1$. We applied the test $\varphi_2$, where we consider $\gamma$ either to be known or we estimate it using the estimator from \eqref{e:negHill} with $k\in\{20,50\}$. Notice that this does make sense only if $k<n/2$. As levels we used $\lambda=0.05$ and $\lambda=0.01$. For each setup we performed 1,000 independent simulation runs and estimated the size by the ratio of simulation runs in which the hypotheses was rejected. Motivated by the results from Subsection \ref{s:A1sim} we approximated $A_1$ by $13.7$. The results are reported in Tables  \ref{t:sim_null_reg_u2} and \ref{t:sim_null_reg_e2}.

We see that the test $\varphi_2$ is very conservative. In order to see whether this improves with increasing sample size, we performed simulations with $n=500$ and $n=1,000$. Due to the large computation times we performed this additional simulations only in the case of uniform errors and only in the case that either $\gamma$ is known or $k=20$. Moreover, the results for $n=1,000$, $k=20$ are based on only 200 simulation runs. These results are reported in Table \ref{t:sim_null_reg_u2}, too. Choosing $k=20$ seems to be the magic trick which produces good results already at small sample sizes. 


We also did some experiments with the test $\varphi_1$. Not surprisingly, its size was always smaller than that of $\varphi_2$ (with equality only if the estimated sizes of both tests were $0$). This means that it does not exploit the level at the price of having problems in detecting alternatives. In conclusion, the properties of $\varphi_1$ are poorer than that of $\varphi_2$. We also tried $k\in\{5,10\}$ in the estimation of $\gamma$. However, with this choice the size of $\varphi_2$ exceeded the level by far which is caused by an overestimation of $ \gamma$.

\begin{table}
\begin{center}
\begin{tabular}{c|cccccc}
 & 20 & 50 & 100 & 200&500&1,000\\\hline
$\gamma$ known & 0 & 0 & 0.003 & 0.003 &0.005 &0.010\\
20 & - & 0.035 & 0.050 & 0.046 &0.066 &0.05\\
50 & - & - & - & 0.017
\end{tabular}\hspace{1cm}
\begin{tabular}{c|cccccc}
& 20 & 50 & 100 & 200 & 500& 1,000\\\hline
$\gamma$ known & 0 & 0 & 0.001 & 0.001&0.002&0.002\\
20 & - & 0.011 & 0.022 & 0.020&0.032&0.02\\
50 & - & - & - & 0.003
\end{tabular}
\end{center}

Estimated sizes of the test $\varphi_2$ for level 5\% (top) and 1\% (bottom) applied to data that is uniformly distributed on
$[-1,0]$. In the first row we present the result for known $\gamma$, while in the further rows we present results for estimated
values of $\gamma$ with different choices of $k$. In different columns we present results for different values of $n$.
\caption{Size of the test $\varphi_2$: Uniform errors}
\label{t:sim_null_reg_u2}
\end{table}

\begin{table}
\begin{center}
\begin{tabular}{c|cccc}
 & 20 & 50 & 100 & 200\\\hline
 $\gamma$ known & 0.044 & 0.025 & 0.019 & 0.015\\
20 & - & 0 & 0.023 & 0.032\\
50 & - & - & - & 0.001
\end{tabular}\hspace{1cm}
\begin{tabular}{c|cccc}
 & 20 & 50 & 100 & 200\\\hline
$\gamma$ known & 0.017 & 0.009 & 0.008 & 0.004\\
20 & - & 0 & 0.007 & 0.017\\
50 & - & - & - & 0
\end{tabular}
\end{center}

The same as Table \ref{t:sim_null_reg_u2} except that the errors are now negative exponential.
\caption{Size of the test $\varphi_2$: Negative exponential errors}
\label{t:sim_null_reg_e2}
\end{table}



\subsection{The behavior under the alternative}

In this subsection we investigate the behavior of our tests under the alternative. We consider the true regression functions
\begin{align*}
g_1(x)&= c\cdot \sin(\alpha\pi x) \\
g_2(x)&= c\cdot (x-x_0)^p \\
g_3(x)&= -c\cdot (x-x_0)^p.
\end{align*}

We carried out experiments with $c\in\{0.1, 0.2, 0.5, 1\}$ and $\alpha\in \{1,2,4\}$ for $g_1$ and with $p\in\{2,6,7\}$, $c\in\{0.1, 0.2, 0.5, 1, 2, 5\}$ and $x_0\in\{0,0.25,0.5\}$ for $g_2$ and $g_3$. We chose the values of $p$ this way since we wanted to study the difference between small and large  values of $p$ as well as the difference between even and odd values. We restrict ourselves to the case that $n=100$, that  $k=20$ or $\gamma$ is known and that the errors are uniformly distributed.

The results for $g_1$ are shown in Tables \ref{T:alt_sin_100_known} and \ref{T:alt_sin_100_20}. Of course, the power increases with increasing value of $c$. The dependency on $\alpha$ varies. In many cases we have the highest power for $\alpha=2$. This is explained by the following fact: The function $\sin(\pi x)$ can be well approximated by the identity and thus has a much smaller ``distance'' (measured in what metric ever) to the space of affine functions than $\sin(2\pi x)$. From $\sin(2\pi x)$ to $\sin(4\pi x)$ we do not have a decrease in the distance to the space of affine functions, but only an increase in the frequency, which makes the detection more difficult. 

In Tables \ref{T:alt_x2_100_known} and \ref{T:alt_x2_100_20} we reported the results for $g_2$ with $p=2$. We see that the powers are much smaller than the powers for $g_1$ which is explained by the fact that for the same value of $c$ the distance of $g_2$ to the space of affine functions is much smaller than the distance of $g_1$. Moreover, we see that the results for $g_2$ with $p=2$ hardly depend on $x_0$. This had to be expected -- changing the value of $x_0$ just corresponds to adding an affine function and the tests we proposed are invariant under addition of affine functions. So the differences for different values of $x_0$ come only from the use of new random numbers. In Tables \ref{T:alt_x6_100_20} and \ref{T:alt_x7_100_20} we see that for $p=6$ and $p=7$ the choice of $x_0$ matters: The closer $x_0$ is to the centre of the interval $[0,1]$, the smaller is the power of the test. This is not surprising, since the distance between the regression function $g$ and the space of affine function is small when $x_0$ is close to the centre of the interval. For $x_0=0$ the power of the test increases with increasing $p$, while for $x_0=0.25$ and $x_0=0.5$ the power decreases with increasing $p$. The results for $p=6$ and $p=7$ are hardly different, so it does not seem to be of importance whether $p$ is even or odd. 

Finally we present the results for $g_3$ in Tables \ref{T:alt_mx2_100_known} -- \ref{T:alt_mx7_100_20}. The results for $g_3$ are similar to those of $g_2$ -- in some situations (in particular, for $p=2$) they are better, in some other situations (in particular for larger $p$) they are worse. 

The behavior of the test $\varphi_2$ under the alternative at sample size $n=100$ is of acceptable quality. It detects alternatives far from the zero hypotheses with high probabilities but it has problems in detecting alternatives at medium distance from the null hypotheses.      

Finally we investigated how the power increases with increasing sample size $n$. We did this only at level $5\%$ with known $\gamma$ and the $\sin$-function as regression function. The results are reported in Table \ref{T:alt_sin_n}. We see an increase of the power function with increasing $n$ that is approximately as strong as we expected it.

\begin{table}
\begin{center}
\begin{tabular}{c|ccccc}
 & 0($H_0$) & 0.1 & 0.2 & 0.5 & 1\\\hline
1 &  & 0.003 & 0.043 & 0.729 & 0.991\\
2 & 0.003 & 0.008 & 0.143 & 0.950 & 1.000\\
4 &  & 0.021 & 0.262 & 0.952 & 1.000\\
\end{tabular}
\hspace{1cm}
\begin{tabular}{c|ccccc}
 & 0($H_0$) & 0.1 & 0.2 & 0.5 & 1\\\hline
1 &  & 0.001 & 0.006 & 0.594 & 0.986\\
2 & 0.001 & 0.002 & 0.038 & 0.913 & 1.000\\
4 &  & 0.003 & 0.099 & 0.874 & 0.999
\end{tabular}
\end{center}

Estimated power function for $g(x)=c \cdot \sin(\alpha \pi x)$ with $n=100$ and known $\gamma$. The results for different values of $c$ are reported in different columns, while the results for different values of $\alpha$ are reported in different rows.  The left table is for level $5\%$ and the right table is for level $1\%$. 
\caption{Power under the alternative: The $\sin$-function with $n=100$ and known $\gamma$}
\label{T:alt_sin_100_known}
\end{table}



\begin{table}
\begin{center}
\begin{tabular}{c|ccccc}
& 0($H_0$) & 0.1 & 0.2 & 0.5 & 1\\\hline
1 &  & 0.087 & 0.265 & 0.893 & 0.998\\
2 & 0.050 & 0.111 & 0.423 & 0.969 & 0.999\\
4 &  & 0.140 & 0.431 & 0.867 & 0.991
\end{tabular}
\hspace{1cm}
\begin{tabular}{c|ccccc}
& 0($H_0$) & 0.1 & 0.2 & 0.5 & 1\\\hline
1 &  & 0.046 & 0.151 & 0.854 & 0.996\\
2 & 0.022 & 0.054 & 0.264 & 0.953 & 0.999\\
4 &  & 0.062 & 0.235 & 0.696 & 0.976
\end{tabular}
\end{center}

Estimated power function for $g(x)=c \cdot \sin(\alpha \pi x)$ with $n= 100 $, $k= 20 $. The further details are the same as in Table \ref{T:alt_sin_100_known}.
\caption{Power under the alternative: The $\sin$-function with $n=100$, $k=20$}
\label{T:alt_sin_100_20}
\end{table}

\begin{table}
\begin{center}
\begin{tabular}{c|ccccccc}
& 0($H_0$) & 0.1 & 0.2 & 0.5 & 1 & 2 & 5\\\hline
0 && 0.002 & 0.002 & 0.003 & 0.073 & 0.785 & 1.000\\
0.25 & 0.003 & 0.005 & 0.002 & 0.014 & 0.056 & 0.781 & 1.000\\
0.5 && 0.003 & 0.003 & 0.008 & 0.061 & 0.772 & 1.000
\end{tabular}
\hspace{1cm}
\begin{tabular}{c|ccccccc}
& 0($H_0$) & 0.1 & 0.2 & 0.5 & 1 & 2 & 5\\\hline
0 && 0 & 0.001 & 0 & 0.019 & 0.576 & 1.000\\
0.25 & 0.001 & 0.001 & 0 & 0.004 & 0.014 & 0.578 & 1.000\\
0.5 && 0 & 0 & 0.001 & 0.014 & 0.575 & 1.000
\end{tabular}
\end{center}

Estimated power function for $g(x)=c \cdot (x-x_0)^2$ with $n=100$ and known $\gamma$. In different rows the results for different values of $x_0$ are presented. The further details are the same as in Table \ref{T:alt_sin_100_known}.
\caption{Power under the alternative: The $x^2$-function with $n=100$ and known $\gamma$}
\label{T:alt_x2_100_known}
\end{table}



\begin{table}
\begin{center}
\begin{tabular}{c|ccccccc}
 & 0($H_0$) & 0.1 & 0.2 & 0.5 & 1 & 2 & 5\\\hline
0 && 0.044 & 0.047 & 0.071 & 0.271 & 0.868 & 1.000\\
0.25& 0.050 & 0.046 & 0.055 & 0.089 & 0.252 & 0.866 & 1.000\\
0.5 && 0.039 & 0.036 & 0.079 & 0.259 & 0.857 & 1.000
\end{tabular}
\hspace{1cm}
\begin{tabular}{c|ccccccc}
& 0($H_0$) & 0.1 & 0.2 & 0.5 & 1 & 2 & 5\\\hline
0 && 0.013 & 0.019 & 0.030 & 0.142 & 0.754 & 1.000\\
0.25 & 0.022 & 0.022 & 0.023 & 0.043 & 0.115 & 0.756 & 1.000\\
0.5 && 0.022 & 0.017 & 0.037 & 0.136 & 0.748 & 1.000
\end{tabular}
\end{center}

Estimated power function for $g(x)=c \cdot (x-x_0)^2$ with $n= 100 $, $k= 20 $. The further details are the same as in Table \ref{T:alt_x2_100_known}.
\caption{Power under the alternative: The $x^2$-function with $n=100$, $k=20$}
\label{T:alt_x2_100_20}
\end{table}

\begin{table}
\begin{center}
\begin{tabular}{c|ccccccc}
 & 0($H_0$) & 0.1 & 0.2 & 0.5 & 1& 2 &5\\\hline
0 && 0.049 & 0.052 & 0.286 & 0.908 & 1.000 & 1.000 \\
0.25 &0.050 & 0.043 & 0.046 & 0.052  & 0.051 & 0.168 & 0.954 \\
0.50 && 0.039 & 0.034 & 0.050 & 0.048 & 0.043 & 0.046
\end{tabular}
\hspace{1cm}
\begin{tabular}{c|ccccccc}
& 0($H_0$) & 0.1 & 0.2 & 0.5 & 1 & 2 & 5\\\hline
0 && 0.018 & 0.022 & 0.140 & 0.837 &0.998 &1.000 \\
0.25 & 0.022 & 0.022 & 0.018 & 0.031 & 0.022 & 0.065 & 0.922 \\
0.50 && 0.021 & 0.017 & 0.018 & 0.019 & 0.021 & 0.017
\end{tabular}
\end{center}

Estimated power function for $g(x)=c \cdot (x-x_0)^6$ with $n= 100 $, $k= 20 $. The further details are the same as in Table \ref{T:alt_x2_100_known}.
\caption{Power under the alternative: The $x^6$-function with $n=100$, $k=20$}
\label{T:alt_x6_100_20}
\end{table}

\begin{table}
\begin{center}
\begin{tabular}{c|ccccccc}
 & 0($H_0$) & 0.1 & 0.2 & 0.5 & 1& 2 &5\\\hline
0 && 0.047 & 0.053 & 0.361 & 0.940 & 1.000 & 1.000\\
0.25 & 0.050 & 0.043 & 0.046 & 0.050 & 0.045 & 0.102 & 0.910\\
0.5 && 0.039 & 0.036 & 0.048 & 0.051 & 0.049 & 0.048
\end{tabular}
\hspace{1cm}
\begin{tabular}{c|ccccccc}
& 0($H_0$) & 0.1 & 0.2 & 0.5 & 1 & 2 & 5\\\hline
0&& 0.017 & 0.024 & 0.193 & 0.905 & 1.000 & 1.000 \\
0.25 & 0.022& 0.023 & 0.019 & 0.031 & 0.019 & 0.037 & 0.834\\
0.5 && 0.021 & 0.017 & 0.020 & 0.021 & 0.023 & 0.019
\end{tabular}
\end{center}

Estimated power function for $g(x)=c \cdot (x-x_0)^7$ with $n= 100 $, $k= 20 $. The further details are the same as in Table \ref{T:alt_x2_100_known}.
\caption{Power under the alternative: The $x^7$-function with $n=100$, $k=20$}
\label{T:alt_x7_100_20}
\end{table}

\begin{table}
\begin{center}
\begin{tabular}{c|ccccccc}
& 0($H_0$) & 0.1 & 0.2 & 0.5 & 1 & 2 & 5\\\hline
0 &  & 0.003 & 0.002 & 0.011 & 0.099 & 0.652 & 0.998\\
0.25 & 0.003 & 0.004 & 0.004 & 0.011 & 0.091 & 0.677 & 0.998\\
0.5 &  & 0.002 & 0.003 & 0.006 & 0.084 & 0.663 & 0.997
\end{tabular}
\hspace{1cm}
\begin{tabular}{c|ccccccc}
& 0($H_0$) & 0.1 & 0.2 & 0.5 & 1 & 2& 5\\\hline
0 &  & 0 & 0.001 & 0.002 & 0.020 & 0.525 & 0.998\\
0.25 & 0.001 & 0.001 & 0 & 0.002 & 0.024 & 0.545 & 0.998\\
0.5 &  & 0 & 0.001 & 0.001 & 0.020 & 0.532 & 0.995
\end{tabular}
\end{center}

Estimated power function for $g(x)=-c \cdot (x-x_0)^2$ with $n=100$ and known $\gamma$. In different rows the results for different values of $x_0$ are presented, while in different columns the results for different values of $c$ are presented. The further details are the same as in Table \ref{T:alt_sin_100_known}.
\caption{Power under the alternative: The negative $x^2$-function with $n=100$ and known $\gamma$}
\label{T:alt_mx2_100_known}
\end{table}



\begin{table}
\begin{center}
\begin{tabular}{c|ccccccc}
& 0($H_0$) & 0.1 & 0.2 & 0.5 & 1 & 2 & 5\\\hline
0 &  & 0.049 & 0.052 & 0.115 & 0.400 & 0.895 & 1.000\\
0.25 & 0.050 & 0.055 & 0.055 & 0.107 & 0.405 & 0.881 & 1.000\\
0.5 &  & 0.045 & 0.050 & 0.110 & 0.403 & 0.871 & 1.000\\
\end{tabular}
\hspace{1cm}
\begin{tabular}{c|ccccccc}
& 0($H_0$) & 0.1 & 0.2 & 0.5 & 1 & 2 & 5\\\hline
0 &  & 0.018 & 0.027 & 0.055 & 0.252 & 0.842 & 0.999\\
0.25 & 0.022 & 0.026 & 0.023 & 0.050 & 0.225 & 0.831 & 0.999\\
0.5 &  & 0.023 & 0.023 & 0.045 & 0.239 & 0.823 & 1.000
\end{tabular}
\end{center}

Estimated power function for $f(x)=-c \cdot (x-x_0)^2$ with $n=100$, $k=20$. The further details are the same as in Table \ref{T:alt_mx2_100_known}.
\caption{Power under the alternative: The negative $x^2$-function with $n=100$, $k=20$}
\label{T:alt_mx2_100_20}
\end{table}

\begin{table}
\begin{center}
\begin{tabular}{c|ccccccc}
 & 0($H_0$) & 0.1 & 0.2 & 0.5 & 1&2&5\\\hline
0 &  & 0.003 & 0.010 & 0.104 & 0.599 & 0.958 & 1.000\\
0.25 &0.003  & 0.004 & 0.003 & 0.009 & 0.009 & 0.041 & 0.528\\
0.5 &  & 0.002 & 0.001 & 0.003 & 0.003 & 0.003 & 0.009
\end{tabular}
\begin{tabular}{c|ccccccc}
 & 0($H_0$) & 0.1 & 0.2 & 0.5 & 1 & 2 & 5\\\hline
0 &  & 0 & 0.003 & 0.025 & 0.463 & 0.949 & 1.000\\
0.25 & 0.001 & 0.001 & 0 & 0.002 & 0.002 & 0.006 & 0.369\\
0.5 &  & 0 & 0.001 & 0.001 & 0 & 0 & 0.001
\end{tabular}
\end{center}

Estimated power function for $g(x)=-c \cdot (x-x_0)^6$ with $n=100$ and known $\gamma$. The further details are the same as in Table \ref{T:alt_mx2_100_known}.
\caption{Power under the alternative: The negative $x^6$-function with $n=100$ and known $\gamma$}
\label{T:alt_mx6_100_known}
\end{table}

\begin{table}
\begin{center}
\begin{tabular}{c|ccccccc}
 & 0($H_0$) & 0.1 & 0.2 & 0.5 & 1 & 2 & 5\\\hline
0 &  & 0.062 & 0.094 & 0.420 & 0.847 & 0.998 & 1.000\\
0.25 & 0.050  & 0.051 & 0.053 & 0.058 & 0.084 & 0.278 & 0.805\\
0.5 &  & 0.040 & 0.037 & 0.048 & 0.054 & 0.070 & 0.109
\end{tabular}
\begin{tabular}{c|ccccccc}
 & 0($H_0$) & 0.1 & 0.2 & 0.5 & 1 & 2 & 5\\\hline
0 &  & 0.020 & 0.038 & 0.274 & 0.796 & 0.996 & 1.000\\
0.25 & 0.022  & 0.026 & 0.022 & 0.029 & 0.051 & 0.154 & 0.735\\
0.5 &  & 0.021 & 0.017 & 0.021 & 0.026 & 0.030 & 0.052
\end{tabular}
\end{center}

Estimated power function for $f(x)=-c \cdot (x-x_0)^6$ with $n=100$, $k=20$. The further details are the same as in Table \ref{T:alt_mx2_100_known}.
\label{T:alt_mx6_100_20}
\caption{Power under the alternative: The negative $x^6$-function with $n=100$, $k=20$} 
\end{table}

\begin{table}
\begin{center}
\begin{tabular}{c|ccccccc}
 & 0($H_0$) & 0.1 & 0.2 & 0.5 & 1 & 2& 5\\\hline
0 &  & 0.003 & 0.010 & 0.112 & 0.592 & 0.953 & 1.000\\
0.25 & 0.003 & 0.004 & 0.002 & 0.008 & 0.008 & 0.018 & 0.283\\
0.5 &  & 0.002 & 0.001 & 0.003 & 0.002 & 0.003 & 0.005
\end{tabular}
\begin{tabular}{c|ccccccc}
 & 0($H_0$) & 0.1 & 0.2 & 0.5 & 1 & 2 & 5\\\hline
0 &  & 0 & 0.002 & 0.028 & 0.473 & 0.943 & 1.000\\
0.25 & 0.001 & 0.001 & 0 & 0.001 & 0.001 & 0.001 & 0.117\\
0.5 &  & 0 & 0.001 & 0.001 & 0 & 0 & 0.001
\end{tabular}
\end{center}

Estimated power function for $g(x)=-c \cdot (x-x_0)^7$ with $n=100$ and known $\gamma$. The further details are the same as in Table \ref{T:alt_mx2_100_known}.
\caption{Power under the alternative: The negative $x^7$-function with $n=100$ and known $\gamma$}
\label{T:alt_mx7_100_known}
\end{table}

\begin{table}
\begin{center}
\begin{tabular}{c|ccccccc}
 & 0($H_0$) & 0.1 & 0.2 & 0.5 & 1 & 2 & 5\\\hline
0 &  & 0.063 & 0.100 & 0.449 & 0.853 & 0.997 & 1.000\\
0.25 & 0.050 & 0.050 & 0.051 & 0.057 & 0.076 & 0.184 & 0.633\\
0.5 &  & 0.040 & 0.035 & 0.049 & 0.049 & 0.055 & 0.054
\end{tabular}
\begin{tabular}{c|ccccccc}
 & 0($H_0$) & 0.1 & 0.2 & 0.5 & 1 & 2 & 5\\\hline
0 &  & 0.020 & 0.043 & 0.291 & 0.801 & 0.995 & 1.000\\
0.25 & 0.022 & 0.026 & 0.021 & 0.029 & 0.041 & 0.086 & 0.515\\
0.5 &  & 0.021 & 0.018 & 0.020 & 0.021 & 0.024 & 0.021
\end{tabular}
\end{center}

Estimated power function for $f(x)=-c \cdot (x-x_0)^7$ with $n=100$, $k=20$. The further details are the same as in Table \ref{T:alt_mx2_100_known}.
\caption{Power under the alternative: The negative $x^7$-function with $n=100$, $k=20$}
\label{T:alt_mx7_100_20}
 \end{table}

\begin{table}
\begin{center}
\begin{tabular}{c|ccccc}
& 0($H_0$) & 0.1 & 0.2 & 0.5 & 1\\\hline
1 &  & 0 & 0 & 0 & 0.003\\
2 & 0 & 0 & 0 & 0 & 0.118\\
4 &  & 0 & 0 & 0.001 & 0.122
\end{tabular}
\hspace{1cm}
\begin{tabular}{c|ccccc}
& 0($H_0$) & 0.1 & 0.2 & 0.5 & 1\\\hline
1 &  & 0 & 0.002 & 0.068 & 0.700\\
2 & 0 & 0.002 & 0 & 0.352 & 0.962\\
4 &  & 0.003 & 0.007 & 0.395 & 0.938
\end{tabular}
\end{center}

\begin{center}
\begin{tabular}{c|ccccc}
& 0($H_0$) & 0.1 & 0.2 & 0.5 & 1\\\hline
1 &  & 0.003 & 0.043 & 0.729 & 0.991\\
2 & 0.003 & 0.008 & 0.143 & 0.950 & 1.000\\
4 &  & 0.021 & 0.262 & 0.952 & 1.000
\end{tabular}
\hspace{1cm}
\begin{tabular}{c|ccccc}
& 0($H_0$) & 0.1 & 0.2 & 0.5 & 1\\\hline
1 &  & 0.062 & 0.504 & 0.986 & 1.000\\
2 & 0.003 & 0.163 & 0.838 & 1.000 & 1.000\\
4 &  & 0.262 & 0.849 & 1.000 & 1.000
\end{tabular}
\end{center}

Estimated power function for  the test $\varphi_2$ with level $5\%$ applied to $g(x)=c \cdot \sin(\alpha \pi x)$ with known $\gamma$ and $n=20$ (top left), $n=50$ (top right), $n=100$ (bottom left) and $n=200$ (bottom right). The results for different values of $c$ are reported in different columns, while the results for different values of $\alpha$ are reported in different rows.
\caption{Power under the alternative: The dependency on $n$ for the $\sin$-function with known $\gamma$}
\label{T:alt_sin_n}
\end{table}

\subsection{Life expectancy data}

Here we want to apply the tests $\varphi_1$ and $\varphi_2$ to life expectancy data.

We consider the following model for the life expectancy: We imagine that in each year $x$ there is an unobserved bound for the life expectancy of a country (far below the bound for the lifetimes of individuals), which is the regression function $g(x)$. Due to random effects the countries do not reach that bound. For each year the observation is the maximal life expectancy observed in any country and the difference between the maximal possible life expectancy and the maximal observed life expectancy is the error. We consider effects like climatic influences or different degrees of health conciseness in different countries not as random effects, but as part of the distribution. So we have only a very small dependency in the observed life expectancy in different years coming from the fact that if in one year many people have died, then in the next year the observation population is smaller.  We consider this dependency to be neglectable. We test the null hypothesis that $g$ is an affine function against the alternative that $g$ is not affine.

At first we apply the two tests we developed to the same data which were considered in \cite{OV02}. This data consists of the maximal life expectancy in each year from 1840 till 2000, where until 1947 the maximal life expectancy is missing in many years. For both tests we chose bandwidths $h=h_1=0.2$. Since we need some observations to the left and to the right of the interval $[0,1]$ in order to apply our tests, we rescaled the time axis in such a way that $1840$ became $-0.2$ and $2000$ became $1.2$. We obtained $n=46$ data points in the interval $[0,1]$. We estimated $\gamma$ by the estimator from  \eqref{e:negHill} with $k=10$. Motivated by the results from Section \ref{s:A1sim} we set $A_1=13.7$. We approximated the constant $C_x$ that is needed for the test $\varphi_2$ by
\[ C_x = \frac{\inf\{ \#\{i\in\{1,\dots, n\}\mid x_i\in [x,x+h/2)\} \mid x\in (-h, 1+h/2)\} }{nh} =\frac{1}{9.2}. \]
We obtained a value of the test statistic $T$ of $4.8$, a critical value of the test $\varphi_1$ on level $5\%$ of $71$ and a critical value of the test $\varphi_2$ of $3.3$. The resulting p-values are $0.46$ for the test $\varphi_1$ and $0.0082$ for the test $\varphi_2$. However, the assumption (A\ref{i:ass_scatter_strong}) needed for the application of the test $\varphi_2$ is not justified, since in many of the early decades $8$ of the $10$ data points are missing, while in the last four decades there is no missing data point at all. So we are not able to judge whether the hypotheses is true for this data set -- the test $\varphi_1$ does not reject the null hypotheses, but it turned out to be to conservative in the simulation study, and the assumptions of the test $\varphi_2$ are not sufficiently justified. 
  
Next we only consider data from the post-war period (from 1948) so that we have absolutely regularly scattered observation points. We added the life expectancy data for the years from 2001 till 2016 from the Human Mortality Database \cite{HMD}. Again, we rescaled the time axis in such a way that the first observation point (now 1948) became $-0.2$ and the last observation point (now 2016) became $1.2$. This time putting $C_x=0.5$ in the test $\varphi_1$ seems to be appropriate. Now we obtain a value of the test statistic $T$ of $1.67$, a critical value of $0.97$ for the test $\varphi_1$ and a critical value of $0.45$ for the test $\varphi_2$. The resulting p-values are $0.0024$ for the test $\varphi_1$ and $5\cdot 10^{-10}$ for the test $\varphi_2$. When changing the bandwidth of the test statistic to $h=0.1$, the value of the test statistic changes to $2.3$, the critical value of the test $\varphi_1$ becomes $0.078$ and the critical value of the test $\varphi_2$ becomes $0.036$. The p-value drops down below $10^{-14}$ for both tests. So we are able to reject the null hypotheses with highest significance. Of course, strictly speaking one has to adjust for multiple testing when applying four tests that test the same hypotheses, but if two of the p-values are below $10^{-14}$, this can be neglected.

\section{Discussion and outlook}\label{s:discuss}

We have proposed two tests for testing whether the regression function in a boundary regression model is an affine function. Using these tests we have been able to show that the boundary of the life expectancy is not an affine function.

It would be good to have a test uniting the benefits of both tests we proposed. A tempting idea is to use inhomogeneous Poisson processes in order to get rid of the strong assumption (A\ref{i:ass_scatter_strong}). However, there is a problem. Without this assumption the processes $\Phi_n^o$ and $\Phi_n^e$ will not converge to inhomogeneous Poisson processes, but to homogeneous Poisson processes whose intensity is determined by the intensity of $\Phi_n^e$ and $\Phi_n^o$ near $0$.  

Of course, it would be desirable to have a goodness-of-fit test for the hypotheses that the true regression function lies in a more general given linear subspace. This is challenging, because then the distribution of the test statistic $T$ will depend on which of the regression functions of this subspace is the true one. 

Moreover, an interesting project is to relax the assumption that $\lim_{t\to 0} (1-F(t))/|t|$ exists to $\lim_{t\to 0} (1-F(t))/|t|^\alpha$ existing for some $\alpha>0$. The difficulty is that for $\alpha\ne 1$ we do not have a decomposition of the test statistic $T$ as $S_1+S_2+S_3$ in which the dependence on $g$ cancels. 
A way out in general might be to change the point of view: We are testing a hypotheses about the location parameter of a Weibull distribution, when the observations are not Weibull distributed but only lie in the domain of attraction of a Weibull distribution. 

In the future one should check whether the proposed test is minimax optimal and, if it is not, propose a minimax optimal test.

Allowing for correlated data would require the estimation of the correlation structure. A further interesting project would be to consider spatial design points.

\begin{appendix}
\section{Proofs}

In this appendix we give the proofs of the theoretical statements of this paper. In Subsection \ref{ss:decomposition} we will prove Proposition \ref{p:T1_explicit} and we will derive an equivalent representation of $T$ that holds at any sample size. This decomposition will be important in the examination of the behavior of the asymptotic behavior of $T$ under the null hypotheses that is carried out in Subsections \ref{ss:moments}, \ref{s:CLT} and \ref{ss:Poisson}. In Subsection \ref{ss:proofs} we give the proofs of Theorem \ref{c:test} and Theorem \ref{c:tests}. The behavior of $T$ under the alternative will be examined in Subsection \ref{ss:consistency}. Finally, we establish the properties of the estimator $\hat \gamma$ from \eqref{e:negHill} in Subsection \ref{s:negHill}.


\subsection{Representations at finite sample size}\label{ss:decomposition}

We abbreviate
\[ K:=\sum_{i=1}^{n/2} \big(\hat g(x_{2i-1}) - f(x_{2i-1})\big) \quad\mbox{ and } \quad L:=\sum_{i=1}^{n/2} \big(\hat g(x_{2i-1})-f(x_{2i-1})\big) \cdot x_{2i-1}. 
\] 

\prf[ of Proposition \ref{p:T1_explicit}]
(i) We have
\begin{align*} T_1 = \min_{\tilde m,\tilde c} \sum_{i=1}^{n/2} \Big( \hat g(x_{2i-1})-\tilde mx_{2i-1}-\tilde c\Big)^2 &= \min_{ m, c} \sum_{i=1}^{n/2} \Big( \hat g(x_{2i-1})-f(x_{2i-1}) - mx_{2i-1}- c\Big)^2  \\
&=\min_{ m, c} \sum_{i=1}^{n/2} \big(\hat g(x_{2i-1})-f(x_{2i-1})\big)^2 + R m^2 + \tfrac{n}{2} c^2 - 2m L - 2 cK + 2Smc.\end{align*}
The minimum is obtained for the solution of
\begin{align*}
2R m - 2 L +  2Sc &=0 \\
 nc  - 2K +2Sm &=0,
\end{align*}
i.e.\ for $c=\frac{RK-SL}{Rn/2-S^2}$ and $m=\frac{nL/2-SK}{Rn/2-S^2}$.  
Hence we get
\begin{align*}
T_1 &= \begin{aligned}[t] \sum_{i=1}^{n/2} \big(\hat g(x_{2i-1})-f(x_{2i-1})\big)^2 + R \Big(\frac{nL/2-SK}{Rn/2-S^2}\Big)^2 + &\frac{n}{2} \Big(\frac{RK-SL}{Rn/2-S^2} \Big)^2 -2\frac{nL/2-SK}{Rn/2-S^2}L \\ &-2\frac{RK-SL}{Rn/2-S^2}K + 2S\Big(\frac{nL/2-SK}{Rn/2-S^2}\Big)\Big(\frac{RK-SL}{Rn/2-S^2}\Big) \end{aligned}
\\
&= \begin{aligned}[t]\sum_{i=1}^{n/2} &\big(\hat g(x_{2i-1})-f(x_{2i-1})\big)^2 + \frac{ R(n/2)^2L^2-RSnKL+ RS^2K^2 + R^2nK^2/2 - RSnKL + S^2nL^2/2 }{\big({Rn/2-S^2}\big)^2} \\
&+ \frac{-2R(n/2)^2L^2 + nS^2L^2 + SRnKL - 2S^3KL}{\big({Rn/2-S^2}\big)^2} \\
&+ \frac{ -R^2nK^2 + 2RS^2K^2 + RSnKL - 2S^3KL}{\big({Rn/2-S^2}\big)^2} \\
&+ \frac{ RSnKL - S^2nL^2 - 2RS^2K^2 + 2S^3KL}{\big({Rn/2-S^2}\big)^2} \end{aligned}\\
&= \begin{aligned}[t] \sum_{i=1}^{n/2} \big(\hat g(x_{2i-1})-f(x_{2i-1})\big)^2 & + \frac{RS^2+ R^2n/2 - R^2n + 2RS^2 - 2RS^2}{\big({Rn/2-S^2}\big)^2} K^2\\
& + \frac{ -RSn - RSn + RSn - 2S^3 + RSn -2S^3 + RSn + 2S^3}{\big({Rn/2-S^2}\big)^2}KL \\
&+\frac{R(n/2)^2 + S^2n/2 - 2R(n/2)^2 + S^2n - S^2n}{\big({Rn/2-S^2}\big)^2}L^2\end{aligned}
\\
&= \sum_{i=1}^{n/2} \big(\hat g(x_{2i-1})-f(x_{2i-1})\big)^2  + \frac{- R^2n/2+RS^2}{\big({Rn/2-S^2}\big)^2} K^2 + \frac{-2S^3+RSn}{\big({Rn/2-S^2}\big)^2}KL +\frac{S^2n/2-R(n/2)^2}{\big({Rn/2-S^2}\big)^2}L^2\\
&= \sum_{i=1}^{n/2} \big(\hat g(x_{2i-1})-f(x_{2i-1})\big)^2  + \frac{- R}{Rn/2-S^2} K^2 + 2\frac{\frac{S}{\sqrt{n/2}}\cdot \sqrt{n/2}}{Rn/2-S^2}KL -\frac{\sqrt{n/2}^2}{Rn/2-S^2}L^2\\
&= \sum_{i=1}^{n/2} \big(\hat g(x_{2i-1})-f(x_{2i-1})\big)^2  + \frac{(-R+\frac{S^2}{n/2})K^2 -\big(\sqrt{n/2}L-\frac{S}{\sqrt{n/2}}K\big)^2}{Rn/2-S^2}\\
&= \begin{aligned}[t]
\sum_{i=1}^{n/2} \big(\hat g(x_{2i-1})-f(x_{2i-1})\big)^2  &- \frac{2}{n} \cdot \Big( \sum_{i=1}^{n/2} \big(\hat g(x_{2i-1}) - f(x_{2i-1}) \big) \Big)^2\\
& - \frac{(n/2) \cdot \Big(\sum_{i=1}^{n/2} \big(\hat g(x_{2i-1}) -f(x_{2i-1}) \big) \cdot \big(x_{2i-1}-\frac{S}{n/2}\big)  \Big)^2}{Rn/2-S^2}.\end{aligned} 
\end{align*}
(ii) Since
\[ R := \sum_{i=1}^{n/2} x_{2i-1}^2= \frac{1}{n^2} (\tfrac{1}{6} n^3 - \tfrac{1}{6} n)= \frac{n^2-1}{6n} \mbox{ and }
S := \sum_{i=1}^{n/2} x_{2i-1} = \frac{1}{n}	(\tfrac{1}{4}n^2) = \frac{1}{4}n,  \]
we get  
\[\frac{S}{n/2}=\frac{1}{2} \quad \mbox{and}\quad \frac{n/2}{Rn/2-S^2} = \frac{n/2}{\tfrac{1}{12}n^2-\tfrac{1}{12} - \tfrac{1}{16}n^2}= \frac{n}{\tfrac{1}{24}n^2 - \tfrac{1}{6} }. \]
Hence the assertion follows. \qed\medskip

In order to analyse the statistic $T$, we consider the alternative form
\begin{align}
 T_f = \sum_{i=1}^{n/2} (\hat g(x_{2i-1})&-f(x_{2i-1}))^2 + \frac{2}{\gamma'}\sum_{i=1}^{n/2} (Y_{2i-1}-f(x_{2i-1}))\mathbf{1}_{\{Y_{2i-1}\ge\hat g(x_{2i-1})\}} \notag\\
 &- \frac{2}{n} \Big( \sum_{i=1}^{n/2} \hat g(x_{2i-1})+ \frac{1}{\gamma'}\sum_{i=1}^{n/2} \mathbf{1}_{\{Y_{2i-1}\ge\hat g(x_{2i-1})\}}-\sum_{i=1}^{n/2} f(x_{2i-1})\Big)^2 \notag\\
& - \frac{ (n/2) \Big(\sum_{i=1}^{n/2} \big(\hat g(x_{2i-1})+ \frac{1}{\gamma'} \mathbf{1}_{\{Y_{2i-1} \ge \hat g(x_{2i-1})\}}-f(x_{2i-1})\big) \big(x_{2i-1}-\frac{S}{n/2}\big) \Big)^2 }{ Rn/2-S^2} \label{e:T_f}
\end{align}
for any affine function $f:\mathbb{R}\to \mathbb{R}$. 

\begin{lemma}\label{l:T_f}
If $f:\mathbb{R}\to\mathbb{R}$ is an affine function, then
\[ T=T_f. \]
\end{lemma}

In particular, under $H_0$ we can write $T=S_1+S_2+S_3$, where
\begin{align*}
S_1 & := \sum_{i=1}^{n/2}  (\hat g(x_{2i-1})-g(x_{2i-1}))^2  + \frac{2}{\gamma'} \sum_{i=1}^{n/2} \big(Y_{2i-1}-g(x_{2i-1})\big)\mathbf{1}_{\{Y_{2i-1} \ge \hat g(x_{2i-1})\}} \\
S_2 & := \frac{2}{n}\big(\sum_{i=1}^{n/2} \big(\hat g(x_{2i-1})-g(x_{2i-1}) + \frac{1}{\gamma'}  \mathbf{1}_{\{Y_{2i-1} \ge \hat g(x_{2i-1})\}}\big) \big)^2 \\
S_3 & :=  \frac{ (n/2) \cdot \big(\sum_{i=1}^{n/2} \big(\hat g(x_{2i-1})-g(x_{2i-1}) + \frac{1}{\gamma'}  \mathbf{1}_{\{Y_{2i-1} \ge \hat g(x_{2i-1})\}}\big) (x_{2i-1}-\tfrac{2S}{n})\big)^2 }{Rn/2-S^2}
\end{align*}
and where $g$ is the true function. Sometimes it will be convenient to split the first summand further as $S_1=S_1'+2/\gamma' S_1''$, where
\begin{align*}
S_1' & := \sum_{i=1}^{n/2}  (\hat g(x_{2i-1})-g(x_{2i-1}))^2\\
S_1'' & := \sum_{i=1}^{n/2} \big(Y_{2i-1}-g(x_{2i-1})\big) \mathbf{1}_{\{Y_{2i-1} \ge \hat g(x_{2i-1})\}}.
\end{align*}

\prf[ of Lemma \ref{l:T_f}] We abbreviate
\[ \delta_i:= \frac{1}{\gamma'} \mathbf{1}_{\{ Y_i \ge \hat g(x_i) \} }, \quad \eta_i = \hat g(x_i)-f(x_i), \quad i=1, \dots, n, \quad  U = \sum_{i=1}^{n/2} \eta_{2i-1} x_{2i-1} \ \mbox{and} \ V = \sum_{i=1}^{n/2} \eta_{2i-1}.\]

Let $m,c\in\mathbb{R}$ be the numbers with $f(x)=mx+c$. Then
\[ \sum_{i=1}^{n/2} f(x_{2i-1})(x_{2i-1} - \tfrac{2S}{n}) = \sum_{i=1}^{n/2} mx_{2i-1}^2 + c x_{2i-1} -\tfrac{2S}{n}mx_{2i-1} -\tfrac{2S}{n}c = m \big(R-\tfrac{2S^2}{n} \big) + c \big(S-\tfrac{n}{2}\cdot \tfrac{2S}{n} \big) = m \big(R-\tfrac{2S^2}{n} \big). \]
Thus
\begin{align*}
T- T_f &= \begin{aligned}[t] \sum_{i=1}^{n/2} &f(x_{2i-1})^2 + 2 \sum_{i=1}^{n/2} f(x_{2i-1}) \eta_{2i-1} + 2 \sum_{i=1}^{n/2} f(x_{2i-1}) \delta_{2i-1}\\
& - \frac{4}{n} \big(\sum_{i=1}^{n/2} \eta_{2i-1} +\frac{1}{\gamma'} \sum_{i=1}^{n/2} \mathbf{1}_{\{Y_{2i-1} \ge \hat g(x_{2i-1}) \} } \big)\cdot \big(\sum_{i=1}^{n/2} f(x_{2i-1}) \big) - \frac{2}{n} \big( \sum_{i=1}^{n/2} f(x_{2i-1}) \big)^2 \\
& - 2\frac{ n\big(\sum_{i=1}^{n/2} \eta_{2i-1}\big(x_{2i-1}-\frac{2S}{n}\big) +\frac{1}{\gamma'} \sum_{i=1}^{n/2} \mathbf{1}_{\{Y_{2i-1} \ge \hat g(x_{2i-1}) \} }\big(x_{2i-1}-\frac{2S}{n}\big) \big)\cdot \big(\sum_{i=1}^{n/2} f(x_{2i-1})\big(x_{2i-1}-\frac{2S}{n}) \big) }{ Rn-2S^2}\\
& - \frac{ n\big( \sum_{i=1}^{n/2} f(x_{2i-1})\big(x_{2i-1}-\frac{2S}{n}\big) \big)^2 }{ Rn-2S^2}\end{aligned} 
\\
&= \begin{aligned}[t] m^2R+&2mcS +c^2n/2 + 2 mU +2cV + 2m\sum_{i=1}^{n/2} x_{2i-1} \delta_{2i-1} +2c \sum_{i=1}^{n/2}\delta_{2i-1} \\
& - \frac{4}{n} \big(V + \sum_{i=1}^{n/2} \delta_{2i-1} \big)\cdot \big(mS+cn/2 \big) - \frac{2}{n} \big( mS+cn/2 \big)^2 \\
& - 2\frac{ n\big( \big(U-\frac{2S}{n}V\big) + \sum_{i=1}^{n/2} \delta_{2i-1} x_{2i-1} - \frac{2S}{n}\sum_{i=1}^{n/2}\delta_{2i-1} \big)\cdot \big( m \big(R-\frac{2S^2}{n} \big)  \big) }{ Rn-2S^2}\\
& - \frac{ n\big( m \big(R-\frac{2S^2}{n} \big) \big)^2 }{ Rn-2S^2}\end{aligned} 
\\
&=\begin{aligned}[t] \big(R - \tfrac{2}{n}S^2 &- \frac{n\big(R-\frac{2S^2}{n} \big)^2}{Rn-2S^2}\big) \cdot m^2 + \big(2S- \frac{4}{n}Sn/2  \big)\cdot mc + \big(n/2 - \frac{2}{n}(n/2)^2 \big) \cdot c^2 \\
&+ \big(2m - 2\frac{n m \big(R-\frac{2S^2}{n} \big) }{Rn-2S^2}\big) \cdot U + \big(2c - \frac{4}{n} \big(mS+cn/2 \big) + \frac{ \frac{4S}{n} nm \big(R-\frac{2S^2}{n} \big)   }{ Rn-2S^2}\big)\cdot V\\
& + \big(2m-2\frac{  nm \big(R-\frac{2S^2}{n} \big)  }{ Rn-2S^2}\big)\cdot \sum_{i=1}^{n/2} x_{2i-1} \delta_{2i-1}\\
& + \big(2c - \frac{4}{n} \big(mS+cn/2 \big) + \frac{ \frac{4S}{n}  nm \big(R-\frac{2S^2}{n} \big) }{Rn-2S^2} \big) \cdot \sum_{i=1}^{n/2} \delta_{2i-1}
\end{aligned} 
\\
&=0. \qquad \qed
\end{align*}

\subsection{Asymptotic behavior of the moments}\label{ss:moments}

In this section we are going to derive asymptotic bounds and orders for the moments of $T$. In order to do this, we first have to derive an asymptotic bound for the moments of $\hat g(x)-g(x)$.

\begin{prop}\label{p:asy_ex}
Assume that \eqref{e:ass_scatter}, \eqref{e:lim=gamma} and \eqref{e:ass_distr_tail} hold. Then we have
\[
\limsup_{n\to\infty} \gamma^k(hn)^{k}\mathbb{E}[|\hat g(x)-g(x)|^k] \le \tfrac{2}{C_x^k} \Gamma\big(1+k \big).
\]
Assume that \eqref{e:ass_scatter_up}, \eqref{e:lim=gamma} and \eqref{e:ass_distr_tail} hold. Then we have 
\[
\liminf_{n\to\infty} \gamma^k(hn)^{k}\mathbb{E}[|\hat g(x)-g(x)|^k] \ge \tfrac{1}{(2C'_x)^k} \Gamma\big(1+k \big)
\]
\end{prop}
  	
\prf Put $Z_1:=\max\{\epsilon_{2i} \mid x_{2i}\in [x-h,x) \}$ and $Z_2:=\max\{ \epsilon_{2i} \mid x_{2i}\in [x, x+h) \}$. Then we have $\hat g(x)-g(x) \ge \min\{Z_1, Z_2 \}$ and hence
\[ \mathbb{E}[|\hat g(x) - g(x)|^k] \le \mathbb{E}[ |Z_1|^k ] + \mathbb{E}[ |Z_2|^k ]. \]
From \cite[Proposition 2.1]{Res87} we get
\[ \lim_{n\to\infty} \gamma^k(\#\{x_i\})^{k}\mathbb{E}[ |Z_j|^k ] = \Gamma\big(1+k \big),\quad j=1,2, \]
which in view of \eqref{e:ass_scatter} implies the first assertion. 

On the other hand we have $\hat g(x)-g(x) \le \max\{Z_1, Z_2 \}$. Again, we get from \cite[Proposition 2.1]{Res87} that 
\[ \lim_{n\to\infty} \gamma^k(\#\{x_i\})^{k}\mathbb{E}[ |\max\{Z_1,Z_2\}|^k ] = \Gamma\big(1+k \big), \] 
which in view of \eqref{e:ass_scatter_up} implies the second assertion. \qed\medskip

\begin{lemma}\label{l:expectation}
Assume \eqref{e:ass_scatter}, \eqref{e:ass_distr} and \eqref{e:ass_distr_tail}. Let $\hat\gamma$ be an estimator for $\gamma$ fulfilling (G\ref{i:G1}) and (G\ref{i:G2}). Then
\[
\limsup_{n\to\infty} n^{3/2}h^{5/2}\mathbb{E}\big[ \sum_{i=1}^{n/2}  (\hat g(x_{2i-1})-g(x_{2i-1}))^2  + \frac{2}{\hat\gamma} \sum_{i=1}^{n/2} \big(Y_{2i-1}-g(x_{2i-1})\big) \mathbf{1}_{\{Y_{2i-1} \ge \hat g(x_{2i-1})\}} \big]  \le  \sqrt{12} \frac{\sqrt{\hat C}}{\gamma C_x^2}. 
\]
\end{lemma}
\prf  We have
\begin{align*}
\mathbb{E}\big[ &(\hat g(x_{i})-g(x_{i}))^2  + \tfrac{2}{\hat\gamma}  \big(Y_{i}-g(x_{i})\big) \mathbf{1}_{\{Y_{i} \ge \hat g(x_i)\}} \big] \\
&= \mathbb{E}\big[  (\hat g(x_{i})-g(x_{i}))^2  + \tfrac{2}{\gamma} \big(Y_{i}-g(x_{i})\big) \mathbf{1}_{\{Y_{i} \ge \hat g(x_{i})\}} \big] + \mathbb{E}\big[ \big( \tfrac{2}{\hat\gamma}-\tfrac{2}{\gamma} \big) \big(Y_{i}-g(x_{i})\big) \mathbf{1}_{\{Y_{i} \ge \hat g(x_{i})\}} \big].
\end{align*}
Using the fact that $\Ex[X \mathbf{1}_{\{X\ge s\}}] = -\int_{s}^0 \mathbb{P}(s\le X <t)  \, dt$ for any $(-\infty,0]$-valued random variable $X$, we get 
\begin{align*}
 \Ex\Big[ (\hat g(x_i) -g(x_i))^2 &+ \frac{2}{\gamma} \big(Y_i-g(x_i)\big)  \mathbf{1}_{\{Y_i \ge \hat g(x_i)\}}   \Big] \\
&=  \Ex\Big[ \Ex\big[ (\hat g(x_i) -g(x_i))^2 + \frac{2}{\gamma} \big(Y_i-g(x_i)\big)  \mathbf{1}_{\{Y_i \ge \hat g(x_i)\}} \mid \hat g \big]  \Big] \\
& =  \Ex\Big[ (\hat g(x_i) -g(x_i))^2 - \int_{\hat g(x_i)-g(x_i)}^{0}\frac{2}{\gamma} \big(F(y)-F(\hat g(x_i)-g(x_i)) \big) \, dy   \Big] \\
& = \begin{aligned}[t]\Ex\Big[ (\hat g(x_i)& - g(x_i))^2+ (\hat g(x_i)-g(x_i))^2 -2(\hat g(x_i)-g(x_i))^2  \\
-& \int_{\hat g(x_i)-g(x_i)}^{0}\frac{2}{\gamma} \big( F(y)-\gamma y-F(\hat g(x_i)-g(x_i)) + \gamma (\hat g(x_i)-g(x_i)) \big) \, dy  \Big] \end{aligned}
\\
& =   \Ex\Big[ -\int_{\hat g(x_i)-g(x_i)}^{0}\frac{2}{\gamma} \big(  (1-F(\hat g(x_i)-g(x_i))+\gamma (\hat g(x_i)-g(x_i)))- (1-F(y)+\gamma y) \big) \, dy\Big].  
\end{align*}
We have
\[ \Big|\int_{\hat g(x_i)-g(x_i)}^{0}\frac{2}{\gamma} (1-F(\hat g(x_i)-g(x_i))+\gamma (\hat g(x_i)-g(x_i))) \, dy \Big| \le \int_{\hat g(x_i)-g(x_i)}^{0}\frac{2}{\gamma} C_F |\hat g(x_i)-g(x_i)|^2 \, dy  = \frac{2C_F}{\gamma}  |\hat g(x_i)-g(x_i)|^{3} \]
and
\[ \Big| \int_{\hat g(x_i)-g(x_i)}^{0}\frac{2}{\gamma}  (1-F(y)+\gamma y) \, dy \Big| \le \int_{\hat g(x_i)-g(x_i)}^{0}\frac{2}{\gamma} C_F |y|^2 \, dy = \frac{2C_F}{3\gamma}  |\hat g(x_i)-g(x_i)|^{3}. \]
So
\begin{equation} 
\Ex\big[ (\hat g(x_i) -g(x_i))^2 + \frac{2}{\gamma} \big(Y_i-g(x_i)\big)  \mathbf{1}_{\{Y_i \ge \hat g(x_i)\}} \mid \hat g  \big] \le  \frac{8C_F}{3\gamma}  |\hat g(x_i)-g(x_i)|^{3}. \label{e:cond_bound} \end{equation}
Hence Proposition \ref{p:asy_ex} yields
\begin{align*}
\limsup_{n\to\infty} n^2h^{3} &\Ex\big[ \sum_{i=1}^{n/2} (\hat g(x_{2i-1}) -g(x_{2i-1}))^2 + \frac{2}{\gamma} \big(Y_{2i-1}-g(x_{2i-1})\big)  \mathbf{1}_{\{Y_{2i-1} \ge \hat g(x_{2i-1})\}}   \big] \\
& \le \limsup_{n\to\infty} n^{-1}(nh)^{3} \sum_{i=1}^{n/2} \frac{8C_F}{3\gamma} \Ex[ |\hat g(x_{2i-1}) - g(x_{2i-1})|^{3} ] \\
& \le  \frac{1}{2}\frac{8C_F}{3\gamma}  \frac{2}{C_x^3\gamma^{3}} \Gamma\big(4 \big) \\
&= \frac{C_F}{C_x^3} \frac{16}{\gamma^4}.
\end{align*}
Moreover, with a similar calculation as above, we get 
\begin{align*}
\big|\mathbb{E}\big[ \big(& \frac{2}{\hat\gamma}-\frac{2}{\gamma} \big) \big(Y_{i}-g(x_{i})\big) \mathbf{1}_{\{Y_{i} \ge \hat g(x_{i})\}} \big]\big| \\
&=\big|\mathbb{E}\big[ \mathbb{E}\big[ \big( \frac{2}{\hat\gamma}-\frac{2}{\gamma} \big) \big(Y_{i}-g(x_{i})\big) \mathbf{1}_{\{Y_{i} \ge \hat g(x_{i})\}} \mid\hat\gamma, \hat g\big] \big]\big|\\
&\le\mathbb{E}\Big[ \big| \frac{2}{\hat\gamma}-\frac{2}{\gamma} \big| \cdot \Big| \frac{\gamma}{2}(\hat g(x_i)-g(x_i))^2 + \int_{\hat g(x_i)-g(x_i)}^{0} \big( (1-F(\hat g(x_i)-g(x_i))+\gamma (\hat g(x_i)-g(x_i))) - (1-F(y)+\gamma y)\big) \, dy\Big|\Big]\\
&\le \mathbb{E}\Big[ \big| \frac{2}{\hat\gamma}-\frac{2}{\gamma} \big| \cdot \Big| \frac{\gamma}{2}(\hat g(x_i)-g(x_i))^2 + \frac{4C_F}{3} |\hat g(x_i)-g(x_i)|^3\Big|\Big]\\
&\le \sqrt{\mathbb{E}\Big[ \big( \frac{2}{\hat\gamma} - \frac{2}{\gamma}\big)^2 \Big]} \cdot \sqrt{ \Ex\Big[  \Big( \frac{\gamma}{2}(\hat g(x_i)-g(x_i))^2 + \frac{4C_F}{3} |\hat g(x_i)-g(x_i)|^3\Big)^2\Big]}.
\end{align*}
Hence Proposition \ref{p:asy_ex} yields
\begin{align*}
\limsup_{n\to\infty} n^{3/2}h^{5/2} &\Ex\Big[  \big( \frac{2}{\hat\gamma}-\frac{2}{\gamma} \big) \sum_{i=1}^{n/2}  \big(Y_{2i-1}-g(x_{2i-1})\big)  \mathbf{1}_{\{Y_{2i-1} \ge \hat g(x_{2i-1})\}}   \Big] \\
& \le \limsup_{n\to\infty} \frac{(nh)^{1/2}}{n}(nh)^{2} \sum_{i=1}^{n/2} \Big(    \sqrt{\mathbb{E}\Big[ \big( \frac{2}{\hat\gamma} - \frac{2}{\gamma}\big)^2 \Big]} \cdot \sqrt{ \Ex\Big[  \Big( \frac{\gamma}{2}(\hat g(x_i)-g(x_i))^2 + \frac{4C_F}{3} |\hat g(x_i)-g(x_i)|^3\Big)^2\Big]} \Big)\\
& \le \tfrac{1}{2} \sqrt{4\hat C}\cdot \sqrt{\frac{2\gamma^2}{4C_x^4\gamma^4}\Gamma(5)}\\
&=  \sqrt{12} \frac{\sqrt{\hat C}}{\gamma C_x^2}. \qquad \qed
\end{align*}

\begin{lemma}\label{l:variances}
Assume that \eqref{e:ass_scatter}, \eqref{e:ass_scatter_up}, \eqref{e:ass_distr} and \eqref{e:ass_distr_tail} hold. Let $k_1=k_1(n)$ and $k_2=k_2(n)$ be sequences in $\{1, \dots, n/2\}$ with $k_1(n)<k_2(n)$ for all $n\in\mathbb{N}$. Then
\begin{align*}
\limsup_{n\to\infty} \frac{(nh)^{3}}{k_2-k_1+1} \Var\Big(  \sum_{i=k_1}^{k_2} (\hat g(x_{2i-1})-g(x_{2i-1}))^2 \Big) & \le  \frac{192C_x'}{\gamma^4C_x^4} \\ 
\limsup_{n\to\infty} n^{2}h^{3} \Var\Big(\sum_{i=1}^{n/2} (\hat g(x_{2i-1}) - g(x_{2i-1}))^2 \Big) & \le  \frac{96C_x'}{\gamma^4C_x^4} \\ 
 \limsup_{n\to\infty} \frac{(nh)^{3}}{k_2-k_1+1} \Var\Big(  \sum_{i=k_1}^{k_2} (Y_{2i-1}- g(x_{2i-1})) \mathbf{1}_{\{Y_{2i-1} \ge \hat g(x_{2i-1})\}} \Big) & \le \frac{48 C_x'}{C_x^4\gamma^2} + \frac{4}{C_x^3\gamma^2} \\ 
\limsup_{n\to\infty} n^{2}h^{3} \Var\Big(\sum_{i=1}^{n/2} (Y_{2i-1}- g(x_{2i-1})) \mathbf{1}_{\{Y_{2i-1} \ge \hat g(x_{2i-1})\}} \Big) & \le \frac{24C_x'}{C_x^4\gamma^2 } + \frac{2}{C_x^3\gamma^2} \\ 
\liminf_{n\to\infty} \frac{(nh)^{3}}{k_2-k_1+1}   \Ex\Big[\Var\Big(  \sum_{i=k_1}^{k_2} (Y_{2i-1}- g(x_{2i-1})) \mathbf{1}_{\{Y_{2i-1} \ge \hat g(x_{2i-1})\}} \mid\hat g \Big) \Big] & \ge \frac{3}{256(C_x')^3\gamma^2}\\
\limsup_{n\to\infty} n^{2}h^{3} \Ex\Big[\Var\Big(  \sum_{i=1}^{n/2} (Y_{2i-1}- g(x_{2i-1})) \mathbf{1}_{\{Y_{2i-1} \ge \hat g(x_{2i-1})\}} \mid \hat g\Big) \Big] & \ge \frac{3}{512(C_x')^3\gamma^2}\\
\end{align*}
\end{lemma}

\prf 
We have
\begin{align*} 
\Var\Big( \sum_{i=k_1}^{k_2}  (&\hat g(x_{2i-1}) -g(x_{2i-1}))^2 \Big) = \sum_{i=k_1}^{k_2} \sum_{j\in\{k_1, \dots, k_2\}:\atop |x_{2j-1}-x_{2i-1}| < 2h}\Cov\big((\hat g(x_{2i-1}) -g(x_{2i-1}))^2, (\hat g(x_{2j-1}) -g(x_{2j-1}))^2\big) \\
&\le (k_2-k_1+1) \max\big\{\#\{i\mid x_{2i-1}\in (x-2h,x+2h) \cap [0,1]\} \mid x\in (0, 1)\big\} \mathbb{E}[(\hat g(x_{2i-1}) - g(x_{2i-1}))^4]. 
\end{align*}
Hence Proposition \ref{p:asy_ex} yields
\[ \limsup_{n\to\infty} \frac{(nh)^{3}}{k_2-k_1+1} \Var\Big(  \sum_{i=k_1}^{k_2} (\hat g(x_{2i-1})-g(x_{2i-1}))^2 \Big) \le 4C'_x \cdot \frac{2\Gamma(5)}{\gamma^4C_x^4} = \frac{192 C_x'}{\gamma^4C_x^4}. \]

In order to show the third assertion, we use 
\begin{align*} \Var\Big(  &\sum_{i=k_1}^{k_2} (Y_{2i-1}- g(x_{2i-1})) \mathbf{1}_{\{Y_{2i-1} \ge \hat g(x_{2i-1})\}} \Big) \\
&= \Var\Big( \Ex\Big[ \sum_{i=k_1}^{k_2} (Y_{2i-1}- g(x_{2i-1})) \mathbf{1}_{\{Y_{2i-1} \ge \hat g(x_{2i-1})\}} \mid \hat g\Big]\Big) + \Ex\Big[\Var\Big(  \sum_{i=k_1}^{k_2} (Y_{2i-1}- g(x_{2i-1})) \mathbf{1}_{\{Y_{2i-1} \ge \hat g(x_{2i-1})\}} \mid \hat g \Big) \Big].
\end{align*}

Similarly to \eqref{e:cond_bound} we get 
\[ \Big| \Ex[(Y_{2i-1}- g(x_{2i-1})) \mathbf{1}_{\{Y_{2i-1} \ge \hat g(x_{2i-1})\}} \mid \hat g \big] \Big| \le \frac{\gamma}{2} |\hat g(x_{2i-1}) - g(x_{2i-1}) |^{2} + \frac{4C_F}{3} |\hat g(x_{2i-1}) - g(x_{2i-1}) |^{3}.\]
 
Thus
\begin{align} 
 \limsup_{n\to\infty} \frac{(nh)^{3}}{k_2-k_1+1}\Var\Big(& \Ex\big[ \sum_{i=k_1}^{k_2} (Y_{2i-1}- g(x_{2i-1})) \mathbf{1}_{\{Y_{2i-1} \ge \hat g(x_{2i-1})\}} \mid \hat g \big] \Big)\notag\\ 
&\le \begin{aligned}[t] \limsup_{n\to\infty} (nh)^4 &\frac{\max\big\{\#\{i\mid x_{2i-1}\in (x-2h,x+2h) \cap [0,1]\} \mid x\in (0, 1)\big\}}{nh} \\ &\max\Big\{\Ex\Big[\Big| \Ex[(Y_{2i-1}- g(x_{2i-1})) \mathbf{1}_{\{Y_{2i-1} \ge \hat g(x_{2i-1})\}} \mid \hat g \big] \Big|^2 \Big] \mid i=1,\dots, \frac{n}{2} \Big\}\end{aligned}\notag\\
& \le  \limsup_{n\to\infty} 4C_x'(nh)^{4} (\gamma/2)^2\max\Big\{\Ex\big[ |\hat g(x_{2i-1}) - g(x_{2i-1}) |^{4} \big] \mid i=1,\dots, \frac{n}{2} \Big\} \notag\\
& \le 4C_x' \frac{\gamma^2}{4} \frac{48}{\gamma^4C_x^4} 
= \frac{48 C_x'}{\gamma^2C_x^4}. \label{e:lim_Var_S12_1}
\end{align}

Moreover,
\begin{align*}
\Var\big( (Y_{2i-1}- g(x_{2i-1})) & \mathbf{1}_{\{Y_{2i-1} \ge \hat g(x_{2i-1})\}}  \mid \hat g\big) \\
&\le \Ex\big[ (Y_{2i-1}- g(x_{2i-1}))^2 \mathbf{1}_{\{Y_{2i-1} \ge \hat g(x_{2i-1})\}}  \mid \hat g\big] \\
&= \int_0^\infty \mathbb{P}\big( (Y_{2i-1}- g(x_{2i-1}))^2 \cdot \mathbf{1}_{\{Y_{2i-1}\ge \hat g(x_{2i-1})\}} >t\mid \hat g\big) \, dt \\
&= \int_0^{(\hat g(x_{2i-1})-g(x_{2i-1}))^2} \mathbb{P}\big(-\sqrt{t} > Y_{2i-1}- g(x_{2i-1}) \ge \hat g(x_{2i-1}) - g(x_{2i-1}) \mid \hat g\big) \, dt \\
&= \int_0^{(\hat g(x_{2i-1})-g(x_{2i-1}))^2} F(-\sqrt{t})-F(\hat g(x_{2i-1}) -g(x_{2i-1})) \, dt \\
&\le \int_0^{(\hat g(x_{2i-1})-g(x_{2i-1}))^2} \gamma |\hat g(x_{2i-1})-g(x_{2i-1}) +\sqrt{t}| +C_F ((\hat g(x_{2i-1})-g(x_{2i-1}))^2 +\sqrt{t}^2)\, dt \\
&\le \gamma\big( |\hat g(x_{2i-1})-g(x_{2i-1})|^3 -\tfrac{2}{3}\big( (\hat g(x_{2i-1})-g(x_{2i-1}))^2 \big)^{3/2} \big) + \tfrac{3}{2}C_F(\hat g(x_{2i-1})-g(x_{2i-1}))^4 \\
&= \frac{\gamma}{3}  |\hat g(x_{2i-1})-g(x_{2i-1})|^3 + \tfrac{3}{2}C_F (\hat g(x_{2i-1})-g(x_{2i-1}))^4.
\end{align*}
So we get from Proposition \ref{p:asy_ex}
\begin{align} 
\limsup_{n\to\infty} \frac{(nh)^{3}}{k_2-k_1+1} \Ex\Big[\Var\Big( &\sum_{i=k_1}^{k_2} (Y_{2i-1}- g(x_{2i-1})) \mathbf{1}_{\{Y_{2i-1} \ge \hat g(x_{2i-1})\}} \mid\hat g \Big) \Big]\notag\\
&\le \limsup_{n\to\infty} \frac{(nh)^{3}}{k_2-k_1+1}\frac{\gamma}{3}  \sum_{i=k_1}^{k_2} \Ex\big[|\hat g(x_{2i-1})-g(x_{2i-1})|^3 \big] \notag \\
& \le \frac{\gamma}{3} \Big( \frac{2}{C_x^3\gamma^3} \Gamma\big(4\big) \Big)\notag\\
&= \frac{4}{C_x^3\gamma^2}. \label{e:lim_Var_S12_2} 
\end{align}
Now \eqref{e:lim_Var_S12_1} and \eqref{e:lim_Var_S12_2} yield the third assertion. 

On the other hand we have
\begin{align*}
\Var\big( (Y_{2i-1}- g(x_{2i-1})) &\mathbf{1}_{\{Y_{2i-1} \ge \hat g(x_{2i-1})\}}  \mid \hat g\big)\\
 &\ge 
\begin{aligned}[t]\frac{(\hat g(x_{2i-1}) - g(x_{2i-1}))^2}{16} \min\big\{ \mathbb{P}\big( |Y_{2i-1} - g(x_{2i-1})| \le \tfrac{1}{4} | \hat g(x_{2i-1})-g(x_{2i-1})|\mid \hat g \big),&\\
 \mathbb{P}\big(\tfrac{3}{4} | \hat g(x_{2i-1})-g(x_{2i-1})| \le  |Y_{2i-1} - g(x_{2i-1})| \le | \hat g(x_{2i-1})-g(x_{2i-1})| \mid \hat g \big)& \big\} \end{aligned}\\
& \ge \frac{(\hat g(x_{2i-1}) - g(x_{2i-1}))^2}{16} \cdot \big(\frac{\gamma}{4} | \hat g(x_{2i-1})-g(x_{2i-1})| - 2C_F | \hat g(x_{2i-1})-g(x_{2i-1})|^2 \big)
\end{align*}
and hence
\begin{align*} \liminf_{n\to\infty} \frac{(nh)^{3}}{k_2-k_1+1}\Ex\Big[ \Var\Big( & \sum_{i=k_1}^{k_2} (Y_{2i-1}- g(x_{2i-1})) \mathbf{1}_{\{Y_{2i-1} \ge \hat g(x_{2i-1})\}} \mid \hat g \Big) \Big]\\
&\ge \liminf_{n\to\infty} \frac{(nh)^{3}}{k_2-k_1+1} \Ex\Big[ \sum_{i=k_1}^{k_2}\frac{(\hat g(x_{2i-1}) - g(x_{2i-1}))^2}{16} \frac{\gamma}{4} |\hat g(x_{2i-1})-g(x_{2i-1})| \Big]\\
& \ge \frac{\gamma}{64}\frac{\Gamma(4)}{(2C_x'\gamma)^3} = \frac{3}{256(C_x')^3\gamma^2}. 
\end{align*}
So the fifth assertion follows. Putting $k_1=1$ and $k_2=n/2$ we obtain the second, forth and sixth assertion. 
\qed\medskip

We get the following corollary of Lemma \ref{l:variances}.
\begin{cor}\label{c:VarS1} Assume \eqref{e:ass_scatter}, \eqref{e:ass_scatter_up}, \eqref{e:ass_distr} and \eqref{e:ass_distr_tail}. Let $\hat\gamma$ be an estimator for $\gamma$ that fulfills (G\ref{i:G1}) and (G\ref{i:G3}). Then we have
\begin{align*} 
\limsup_{n\to\infty} n^{2}h^{3} \Var(S_1)  \le  \frac{8}{C_x^3\gamma^4} . 
\end{align*}
\end{cor}
\prf For two random variables $A$ and $B$ and $\delta>0$ we have
\[ \Var(A+B)=\Var(A)+2\Cov(A,B)+\Var(B) \le \Var(A) + 2 \sqrt{\delta\Var(A)}\sqrt{\tfrac{1}{\delta}\Var(B)}+ \Var(B) \le (1+\delta)\Var(A) + \big(1+\tfrac{1}{\delta}\big) \Var(B). \]
So let $\delta>0$. We have
\begin{align*}
\Var(S_1) &= \Var\big(S_1'+\tfrac{2}{\gamma} S_1''+(\tfrac{2}{\hat\gamma}-\tfrac{2}{\gamma})S_1''\big)\\
&\le \begin{aligned}[t](1+\delta)\cdot \Ex\big[\Var(S_1'+\tfrac{2}{\gamma}S_1''\mid \hat g)\big] + (1+\delta)\cdot \Var\big(\Ex[S_1'+\tfrac{2}{\gamma}S_1''\mid\hat g]\big)  & + \big(1+\tfrac{1}{\delta}\big) \cdot \Ex\big[\Var\big((\tfrac{2}{\hat\gamma}-\tfrac{2}{\gamma})\cdot S_1'' \mid \hat g, \hat\gamma\big)\big] \\&+ \big(1+\tfrac{1}{\delta}\big)\cdot \Var\big( \Ex\big[ (\tfrac{2}{\hat\gamma}-\tfrac{2}{\gamma})\cdot S_1'' \mid \hat g, \hat\gamma\big]\big). \end{aligned}
\end{align*} 
Now 
\[ \Var\big(S_1'+\tfrac{2}{\gamma}S_1''\mid \hat g\big) = \frac{4}{\gamma^2}\Var(S_1'' \mid\hat g)\]
and hence 
\[ \limsup_{n\to\infty} n^{-1}(nh)^{3} \Ex\big[\Var\big(S_1'+\tfrac{2}{\gamma}S_1''\mid \hat g\big) \big] \le  \frac{4}{\gamma^2}\cdot \frac{2}{C_x^3\gamma^2} = \frac{8}{C_x^3\gamma^4}, \]
can be obtained as in the proof of Lemma \ref{l:variances}. From \eqref{e:cond_bound} we get
\[\Var\big(\Ex[S_1'+\tfrac{2}{\gamma}S_1''\mid\hat g]\big) \in O(n^2h(nh)^{-6}).\]
Moreover, we get
\begin{align*}
\limsup_{n\to\infty} n^{2}h^{3}& \Ex\big[\Var\big((\tfrac{2}{\hat\gamma}-\tfrac{2}{\gamma})\cdot S_1'' \mid \hat g, \hat\gamma\big)\big] \\
&\le \limsup_{n\to\infty} n^{-1}(nh)^{3} \Ex\Big[(\tfrac{2}{\hat\gamma}-\tfrac{2}{\gamma})^2 \cdot \sum_{i=1}^{n/2}\tfrac{\gamma}{3} |\hat g(x_{2i-1})-g(x_{2i-1})|^3\Big]\\
&\le \limsup_{n\to\infty} (nh)^{3} \tfrac{1}{2} \sqrt{\Ex\big[(\tfrac{2}{\hat\gamma}-\tfrac{2}{\gamma})^4\big]} \cdot \max\Big\{\sqrt{\Ex\big[\tfrac{\gamma^2}{9} |\hat g(x_{2i-1})-g(x_{2i-1})|^6\big]} \mid i=1, \dots, \frac{n}{2} \Big\}\\
&=0
\end{align*}
and
\begin{align*}
\limsup_{n\to\infty} n^{2}h^{3}& \Var\big(\Ex\big[(\tfrac{2}{\hat\gamma}-\tfrac{2}{\gamma})\cdot S_1'' \mid \hat g, \hat\gamma\big]\big) \\
&=\limsup_{n\to\infty} n^{-1}(nh)^{3} \Var\Big((\tfrac{2}{\hat\gamma}-\tfrac{2}{\gamma})\cdot \tfrac{\gamma}{2} \sum_{i=1}^{n/2} (\hat g(x_{2i-1})- g(x_{2i-1}))^2\Big) \\
&\le\limsup_{n\to\infty} \begin{aligned}[t](nh)^{3} \tfrac{1}{2}&\max\big\{ \#\{i\in\{1,\dots, n/2\} \mid x_{2i-1}\in (x-2h, x+2h) \} \mid x\in[0,1] \big\} \\ &\cdot \max\Big\{ \Cov\Big( (\tfrac{2}{\hat\gamma}-\tfrac{2}{\gamma})\cdot \tfrac{\gamma}{2}  (\hat g(x_{2i-1})- g(x_{2i-1}))^2, (\tfrac{2}{\hat\gamma}-\tfrac{2}{\gamma})\cdot \tfrac{\gamma}{2}  (\hat g(x_{2j-1})- g(x_{2j-1}))^2\Big) \mid\\& i,j\in\{1,\dots, \frac{n}{2} \} \Big\} \end{aligned} \\ 
&\le\limsup_{n\to\infty} (nh)^{4} 2C_x \max\Big\{\Ex\big[(\tfrac{2}{\hat\gamma}-\tfrac{2}{\gamma})^2\cdot \tfrac{\gamma^2}{4} (\hat g(x_{2i-1})- g(x_{2i-1}))^4\big] \mid i=1, \dots, \frac{n}{2} \Big\} \\
&\le\limsup_{n\to\infty} (nh)^{4} 2C_x \sqrt{\Ex\big[(\tfrac{2}{\hat\gamma}-\tfrac{2}{\gamma})^4\big]}\cdot \max\Big\{\sqrt{\Ex[\tfrac{\gamma^4}{16} (\hat g(x_{2i-1})- g(x_{2i-1}))^8\big]} \mid i=1, \dots, \frac{n}{2} \Big\} \\
&=0.
\end{align*}
Hence we get
\[ \limsup_{n\to\infty} n^{-1}(nh)^{3} \Var(S_1)  \le (1+\delta)\cdot \frac{8}{C_x^3\gamma^4} \] 
Since $\delta>0$ was arbitrary, the assertion follows. \qed

\begin{prop}\label{p:T1prop}
Assume that \eqref{e:ass_scatter}, \eqref{e:ass_scatter_up}, \eqref{e:lim=gamma} and \eqref{e:ass_distr_tail} hold. Then
\begin{enumerate}[(i)]
\item $\mathbb{E}\, T_1 \in O(n^{-1}h^{-2})$
\item $\Var(T_1) \in O(n^{-2}h^{-3})$.
\end{enumerate}
\end{prop}

\prf (i) From Proposition \ref{p:asy_ex} we get
\[\Ex\big[ \sum_{i=1}^{n/2} (\hat g(x_{2i-1})-g(x_{2i-1}))^2 \big] \in O(n^{-1}h^{-2}). \]
Since $0\le T_1 \le \sum_{i=1}^{n/2} (\hat g(x_{2i-1})-g(x_{2i-1}))^2$, the first assertion follows. 

(ii) We use that by Proposition \ref{p:T1_explicit} we have
\begin{equation}
T_1= \begin{aligned}[t]
\sum_{i=1}^{n/2} \big(\hat g(x_{2i-1})-g(x_{2i-1})\big)^2  &- \frac{2}{n} \cdot \Big( \sum_{i=1}^{n/2} \big(\hat g(x_{2i-1}) - g(x_{2i-1}) \big) \Big)^2\\
& - \frac{(n/2) \cdot \Big(\sum_{i=1}^{n/2} \big(\hat g(x_{2i-1}) -g(x_{2i-1}) \big) \cdot \big(x_{2i-1}-\frac{S}{n/2}\big)  \Big)^2}{Rn/2-S^2}.\end{aligned} \label{e:T1}\end{equation}
Observe that
\begin{align*}
 \Var\Big( \sum_{i=1}^{n/2} (\hat g(x_{2i-1}) - g(x_{2i-1}))^2 \Big) &= \sum_{i=1}^{n/2} \sum_{j=1}^{n/2} \Cov\big( (\hat g(x_{2i-1})-g(x_{2i-1}))^2, (\hat g(x_{2j-1})-g(x_{2j-1}))^2 \big). \\
&= \sum_{i=1}^{n/2} \#\{j \mid |x_{2i-1}-x_{2j-1}|<2h\} \cdot \max\{ \Ex[(\hat g(x_{2\iota-1})-g(x_{2\iota-1}))^4] \mid \iota=1, \dots, \tfrac{n}{2} \} \\
& \in O(n^2h(nh)^{-4})= O(n^{-2}h^{-3}). \end{align*}
In order to treat the second and third summand of \eqref{e:T1}, we let $w:[0,1]\to[-1,1]$ be a measurable function, e.g.\ $w\equiv 1$ or $w(x)=x-2S/n$. We have
\begin{align*}
 \Var\Big( \Big( \sum_{i=1}^{n/2}& (\hat g(x_{2i-1}) - g(x_{2i-1})) \cdot w(x_{2i-1}) \Big)^2 \Big)\\
& = \sum_{{i_1}=1}^{n/2} \sum_{{i_2}=1}^{n/2} \sum_{{j_1}=1}^{n/2} \sum_{{j_2}=1}^{n/2} \begin{aligned}[t]\Cov \big( &(\hat g(x_{2i_1-1})-g(x_{2i_1-1})) \cdot w(x_{2i_1-1}) \cdot (\hat g(x_{2i_2-1})-g(x_{2i_2-1})) \cdot w(x_{2i_2-1}), \\&(\hat g(x_{2j_1-1})-g(x_{2j_1-1})) \cdot w(x_{2j_1-1}) \cdot (\hat g(x_{2j_2-1})-g(x_{2j_2-1})) \cdot w(x_{2j_2-1}) \big). \end{aligned}
\end{align*}
Now
\begin{align*}
\Cov \big( &(\hat g(x_{2i_1-1})-g(x_{2i_1-1})) \cdot w(x_{2i_1-1}) \cdot (\hat g(x_{2i_2-1})-g(x_{2i_2-1})) \cdot w(x_{2i_2-1}),\\
& (\hat g(x_{2j_1-1})-g(x_{2j_1-1})) \cdot w(x_{2j_1-1}) \cdot (\hat g(x_{2j_2-1})-g(x_{2j_2-1})) \cdot w(x_{2j_2-1}) \big) =0 
\end{align*}
if $|x_{2i_1-1}-x_{2j_1-1}|>2h$, $|x_{2i_1-1}-x_{2j_2-1}|>2h$, $|x_{2i_2-1}-x_{2j_1-1}|>2h$ and $|x_{2i_2-1}-x_{2j_2-1}|>2h$. Since, moreover
\begin{align*}
\Cov \big( (\hat g(x_{2i_1-1})&-g(x_{2i_1-1})) \cdot w(x_{2i_1-1}) \cdot (\hat g(x_{2i_2-1})-g(x_{2i_2-1})) \cdot w(x_{2i_2-1}),\\
& (\hat g(x_{2j_1-1})-g(x_{2j_1-1})) \cdot w(x_{2j_1-1}) \cdot (\hat g(x_{2j_2-1})-g(x_{2j_2-1})) \cdot w(x_{2j_2-1}) \big)  \\
&\le \begin{aligned}[t]\max\big\{ \Ex[(\hat g(x_{2i_1-1}) - g(x_{2i_1-1}))^4], \Ex[(\hat g(x_{2i_2-1}) - g(x_{2i_2-1}))^4], \Ex[(\hat g(x_{2j_1-1}) - g(x_{2j_1-1}))^4],\\ \Ex[(\hat g(x_{2j_2-1}) - g(x_{2j_2-1}))^4] \big\}\end{aligned},
\end{align*}
we get 
\[ \Var\Big( \Big( \sum_{i=1}^{n/2} (\hat g(x_{2i-1}) - g(x_{2i-1})) \cdot w(x_{2i-1}) \Big)^2 \Big) \in O(n^4h(nh)^{-4})=O(h^{-3}). \]
Moreover we have
\begin{align*} 
R \frac{n}{2}-S^2=\sum_{i=1}^{n/2}\sum_{j=1}^{n/2} x_{2i-1}^2 - x_{2i-1}x_{2j-1} &= \sum_{i=1}^{n/2} \sum_{j=i+1}^{n/2} x_{2i-1}^2+x_{2j-1}^2-2x_{2i-1}x_{2j-1} \\
&= \sum_{i=1}^{n/2} \sum_{j=i+1}^{n/2} (x_{2i-1} - x_{2j-1})^2\\
& \ge \#\big\{i \mid x_{2i-1}\le \tfrac{1}{4}\big\} \cdot  \#\big\{j \mid x_{2j-1}\ge \tfrac{3}{4}\big\} \cdot \frac{1}{4}\\
& \ge \frac{C_x^2}{64} n^2.
\end{align*}
So
\[\Var\Big(\frac{2}{n} \cdot \Big( \sum_{i=1}^{n/2} \big(\hat g(x_{2i-1}) - g(x_{2i-1}) \big) \Big)^2 \Big) \in O(n^{-2}h^{-3}) \]
and
\[ \Var\Big(\frac{(n/2) \cdot \Big(\sum_{i=1}^{n/2} \big(\hat g(x_{2i-1}) -g(x_{2i-1}) \big) \cdot \big(x_{2i-1}-\frac{S}{n/2}\big)  \Big)^2}{Rn/2-S^2} \Big) \in O(n^{-2}h^{-3}). \] 
Putting the pieces together we obtain the second assertion. \qed

\subsection{Central limit theorem}\label{s:CLT}

Our aim in this subsection is to derive a central limit theorem for the test statistic $T$. We will derive a central limit theorem for $S_1$,  and show that the other two summands in \eqref{e:T_f} are asymptotically neglectable. Our method is to employ \cite[Theorem~3]{Ph69}, which is a Lindeberg-Feller type central limit theorem for mixing sequences. 

\begin{theorem}\label{T:gCLT}
If the sample points fulfill \eqref{e:ass_scatter} and \eqref{e:ass_scatter_up}, the errors satisfy \eqref{e:ass_distr} and \eqref{e:ass_distr_tail}, and the sequences $\Var(S_1')n^2h^3$, $\Cov(S_1',S_1'')n^2h^3$ and $\Var(S_1'')n^2h^3$ converge, then
\[ \frac{ (S_1', S_1'') - \Ex \big[(S_1', S_1'')\big]}{\sqrt{n^{-2}h^{-3}}} \longrightarrow \mathcal{N}(0,\Sigma) \]
in distribution for the matrix $\Sigma$ of the variance and covariance limits. 
\end{theorem}  

\begin{prop}\label{p:S2S3negl}
If the sample points fulfill \eqref{e:ass_scatter}, if the errors satisfy \eqref{e:ass_distr} and if $\hat\gamma$ is an estimator for $\gamma$ with (G\ref{i:G1}) and (G\ref{i:G4}), then there is a sequence $b_n\in O(n^{-1/2}h^{-3/2})$ with 
\begin{align*}
S_2-\frac{b_n^2}{n} &\in o_P(n^{-1}h^{-3/2}) = o_P(\sqrt{\Var S_1})\\
S_3-\frac{b_n^2}{n} &\in o_P(n^{-1}h^{-3/2}) = o_P(\sqrt{\Var S_1}).
\end{align*}
\end{prop}

\begin{cor}\label{c:CLT_T}
Assume that the sample points fulfill \eqref{e:ass_scatter} and \eqref{e:ass_scatter_up}, the errors satisfy \eqref{e:ass_distr} and \eqref{e:ass_distr_tail}, that $\hat\gamma$ is an estimator for $\gamma$ with (G\ref{i:G1}) and (G\ref{i:G4}) and that $\lim_{n\to\infty} nh^2=\infty$. Then we have 		
\[ \frac{ T }{\sqrt{\Var(S_1)}} \longrightarrow \mathcal{N}(0,1) \]
as $n\to\infty$ in distribution. 
\end{cor}		
\prf From Theorem \ref{T:gCLT} and Slutzky's theorem we get that every subsequence contains a further subsequence along which
\[ \frac{ S_1 - \Ex [S_1]}{\sqrt{\Var(S_1)}} \longrightarrow \mathcal{N}(0,1).\]
Hence convergence holds in total.
Together with Proposition \ref{p:S2S3negl} this implies 
\[ \frac{ T - \Ex [S_1]}{\sqrt{\Var(S_1)}} \longrightarrow \mathcal{N}(0,1). \]
By Lemma \ref{l:expectation} and Lemma \ref{l:variances} the assertion follows. 
\qed\medskip

\prf[ of Proposition \ref{p:S2S3negl}]
In order to show the first assertion
abbreviate
\[ R_2:= \sum_{i=1}^{n/2} \hat g(x_{2i-1}) - g(x_{2i-1}) + \frac{1}{\hat \gamma} \mathbf{1}_{\{Y_{2i-1} \ge \hat g(x_{2i-1})\}} \]
so that $S_2=2R_2^2/n$. Put $a_n:=\sqrt{\Var(R_2)}$ and $b_n:=\Ex[R_2]$. Then we have
\begin{align*}
|b_n|& = \Big| \Ex\Big[ \sum_{i=1}^{n/2} \Ex\big[ \hat g(x_{2i-1}) - g(x_{2i-1}) + \frac{1}{\gamma} \mathbf{1}_{\{Y_{2i-1} \ge \hat g(x_{2i-1})\}} \mid \hat g \big] + \sum_{i=1}^{n/2} \Ex\big[\big(\tfrac{1}{\hat\gamma} - \tfrac{1}{\gamma}\big) \cdot  \mathbf{1}_{\{Y_{2i-1} \ge \hat g(x_{2i-1})\}} \mid \hat\gamma, \hat g \big] \Big]\Big| \\
&=\Big| \Ex\Big[ \sum_{i=1}^{n/2} \hat g(x_{2i-1}) - g(x_{2i-1}) + \frac{1}{\gamma} (1- F( \hat g(x_{2i-1})-g(x_{2i-1})))  + \sum_{i=1}^{n/2} \big(\tfrac{1}{\hat\gamma} - \tfrac{1}{\gamma}\big) \cdot  (1-F(\hat g(x_{2i-1})-g(x_{2i-1}))) \Big]\Big| \\
&\le \sum_{i=1}^{n/2} \Ex\Big[\frac{C_F}{\gamma} (\hat g(x_{2i-1}) - g(x_{2i-1}))^2+|\tfrac{1}{\hat\gamma}-\tfrac{1}{\gamma}| \cdot (\gamma|\hat g(x_{2i-1}) - g(x_{2i-1})| + C_F(\hat g(x_{2i-1})- g(x_{2i-1}))^2 )\Big]\\
&\in O(n(nh)^{-2}+n\sqrt{(nh)^{-1}}\sqrt{(nh)^{-2}}) \subseteq o(h^{-1}) 
\end{align*}
due to the Cauchy-Schwarz inequality and, moreover, a similar calculation yields
\begin{align*} 
a_n^2 &= \Var\big( \Ex[R_2\mid \hat\gamma, \hat g] \big) + \Ex\big[ \Var(R_2\mid \hat\gamma, \hat g) \big]\\
 &\le \begin{aligned}[t] &\Var\Big( \sum_{i=1}^{n/2} \hat g(x_{2i-1}) - g(x_{2i-1}) + \frac{1}{\gamma} (1- F( \hat g(x_{2i-1})-g(x_{2i-1})))  + \sum_{i=1}^{n/2} \big(\tfrac{1}{\hat\gamma} - \tfrac{1}{\gamma}\big) \cdot  (1-F(\hat g(x_{2i-1})-g(x_{2i-1})))\Big) \\&+ \sum_{i=1}^{n/2} \Ex\Big[ \frac{1}{\gamma^2} \cdot (\gamma |\hat g(x_{2i-1}) - g(x_{2i-1})| + C_F(\hat g(x_{2i-1})- g(x_{2i-1}))^2) \Big] \\  &+ \sum_{i=1}^{n/2}\Ex\Big[\big(\tfrac{1}{\hat\gamma}-\tfrac{1}{\gamma}\big)^2 \cdot\Var\big(  \mathbf{1}_{\{\hat g(x_{2i-1}) \ge Y_{2i-1}\}} \mid \hat\gamma, \hat g \big) \Big] \end{aligned}\\
&\in O(n^2h((nh)^{-4}+ \sqrt{(nh)^{-2}}\cdot\sqrt{(nh)^{-4}})+ n (nh)^{-1}+n\sqrt{(nh)^{-2}}\cdot\sqrt{(nh)^{-2}})= O(h^{-1}). \end{align*}
We get
\begin{align*}
 \mathbb{P}\Big( \big| \frac{R_2^2-b_n^2}{h^{-3/2}} \big|\ge \delta\Big) &= \mathbb{P}\Big( \big| \frac{(R_2-b_n)^2}{h^{-3/2}} + \frac{2(R_2-b_n)b_n}{h^{-3/2}}\big|\ge \delta\Big)\\
&\le \mathbb{P}\Big( \big| \frac{(R_2-b_n)^2}{h^{-3/2}} \big|\ge \delta/2\Big) + \mathbb{P}\Big( \big| \frac{2(R_2-b_n)b_n}{h^{-3/2}} \big|\ge \delta/2\Big)\\
&\le \frac{\Var\big(\frac{R_2-b_n}{h^{-1/2}}\big)}{\frac{\delta}{2h^{1/2}}}   +  \frac{\Var\big(\frac{R_2-b_n}{h^{-1/2}}\big)}{\frac{\delta^2}{16(hb_n)^2}}\\
&\longrightarrow 0,
\end{align*}
as $n\to\infty$. Hence $(S_22-b_n^2/n)/(n^{-1}h^{-3/2})\to 0$ in probability as $n\to\infty$. 

The second assertion is obtained the same way as the first one. 
\qed\medskip

For the proof of Theorem \ref{T:gCLT} we need a modified version of \cite[Theorem 3.1]{JMR14}. 

\begin{prop}\label{p:JMR3.1} Let $0\le a < b \le 1$. Assume that $g$ is affine. Then
\[ |\hat g(x) - g(x) | \le \max\{|Z_j(h)| \mid j=0,\dots, 2\lceil \tfrac{b-a}{h}\rceil +1\} \]
for all $x\in (a,b]$, where
\[ Z_j(h) := \max\{ \epsilon_{2i} \mid x_{2i} \in (a+(j-1)h/2, a+jh/2] \}. \]
\end{prop}
The main difference to \cite[Theorem 3.1]{JMR14} is that our bound is uniform in $x$. Besides, we only treat the special case that the function $g$ is affine and that the parameter $\beta^*$ of \cite[Theorem 3.1]{JMR14} takes the value $1$, thereby deriving explicit constants with a much shorter proof. 

\prf Let $x\in(a,b]$. There are $j_1, j_2\in\{0,\dots, 2\lceil \frac{b-a}{h}\rceil +1\}$ with $(a+(j_1-1)h/2,a+j_1h/2] \subseteq (x-h,x]$ and $(a+(j_2-1)h/2, a+j_2h/2]\subseteq (x,x+h]$. Since 
\[ \hat g(x) \ge \frac{x_{2i_2}-x}{x_{2i_2}-x_{2i_1}} (g(x_{2i_1})+\epsilon_{2i_1}) +  \frac{x-x_{2i_1}}{x_{2i_2}-x_{2i_1}} (g(x_{2i_2})+\epsilon_{2i_2})\]
for all $i_1,i_2$ with $x_{2i_1}\in(x-h,x]$ and $x_{2i_2}\in (x,x+h]$, we get 
\[ \hat g(x)-g(x) \ge  \frac{x_{2i_2}-x}{x_{2i_2}-x_{2i_1}} \epsilon_{2i_1} + \frac{x-x_{2i_1}}{x_{2i_2}-x_{2i_1}} \epsilon_{2i_2} \ge \min\{\epsilon_{2i_1},\epsilon_{2i_2}\}\]
using the assumption that $g$ is affine. Hence
\[ \hat g(x)-g(x) \ge \min\{Z_{j_1}(h), Z_{j_2}(h)\} \]
and, in particular, the assertion follows. \qed\medskip

\prf[ of Theorem \ref{T:gCLT}]
As mentioned in the beginning of this section the strategy is to apply \cite[Theorem 3]{Ph69}. To do this, we put $S(\lambda):= \lambda_1 S_1'+ \lambda_2S_1''$ for $\lambda=(\lambda_1,\lambda_2)$ with $\lambda_1\in(0,\infty)$, $\lambda_2\in(-\infty,0)$. We decompose the sum $S(\lambda)$ in a sequence of finitely many random variables. For every $n\in\mathbb{N}$ we put
\begin{align*}
& K(n):=\lceil \frac{1}{h} \rceil, \\
& \mathcal{X}'_{n,k} = \sum_{i:x_{2i-1}\in ((k-1)h,kh]} \lambda_1\big(\hat g(x_{2i-1})-g(x_{2i-1}) \big)^2 + \lambda_2 (Y_{2i-1}-g(x_{2i-1}))\mathbf{1}_{\{Y_{2i-1}\ge \hat g(x_{2i-1}) \}},\ k=1, \dots, K(n)-1,  \\
& \mathcal{X}'_{n,K(n)} = \sum_{i: x_{2i-1}\in ((K(n)-1)\cdot h, 1]} \lambda_1 \big(\hat g(x_{2i-1})-g(x_{2i-1}) \big)^2 + \lambda_2 (Y_{2i-1}-g(x_{2i-1}))\mathbf{1}_{\{Y_{2i-1}\ge \hat g(x_{2i-1}) \}},\\
& \mathcal{X}_{n,k}=\frac{\mathcal{X}'_{n,k} -\Ex[\mathcal{X}'_{n,k}]}{\sqrt{\Var(S(\lambda))}}, \, k=1,\dots, K(n). 
\end{align*}
 Now  
\[ \Ex[\mathcal{X}_{n,k}]=0, \, k=1,\dots, K(n),\qquad  \Ex\Big[ \big(\sum_{k=1}^{K(n)} \mathcal{X}_{n,k}\big)^2 \Big] =1,\]
and the sequence $\mathcal{X}_{n,k},\, k=1,\dots, K(n),$ is $m$-dependent with $m=2$ and in particular $\varphi$-mixing with $\varphi(\mathfrak{n})=0$ for $\mathfrak{n}\ge 3$. Assume that $\lambda$ is such that $\Var(S(\lambda)) \asymp n^{-2}h^{-3}$. By Lemma 5 this is the case for all $\lambda\in\mathbb{R}^2$ except for one line.  From Lemma \ref{l:variances} we further conclude $\max_k \Var(\mathcal{X}'_{n,k})\in O( (hn)^{-2})$, which yields $\max_k\Var(\mathcal{X}_{n,k})\in O(h)$. 

Most of the assumptions of \cite[Theorem 3]{Ph69} follow directly from these properties. For the remaining assumptions we have to build blocks
\[ \mathcal{Y}_j= \mathcal{X}_{n,\rho_j+1} + \dots + \mathcal{X}_{n,\rho_{j+1}}, \quad j=1,\dots,l, \]
where $0=\rho_1<\rho_2<\dots<\rho_{l+1}\le K(n)$, and a (perhaps empty) remainder block
\[ \mathcal{Y}_{l+1}= \mathcal{X}_{n,\rho_{l+1}+1}+\dots+ \mathcal{X}_{n,K(n)}.\]
These blocks have to meet the following assumptions: 
\begin{enumerate}[{(A}1{)}]
\item\label{A1} There is a sequence $V(n)$ with
\[ \lim_{n\to\infty} V(n)=0, \quad \lim_{n\to\infty} \frac{V(n)}{\Ex[\mathcal{X}_{n,1}^2]} =\infty \]
such that
\[ \max_{j=1,\dots,l} \Big| \frac{\Ex[\mathcal{Y}_j^2]}{V(n)}- 1 \Big| \longrightarrow 0,\quad n\to\infty,\]
and
\[ \limsup_{n\to \infty} \frac{\Ex[\mathcal{Y}_{l+1}^2]}{V(n)} \le 1. \]
\item\label{A2} For every $\epsilon>0$
\[ \lim_{n\to\infty}\sum_{j=1}^l \Ex\big[\mathcal{Y}^2_{j}\mathbf{1}_{\{|\mathcal{Y}_j|>\epsilon\}}\big]=0.\]
\end{enumerate}
We choose $V(n):=h^{2/3}$ and define inductively $\rho_1:=0$,  
\[\rho_{j+1} := \min\Big\{\rho\in\mathbb{N} \mid \Ex\Big[ \big(\sum_{k=\rho_j+1}^{\rho} \mathcal{X}_{n,k}\big)^2 \Big]> V(n) \Big\}, \quad j \ge 1, \]
let $l\in\mathbb{N}$ denote the largest index for which this definition makes sense (i.e.\ the minimum is taken over a non-empty set). Now assumption (A\ref{A1}) is easily seen to be fulfilled, since
\begin{align*}
\Ex[\mathcal{Y}_j^2]= \Ex\Big[ \big(\sum_{k=\rho_j+1}^{\rho_{j+1}} \mathcal{X}_{n,k}\big)^2 \Big] &= \Ex\Big[ \big(\sum_{k=\rho_j+1}^{\rho_{j+1}-1} \mathcal{X}_{n,k}\big)^2 \Big] + 2\Ex\Big[ \sum_{k=\rho_j+1}^{\rho_{j+1}-1} \mathcal{X}_{n,k}\cdot \mathcal{X}_{n,\rho_{j+1}} \Big] + \Ex\Big[ \mathcal{X}_{n,\rho_{j+1}}^2 \Big]\\
& \le V(n) + 2 \sqrt{V(n)} \sqrt{\Ex\Big[ \mathcal{X}_{n,\rho_{j+1}}^2 \Big]} + \Ex\Big[ \mathcal{X}_{n,\rho_{j+1}}^2 \Big] \end{align*}
and
\[ \lim_{n\to\infty} \frac{\Ex[\mathcal{X}_{n,\rho_{j+1}}^2]}{V(n)}=0.\]

As a preparatory step for checking (A\ref{A2}) we show
\begin{equation} \rho_{j+1}-\rho_j \le \frac{h^{2/3}}{\min\{\Var(\mathcal{X}_{n,k}) \mid k=1, \dots, K(n)-1\}}+1 \in O(h^{-1/3}).\label{e:rho_diff}\end{equation}
Indeed,
\[
\Cov(\mathcal{X}_{n,k},\mathcal{X}_{n,k'}) = \Ex\big[ \Cov(\mathcal{X}_{n,k}, \mathcal{X}_{n,k'}\mid \hat g)\big] + \Cov\big(\Ex[ \mathcal{X}_{n,k} \mid \hat g], \Ex[ \mathcal{X}_{n,k'}\mid\hat g] \big)
\ge 0
\]
since $\Cov(\mathcal{X}_{n,k},\mathcal{X}_{n,k'}\mid \hat g)=0$ a.s.\ for $k\ne k'$ and the random variables $\Ex[\mathcal{X}_{n,k} \mid \hat g]$ and $\Ex[\mathcal{X}_{n,k'} \mid \hat g]$ are component-wise monotonically decreasing function of $\hat g(x_{2i-1})$, $i=1, \dots, n/2$ and the latter are, in turn, component-wise monotonically increasing functions of $Y_{2i}$, $i=1, \dots, n/2$. Hence
\[ V(n) \ge  \Var\big( \sum_{k=\rho_j+1}^{\rho_{j+1}-1} \mathcal{X}_{n,k}\big) \ge   (\rho_{j+1}-1-\rho_j) \cdot \min\{\Var(\mathcal{X}_{n,k}) \mid k=1, \dots, K(n)-1\}. \]
Now Lemma \ref{l:variances} yields \eqref{e:rho_diff}.

In order to check assumption (A\ref{A2}) observe that
\begin{align*}
 \mathcal{Y}_j &=  \begin{aligned}[t]  \lambda_1\frac{1}{\sqrt{\Var\big(S(\lambda) \big)}}&\Big(\sum_{i:x_{2i-1}\in [\rho_j h,\rho_{j+1}h)} \big(\hat g(x_{2i-1})-g(x_{2i-1}) \big)^2 \Big) \\
& + \lambda_2 \frac{1}{\sqrt{\Var\big(S(\lambda) \big)}}\Big(\sum_{i: x_{2i-1}\in[\rho_j\cdot h, \rho_{j+1}\cdot h)} (Y_{2i-1}-g(x_{2i-1})) \mathbf{1}_{\{Y_{2i-1} \ge \hat g(x_{2i-1})\}} \Big)\\
& - \frac{1}{\sqrt{\Var\big(S(\lambda) \big)}}\Ex\Big[\sum_{i: x_{2i-1} \in [\rho_j\cdot h, \rho_{j+1}\cdot h)} \lambda_1 \big(\hat g(x_{2i-1})-g(x_{2i-1}) \big)^2 +\lambda_2 (Y_{2i-1}-g(x_{2i-1})) \mathbf{1}_{\{Y_{2i-1} \ge \hat g(x_{2i-1})\}}\Big] \end{aligned} \\
&=: \lambda_1\mathcal{Y}_j^{(1)} + \lambda_2\mathcal{Y}_j^{(2)} - \mathcal{Y}_j^{(3)}. 
\end{align*}
Let $\epsilon>0$. Then
\begin{align*}
\Ex\big[ \mathcal{Y}_j^2 \mathbf{1}_{\{|\mathcal{Y}_j|> \epsilon\}} \big]  &\le 9 \Ex\big[ \max\big\{ (\lambda_1\mathcal{Y}_j^{(1)})^2,  (\lambda_2\mathcal{Y}_j^{(2)})^2, (\mathcal{Y}_j^{(3)})^2 \big\} \mathbf{1}_{\{ 3\max\{ |\lambda_1\mathcal{Y}_j^{(1)}|, |\lambda_2\mathcal{Y}_j^{(2)}|, |\mathcal{Y}_j^{(3)}|\} > \epsilon \}} \big] \\
&\le 9 \lambda_1^2\Ex[(\mathcal{Y}_j^{(1)})^2\mathbf{1}_{\{|\lambda_1\mathcal{Y}_j^{(1)}|>\epsilon/3\}}] + 9 \lambda_2^2 \Ex[(\mathcal{Y}_j^{(2)})^2\mathbf{1}_{\{|\lambda_2\mathcal{Y}_j^{(2)}|>\epsilon/3\}}] + 9 (\mathcal{Y}_j^{(3)})^2\mathbf{1}_{\{|\mathcal{Y}_j^{(3)}|>\epsilon/3\}}. 
\end{align*}
Put $\nu^-(n):=\inf\{ \#\{i\in\{1, \dots, n\}\mid x_i \in [x,x+h/2)\} \mid x\in[-h/2,1] \}$ and $\nu^+(n):=\sup\{ \#\{i\in\{1, \dots, n\}\mid x_i \in [x,x+h/2)\} \mid x\in[-h/2,1] \}$. Since
\begin{align*}
\Ex\Big[&\sum_{i: x_{2i-1}\in[ \rho_j \cdot h, \rho_{j+1}\cdot h)} \lambda_1\big(\hat g(x_{2i-1})-g(x_{2i-1}) \big)^2 +  \lambda_2 (Y_{2i-1}-g(x_{2i-1})) \mathbf{1}_{\{Y_{2i-1} \ge \hat g(x_{2i-1})\}} \Big] \\
& \le (\rho_{j+1}-\rho_j) \cdot \nu^+(n) \cdot \max\{ \mathbb{E}\big[ \lambda_1 (\hat g(x_{2i-1})- g(x_{2i-1}))^2 + \lambda_2 (Y_{2i-1}-g(x_{2i-1})) \mathbf{1}_{\{Y_{2i-1} \ge \hat g(x_{2i-1})\}}\big] \mid i=1,\dots, n/2\big\}\\
&\in O(h^{-1/3}hn(hn)^{-2})\\
&= O(h^{-4/3}n^{-1})
\end{align*}
by Proposition \ref{p:asy_ex} and we assumed
\[ \Var(S(\lambda)) \asymp n^{-2}h^{-3}, \]
we have $|\mathcal{Y}_j^{(3)}|< \epsilon/3$ provided that $n$ is sufficiently large which will be assumed in the sequel.

By Proposition \ref{p:JMR3.1} we have
\[ |\hat g(x)-g(x)| \le \max\big\{ \big|Z_\mathfrak{j}(h,j) \big| \mid \mathfrak{j} =0, \dots, 2 (\rho_{j+1} - \rho_j )+1\big\} \] 
for all $x\in [\rho_j h, \rho_{j+1} h ]$, where
\[ Z_\mathfrak{j}(h,j)=\max\big\{\epsilon_{2i}\mid x_{2i}\in \big(\rho_j  \cdot h+ (\mathfrak{j}-1) h/2, \rho_{j} \cdot h + \mathfrak{j}h/2 \big) \big\}. \]
Put $V:=\sqrt{\Var(S(\lambda))}$. We have for any $t>0$
\begin{align*}
 \mathbb{P}(\mathcal{Y}_j^{(1)} > t) &\le \mathbb{P}\Big(\sum_{i: x_{2i-1} \in [\rho_j \cdot h, \rho_{j+1} \cdot h]} (\hat g(x_{2i-1})-g(x_{2i-1}))^2 >tV \Big) \\
&\le \mathbb{P}\Big( (\rho_{j+1}-\rho_j) \cdot \nu^+(n) \cdot \max_{\mathfrak{j}}  \big|Z_\mathfrak{j} (h, j)\big|^2 > tV \Big)\\
& \le \sum_{\mathfrak{j}=0}^{2(\rho_{j+1} - \rho_j)+1 } \mathbb{P}\Big( \max_i \epsilon_i< - \sqrt{\frac{tV}{ (\rho_{j+1}-\rho_j)\cdot \nu^+(n)}} \Big)\\
&\le \big(2\rho_{j+1} - 2\rho_j +2\big) \cdot F\Big(- \sqrt{\frac{tV}{ (\rho_{j+1}-\rho_j)\cdot \nu^+(n)}} \Big)^{\lfloor \nu^-(n)/2 \rfloor}. 
\end{align*}

Thus, with $\epsilon':=\epsilon/(3\lambda_1)$, 
\[ \Ex\big[ (\mathcal{Y}_j^{(1)})^2\mathbf{1}_{\{ |\mathcal{Y}_j^{(1)}|>\epsilon' \}} \big] = \int_{\epsilon'}^\infty \mathbb{P}(\mathcal{Y}_j^{(1)}>t)2t \, dt \le \int_{\epsilon'}^\infty (2\rho_{j+1} - {2\rho_j }+2)   \cdot F\Big(- \sqrt{\frac{tV}{ (\rho_{j+1}-\rho_j)\cdot \nu^+(n)}} \Big)^{\lfloor \nu^-(n)/2 \rfloor} 2t \, dt. \] 
By \eqref{e:lim=gamma} there is some $s\in(-\infty, 0)$ with
\[ F(u) \le \exp\big\{ -\tfrac{\gamma}{2} |u|\big\}, \quad u\in(s,0).\]
Abbreviate $W(n):=  V/((\rho_{j+1}-\rho_j)\nu^+(n)) $. Then
\begin{align*} \int_{\epsilon'}^\infty  &F\big(- \sqrt{tW(n)} \big)^{\lfloor \nu^-(n)/2 \rfloor} 2t \, dt\\
& \le  \int_{\epsilon'}^{s^2/W(n) \vee \epsilon'} \exp\big\{- \tfrac{\gamma}{2}\sqrt{tW(n)} \big\}^{\lfloor \nu^-(n)/2 \rfloor} 2t \, dt +  \int_{s^2/W(n)}^\infty  F\big(- \sqrt{tW(n)} \big)^{\lfloor \nu^-(n)/2 \rfloor} 2t \, dt.
\end{align*}
While the second summand may be infinite for small $n$, it is finite for large enough $n$, since we assume \eqref{e:ass_distr_tail}. So choose some $n_0$ for which this summand is finite. Then for all $n$ with $\nu^-(n) \ge \nu^-(n_0)$ we have 
\begin{align*} 
\int_{s^2/W(n)}^\infty & F\big(- \sqrt{tW(n)} \big)^{\lfloor \nu^-(n)/2 \rfloor} 2t \, dt \\
&\le F(s)^{\lfloor \nu^-(n)/2 \rfloor - \lfloor \nu^-(n_0)/2 \rfloor} \int_{s^2/W(n_0)}^\infty  F\big(- \sqrt{uW(n_0)} \big)^{\lfloor \nu^-(n_0)/2 \rfloor} 2u\cdot  \big(\frac{W(n_0)}{W(n)}\big)^2\, du.  
\end{align*}
So this summand converges to zero with exponential speed. 

For the first summand we get
\begin{align*}
\int_{\epsilon'}^{s^2/W(n)\vee \epsilon'} \exp\big\{- \tfrac{\gamma}{2}\sqrt{tW(n)} \big\}^{\lfloor \nu^-(n)/2 \rfloor} 2t \, dt
&\le \int_{\epsilon'}^\infty \exp\big\{-\tfrac{\gamma}{2} {\lfloor \nu^-(n)/2\rfloor}\cdot \big( tW(n)\big)^{1/2} \big\} 2t \, dt \\
& = \int_{\epsilon'\cdot \widetilde W}^\infty \exp\big\{-\tfrac{\gamma}{2} u^{1/2}\} \frac{2u }{\widetilde W^2 } \, du,
\end{align*}
where we have put $\widetilde W:=W(n) \cdot \lfloor \nu^-(n)/2 \rfloor^{2}$. 

We see
\[ \widetilde W \asymp \frac{(hn)^{2}}{\sqrt{n^2h^3} \cdot h^{-1/3}\cdot  nh} = h^{-1/6}. \]
For sufficiently large $u$ we have $\exp\{-\gamma/2 \cdot u^{1/2}\}u\le \tfrac{1}{5u^6}$ and thus
\[ \int_{\epsilon' \cdot \widetilde W}^\infty  \exp\big\{-\tfrac{\gamma}{2} u^{1/2}\}  u \, du \le \frac{1}{(\epsilon' \cdot \widetilde W)^5}. \]
So we have (up to an exponentially decaying term)
\begin{equation}
 \Ex\big[ (\mathcal{Y}_j^{(1)})^2 \mathbf{1}_{\{|\mathcal{Y}_j^{(1)}|>\epsilon'\}} \big] \le \big(2\rho_{j+1} - 2\rho_j +2 \big) \cdot \frac{1}{(\epsilon \cdot \tilde W /(3\lambda_1))^5} \cdot \frac{1}{\tilde W^2} \asymp h^{5/6}. \label{e:Y1_order} \end{equation}

We now turn to $\mathcal{Y}_j^{(2)}$. We have
\begin{align*} 
\mathbb{P}\big(|\mathcal{Y}_j^{(2)}|>t\big) &= \mathbb{P}\Big( \sum_{i:x_{2i-1}\in[\rho_j\cdot h,\rho_{j+1}\cdot h)} |Y_{2i-1}-g(x_{2i-1})| \mathbf{1}_{\{Y_{2i-1} \ge \hat g(x_{2i-1})\}} > tV \Big)\\
&\le \mathbb{P}\Big( \sum_{i:x_{2i-1}\in[\rho_j\cdot h,\rho_{j+1}\cdot h)} \max_{\mathfrak{j}} Z_\mathfrak{j}(h,j) \mathbf{1}_{\{|Y_{2i-1}- g(x_{2i-1})| \le \max_{\mathfrak{j}} Z_\mathfrak{j}(h,j)\}} > tV \Big)
\\
& \le \begin{aligned}[t] \mathbb{P}&\Big( \max_{\mathfrak{j}} Z_\mathfrak{j}(h,j) >\frac{1}{2}\sqrt{\frac{tV}{\bar\gamma\cdot (\rho_{j+1}-\rho_j) \cdot \nu^+(n)}} \Big) \\
 +  &  \mathbb{P}\Big( \Big\{ \sum_{i: x_{2i-1}\in[\rho_j\cdot h, \rho_{j+1}\cdot h)}  \max_{\mathfrak{j}} Z_\mathfrak{j}(h,j) \mathbf{1}_{\{|Y_{2i-1}- g(x_{2i-1})| \le \max_{\mathfrak{j}} Z_\mathfrak{j}(h,j)\}} > tV \Big\}\\
& \cap \Big\{\max_{\mathfrak{j}} Z_\mathfrak{j}(h,j) \le \frac{1}{2}\sqrt{\frac{tV}{\bar\gamma\cdot (\rho_{j+1}-\rho_j) \cdot \nu^+(n)}} \Big\} \Big),\end{aligned}
\end{align*}
where $\bar\gamma$ is a constant such that $1-F(t) \le \bar\gamma|t|$ for all $t\le 0$. Put $\epsilon'':=\epsilon/(3 \lambda_2)$. Then
\begin{align*}
\Ex[(\mathcal{Y}_j^{(2)})^2\mathbf{1}_{\{|\mathcal{Y}_j^{(2)}|>\epsilon''\}}]
 &\le \begin{aligned}[t] \int_{\epsilon''}^\infty  \mathbb{P}\Big( &\max_{\mathfrak{j}} Z_\mathfrak{j}(h,j) > \frac{1}{2} \sqrt{\frac{tV}{\bar\gamma\cdot  (\rho_{j+1}-\rho_j) \cdot \nu^+(n)}} \Big) \cdot 2t \, dt\\
+\int_{\epsilon''}^\infty \mathbb{P}\Big( &\Big\{ \sum_{i:x_{2i-1}\in[\rho_j\cdot h, \rho_{j+1}\cdot h)} \max_{\mathfrak{j}} Z_\mathfrak{j}(h,j) \mathbf{1}_{\{|Y_{2i-1}- g(x_{2i-1})| \le \max_{\mathfrak{j}} Z_\mathfrak{j}(h,j)\}} > tV \Big\}\\
& \cap \Big\{\max_{\mathfrak{j}} Z_\mathfrak{j}(h,j) \le \frac{1}{2} \sqrt{\frac{tV}{\bar\gamma\cdot (\rho_{j+1}-\rho_j) \cdot \nu^+(n)}} \Big\} \Big) \cdot 2t \, dt \end{aligned}
\\
& =: A_I+A_{II}.
\end{align*}

We get $A_I \asymp h^{5/6}$ with exactly the same arguments that were used to derive \eqref{e:Y1_order}.

So let us turn to $A_{II}$. Using Lemma \ref{l:binom}, \ref{l:poisson} and \ref{l:stirling} we can estimate 
\begin{align*}
\mathbb{P}\bigg( \bigg\{ &\begin{aligned}[t]  \sum_{i: x_{2i-1}\in[\rho_j\cdot h, \rho_{j+1}\cdot h)}& \max_{\mathfrak{j}} Z_\mathfrak{j}(h,j) \mathbf{1}_{\{|Y_{2i-1}- g(x_{2i-1})| \le \max_{\mathfrak{j}} Z_\mathfrak{j}(h,j)\}} >tV \bigg\} \\& \cap \bigg\{ \max_{\mathfrak{j}} Z_\mathfrak{j}(h,j) \le \frac{1}{2} \sqrt{\frac{tV}{\bar\gamma\cdot (\rho_{j+1}-\rho_j) \cdot \nu^+(n)}} \bigg\} \bigg)\end{aligned}\\
& \le \mathbb{P}\Big(  \sum_{i: x_{2i-1} \in[\rho_j\cdot h, \rho_{j+1}\cdot h)} \mathbf{1}_{\big\{|Y_{2i-1}- g(x_{2i-1})| \le \frac{1}{2\bar\gamma}\sqrt{\frac{\bar\gamma tV}{(\rho_{j+1}-\rho_j) \cdot \nu^+(n)}}\big\}} >2\sqrt{\bar\gamma tV (\rho_{j+1}-\rho_j) \cdot \nu^+(n)} \Big)\\
& \le 2\binom{m}{k} p^k (1-p)^{m-k} \le 2\frac{(mp)^k}{k!} \le 2\frac{(mp)^k}{\sqrt{\frac{\pi k}{2}} \big(\frac{k}{e}\big)^{k}} = \sqrt{\frac{8}{\pi k}} \big(\frac{mpe}{k}\big)^k\le \sqrt{\tfrac{8}{\pi k}} \big(\tfrac{e}{4}\big)^{k},
\end{align*}
where
\[
k:= \big\lceil 2\sqrt{\bar\gamma tV (\rho_{j+1}-\rho_j) \cdot \nu^+(n)} \big\rceil, \quad
m:= (\rho_{j+1}-\rho_j) \cdot \nu^+(n) \quad \mbox{and} \quad
p:= \frac{1}{2} \sqrt{\frac{\bar\gamma t V}{(\rho_{j+1}-\rho_j) \cdot \nu^+(n)}}, 
\]
since the indicators $\mathbf{1}_{\big\{|Y_{2i-1}- g(x_{2i-1})| \le \frac{1}{2\bar\gamma}\sqrt{\frac{\bar\gamma tV}{(\rho_{j+1}-\rho_j) \cdot \nu^+(n)}}\big\}}$ are independent with success probability 
\[1- F\Big( -\frac{1}{2\bar\gamma}\sqrt{\frac{\bar\gamma tV}{(\rho_{j+1}-\rho_j) \cdot \nu^+(n)}} \Big) \le \frac{1}{2}\sqrt{\frac{\bar\gamma tV}{(\rho_{j+1}-\rho_j) \cdot \nu^+(n)}}. \]
Hence 
\begin{align*}
A_{II} \le \int_{\epsilon''}^\infty \frac{2}{\sqrt{\pi\sqrt{tW[n]} }} \big(\frac{e}{4}\big)^{2\sqrt{tW[n]} } 2t \, dt &= \int_{\epsilon'' W[n]}^\infty \frac{2}{\sqrt{\pi\sqrt{u} }} \big(\frac{e}{4}\big)^{2\sqrt{u} } 2\frac{u}{W[n]^2} \, du\\
&\le \frac{1}{W[n]^2} \int_0^\infty \frac{2}{\sqrt{\pi\sqrt{u} }} \big(\frac{e}{4}\big)^{2\sqrt{u}} 2u\, du,
\end{align*}
where we have put
\[ W[n]:=  \bar\gamma V (\rho_{j+1}-\rho_j) \cdot \nu^+(n)\asymp \sqrt{n^{-2}h^{-3}}h^{-1/3}nh= h^{-5/6}. \]
Hence we get
\[ \Ex[(\mathcal{Y}_j^{(2)})^2\mathbf{1}_{\{|\mathcal{Y}_j^{(2)}|>\epsilon''\}}] \asymp h^{5/6}+  (h^{-5/6})^{-2} = h^{5/6}. \]
In conclusion, 
\begin{align*}
\sum_{j=1}^l \Ex\big[\mathcal{Y}_j^2\mathbf{1}_{\{|\mathcal{Y}_j|>\epsilon\}}\big] &\le \sum_{j=1}^l 9\lambda_1^2\Ex\big[(\mathcal{Y}_j^{(1)})^2\mathbf{1}_{\{|\mathcal{Y}_j^{(1)}|>\epsilon/(3\lambda_1)\}}\big] + 9\lambda_2^2\Ex\big[(\mathcal{Y}_j^{(2)})^2\mathbf{1}_{\{|\mathcal{Y}_j^{(2)}|>\epsilon/(3\lambda_2)\}}\big] \asymp h^{-2/3} \cdot h^{5/6} = h^{1/6}
\end{align*}
which shows that assumption (A\ref{A2}) is fulfilled. 

Thus we can apply \cite[Theorem 3]{Ph69}, which yields
\[ \frac{S(\lambda) - \Ex\big[S(\lambda) \big]}{\sqrt{\Var(S(\lambda))}} =\sum_{k=1}^{K(n)} \mathcal{X}_{n,k} \longrightarrow \mathcal{N}(0,1) \] 
in distribution as $n\to\infty$. Hence
\[ \frac{S(\lambda) - \Ex\big[S(\lambda) \big]}{\sqrt{n^{-2}h^{-3}}} \longrightarrow \mathcal{N}(0,\sigma^2), \]
since $\sigma^2:=\lim_{n\to\infty} \Var(S(\lambda))/n^{-2}h^{-3}$ is assumed to exist. 
Now a sharp version of the Cram\'er-Wold theorem \cite{CFR07} yields the assertion.\qed \medskip

\begin{lemma}\label{l:binom}
Let $X\sim\mathsf{Bin}(m,p)$ be a binomial distributed random variable with parameters $m\in\mathbb{N}$ and $p\in(0,1)$ and let $k\in\{0,\dots, m\}$. If $k\ge 4mp$, then
\[\mathbb{P}(X \ge k) \le 2 \mathbb{P}(X=k). \]
\end{lemma}
\prf We have
\[ \mathbb{P}(X\ge k) = \sum_{l=k}^m \binom{m}{l}p^l(1-p)^{m-l}. \]
Since for $l\ge k$ we have
\[\frac{\binom{m}{l+1}p^{l+1}(1-p)^{m-l-1}}{\binom{m}{l}p^l(1-p)^{m-l}} = \frac{m-l}{l+1} \frac{p}{1-p} \le \frac{mp}{k \cdot \tfrac{1}{2}} \le \frac{1}{2}, \]
we get
\[ \mathbb{P}(X \ge k) = \sum_{l=k}^m \binom{m}{l}p^l(1-p)^{m-l} \le \sum_{l=k}^m \big(\frac{1}{2}\big)^{l-k}\binom{m}{k}p^k(1-p)^{m-k}\le 2\mathbb{P}(X=k). \qed\]
 
\begin{lemma}\label{l:poisson}
Let $X\sim\mathsf{Bin}(m,p)$ be a binomial distributed random variable with parameters $m\in\mathbb{N}$ and $p\in(0,1)$ and let $k\in\{0,\dots, m\}$. Then
\[ \mathbb{P}(X=k) = \binom{m}{k} p^k (1-p)^{m-k} \le \frac{(mp)^k}{k!}. \]
\end{lemma}
\prf We have
\[\binom{m}{k} p^k (1-p)^{m-k} = \frac{\prod_{i=m-k+1}^m i}{k!} p^k (1-p)^{m-k} \le \frac{(mp)^k}{k!}. \qquad \qed\]

\begin{lemma}\label{l:stirling}
For $k\in\mathbb{N}$ we have
\[k! \ge \sqrt{\tfrac{\pi k}{2}}\big(\tfrac{k}{e}\big)^k. \]
\end{lemma}
\prf The proof uses some ideas of the proof of Stirling's formula $k! \sim \sqrt{2\pi k}(k/e)^k$. 

We have 
\begin{align*}
k!= \int_0^\infty x^k e^{-x} \, dx 
= \int_0^\infty e^{k \ln x - x} \, dx 
= \int_0^\infty e^{k \ln (ky) -ky} k \, dy 
&= k \cdot e^{k \ln k} \int_0^\infty  e^{k(\ln y-y)} \, dy\\
&= k^{k+1} \int_0^\infty  e^{k(\ln y-y)} \, dy. 
\end{align*}
We are now going to find a lower bound for this integral using a Taylor expansion for
\[ f(y):= \ln y-y. \]
We have
\[ f'(y)= \frac{1}{y}-1, \qquad f''(y) = \frac{-1}{y^2}, \qquad f'''(y)= \frac{2}{y^3} , \qquad y>0,\]
and hence
\[ f(y) = f(1)+ f'(1)(y-1)+\tfrac{1}{2} f''(1)(y-1)^2+ \tfrac{1}{6} f'''(\tilde y)(y-1)^3 \ge -1- \tfrac{1}{2}(y-1)^2\]
for all $y\in[1,\infty)$ with some $\tilde y$ depending on $y$. 
Therefore
\[\int_0^\infty e^{k(\ln y-y)}\, dy \ge \int_1^\infty e^{k(-1-\frac{1}{2}(y-1)^2)}\, dy = e^{-k} \int_0^\infty e^{-\frac{1}{2}x^2} \frac{1}{\sqrt{k}}\, dx = e^{-k} \sqrt{\tfrac{\pi}{2k}}. \]
So
\[ k! \ge k^{k+1} e^{-k} \sqrt{\tfrac{\pi}{2k}} = \big(\tfrac{k}{e}\big)^k \sqrt{\tfrac{\pi k}{2}}. \qquad \qed \]

\subsection{An approximation of the variance involving Poisson processes}\label{ss:Poisson}

Here we find an asymptotically equivalent expression for the variance of the test statistic $T$ that only depends on $\gamma$ and not on the entire distribution of the errors. We use the fact that only observations close to the border of the support play a role in computing the test statistic $T$ and these observations asymptotically form a Poisson process.  

In order to make this asymptotic behavior precise, we use the vague topology on the set of point measures. 
The set of all point measures on a closed set $E\subseteq\mathbb{R}^d$ will be denoted by $M_p(E)$. It is equipped with the vague topology. We say that a sequence $(\phi_n)_{n\in\mathbb{N}}$ in $M_p(E)$ converges vaguely to $\phi\in M_p(E)$ if 
\[\lim_{n\to\infty} \int_E f(x) \, d\phi_n(x) = \int_E f(x) \, d\phi(x) \]
for all continuous functions $f:E\to [0,\infty)$ with compact support. For further information on vague convergence, see e.g. \cite[Section 3.4]{Res87}.  We will rescale the system of observation points such that it converges in the vague topology on the set of point measures. At the same time the functional $G$, which is applied to the point measures in order to obtain the test statistic, should be a fixed functional independent of the sample size. Since in the computation of $G$ it is crucial whether two points of the data set have a distance in $x$-direction of more or less than $h$, we have to rescale in such a way that $h$ becomes fixed -- we will rescale by a factor $1/h$ so that $h$ becomes $1$. Thus the distance of two neighboring points in $x$-direction is $2(nh)^{-1}$. So we have to rescale by a factor $nh$ in $y$-direction in order to get a non-degenerate limit of the point processes. Hence we define point processes on $[0,1]\times (-\infty, 0]$ by
\[\Phi_{k,n}^o = \Big\{ \big(\frac{x_{2i-1}}{h}-(k-1), hn \cdot \epsilon_{2i-1}\big) \mid \frac{x_{2i-1}}{h} \in(k-1,k] \Big\}, \quad k=1, \dots, \big\lfloor \frac{1}{h} \big\rfloor, \]
and point processes on $[-1,2]\times (-\infty,0]$ by
\[ \Phi_{k,n}^e = \Big\{ \big(\frac{x_{2i}}{h}-(k-1), hn \cdot \epsilon_{2i}\big) \mid \frac{x_{2i}}{h} \in(k-2,k+1] \Big\}, \quad k=1, \dots, \big\lfloor \frac{1}{h} \big\rfloor \]
-- the ``e'' and ``o'' stand for ``even'' and ``odd'' (as already said, we will not be concerned with boundary effects in this paper but we will assume that some additional observations are available).   
We let $\Phi_k^o, \, k=1, \dots, \lfloor \frac{1}{h} \rfloor,$ be Poisson processes on $[0,1]\times (-\infty, 0]$ of intensity $1/(2\gamma)$ and we put
\[ \Phi_k^e := \big\{ (\mathfrak{X}_i-(k-1), \mathfrak{Y}_i) \mid (\mathfrak{X}_i,\mathfrak{Y}_i) \in \Phi^e,\ \mathfrak{X}_i\in(k-2,k+1] \big\}, \]
where $\Phi^e$ is a Poisson process of intensity $1/(2\gamma)$ on $\mathbb{R} \times (-\infty, 0]$. So the point processes $\Phi_k^e, \, k=1,\dots, \lfloor \frac{1}{h}\rfloor,$ are Poisson processes on $[-1,2] \times (-\infty,0]$ of intensity $1/(2\gamma)$ with a certain dependency in between them. 

Similar to \cite[Corollary 4.19(iii)]{Res87} we get the following proposition. 

\begin{prop}\label{p:PPcon} Let \eqref{e:lim=gamma} and (A\ref{i:ass_scatter_strong}) be fulfilled and let $(k_n)_{n\in\mathbb{N}}$ be a sequence with $k_n\in\{1, \dots, \lfloor \frac{1}{h}\rfloor \}$ for all $n\in\mathbb{N}$ and let $k\in\mathbb{N}$. Then we have
\[ (\Phi_{n,k_n}^o)_{k=1, \dots, K} \to (\Phi_k^o)_{k=1,\dots, K}, \qquad (\Phi_{n,k_n}^e)_{k=1,\dots, K} \to (\Phi_k^e)_{k=1,\dots, K}, \qquad n\to \infty\]
in distribution w.r.t.\ the vague topology on $M_p([0,1] \times (-\infty, 0])$ or $M_p([-1,2]\times (-\infty,0])$ respectively. \end{prop}

We are now going to find an asymptotically equivalent expression for the variance of $S_1$ based on the limiting Poisson processes. For this we decompose
\[ S_1= \sum_{k=1}^{\lfloor 1/h \rfloor} \frac{1}{nh} S_1(h,k) + \frac{1}{nh} \mathfrak{R}, \]
where
\begin{align*}
 S_1(h,k)&= nh \sum_{i:x_{2i-1}\in((k-1)h,kh]} (\hat g(x_{2i-1})-g(x_{2i-1}))^2 + nh \frac{2}{\hat \gamma} \sum_{i: x_{2i-1}\in((k-1)h, kh]} (Y_{2i-1}- g(x_{2i-1})) \mathbf{1}_{\{Y_{2i-1}\ge \hat g(x_{2i-1})\}} \\
&= S_1^1(h,k) + S_1^2(h,k)
\end{align*}
and 
\[ \mathfrak{R} = nh \sum_{i:x_{2i-1}\in(\lfloor 1/h\rfloor h,1]} (\hat g(x_{2i-1})-g(x_{2i-1}))^2 + nh \frac{2}{\hat\gamma} \sum_{i: x_{2i-1}\in (\lfloor 1/h\rfloor h,1]} (Y_{2i-1}- g(x_{2i-1})) \mathbf{1}_{\{Y_{2i-1}\ge \hat g(x_{2i-1})\}}\]
with
\[ S_1^1(h,k)= \frac{(nh)^2}{nh} \sum_{i:x_{2i-1}\in((k-1)h,kh]} (\hat g(x_{2i-1})-g(x_{2i-1}))^2 - \frac{(nh)^2}{2h}\int_{(k-1)h}^{kh} (\hat g(x)-g(x))^2 \, dx \]
and
\[ S_1^2(h,k)= \frac{(nh)^2}{2h} \int_{(k-1)h}^{kh} (\hat g(x)-g(x))^2 \, dx + nh\frac{2}{\hat\gamma} \sum_{i: x_{2i-1}\in((k-1)h,kh]} (Y_{2i-1}- g(x_{2i-1})) \mathbf{1}_{\{Y_{2i-1}\ge \hat g(x_{2i-1})\}}. \]
While $S_1^1(h,k)$ will be shown to converge to zero, $S_1^2(h,k)$ can be represented as a certain functional of the point processes $\Phi_{n,k}^e$ and $\Phi_{n,k}^o$ and thus converges to the functional applied to the limit processes $\Phi_k^e$ and $\Phi_k^o$.
 We recall 
\[  G(\phi^o, \phi^e, \gamma') := G_1(\phi^e)+\frac{2}{\gamma'}G_2(\phi^o,\phi^e) := \frac{1}{2}\int_0^1\tilde g(x)^2 \, dx + \frac{2}{\gamma'} \sum_{i=1}^\infty \mathfrak{Y}_i  \mathbf{1}_{\{ \mathfrak{Y}_i \ge \tilde g(\mathfrak{X}_i)\}} \]
where $\phi^o=\{ (\mathfrak{X}_i, \mathfrak{Y}_i) \mid i=1,\dots\}$ and $\tilde g$ is the estimator $\hat g$ at bandwidth $h=1$ applied to the data set $\phi^e$. If $\phi^e$ is of the form $\Phi_{k,n}^e$ then we have $\tilde g(x)=(\hat g(h \cdot (x+k-1))-g(h \cdot (x+k-1))\cdot nh$, $x\in[0,1]$, where the estimator $\hat g$ is applied with bandwidth $h$ to the observations $(x_{2i},Y_{2i})$. Thus we have in particular    
\begin{align*}  G(\Phi^o_{n,k}, \Phi^e_{n,k}, \hat\gamma) & = \frac{1}{2h} \int_{(k-1)h}^{kh} \tilde g\big(\tfrac{x}{h}-(k-1)\big)^2 \, dx + \frac{2}{\hat\gamma} \sum_{i:\frac{x_{2i-1}}{h} \in(k-1,k]} \epsilon_{2i-1}\cdot hn \cdot \mathbf{1}_{\{\epsilon_{2i-1}\cdot hn \ge (\hat g(x_{2i-1})-g(x_{2i-1})) \cdot hn\}} \\ 
&= \frac{(hn)^2}{2h} \int_{(k-1)h}^{kh} \big(g(x) - \hat g(x)\big)^2 \, dx+ \frac{2}{\hat\gamma}  hn \sum_{i:\frac{x_{2i-1}}{h} \in(k-1,k]} \big(Y_{2i-1}-g(x_{2i-1})\big) \mathbf{1}_{\{Y_{2i-1}\ge \hat g(x_{2i-1})\}}\\
&=S_1^2(h,k). \end{align*} 

We would like to show the following proposition. 

\begin{prop}\label{p:Gcon} Assume \eqref{e:lim=gamma} and (A\ref{i:ass_scatter_strong}), let $\hat \gamma$ be a weakly consistent estimator for $\gamma$ and let $(k_n)_{n\in\mathbb{N}}$ be a sequence with $k_n\in\{1, \dots, \lfloor \frac{1}{h}\rfloor \}$ for all $n\in\mathbb{N}$. Then we have
\[ \big(G(\Phi^o_{n,k+k_n}, \Phi^e_{n,k+k_n}, \hat\gamma)\big)_{k=1, \dots, K} \to \big(G(\Phi^o_{k}, \Phi^e_{k}, \gamma)\big)_{k=1, \dots, K}, \qquad n\to\infty, \]
for any $K\in\mathbb{N}$ in distribution. \end{prop}

This proposition follows from Proposition \ref{p:PPcon} and the following lemma by a sharp version of the continuous mapping theorem (see e.g.\ \cite[p.\ 152]{Res87}). 

\begin{lemma}\label{l:Gparts_continuous} Let
\begin{align*}
 \tilde M = \big\{ (\phi^o, \phi^e) \in &M_p\big([0,1]\times (-\infty,0]\big)\times M_p\big([-1,2]\times (-\infty,0]\big) \mid  \mbox{there is no $(\mathfrak{X}, \mathfrak{Y}) \in \phi^o$ with $\tilde g(\mathfrak{X})=\mathfrak{Y}$} \\
& \mbox{ and } \mbox{there are no $(\mathfrak{X}^o, \mathfrak{Y}^o) \in \phi^o$ and $(\mathfrak{X}^e, \mathfrak{Y}^e)\in \phi^e$ with $|\mathfrak{X}^o-\mathfrak{X}^e| \in\mathbb{Z}$} \\
& \mbox{ and } \mbox{for all $k\in\{1,\dots, 24\}$ there is $(\mathfrak{X},\mathfrak{Y})\in \phi^e$ with $\mathfrak{X}\in(-1+\frac{k-1}{8}, -1+\frac{k}{8}]$} \\
& \mbox{ and } \mbox{there is no $x\in(-1,2]$ with $\phi^e(\{x\}\times (-\infty, 0])>1$}\\
& \mbox{ and } \mbox{$\phi^e(\{-1\}\times (-\infty,0])=0$ and $\phi^o(\{0\} \times (-\infty,0])=0$}\big\}.
\end{align*}
Then we have
\begin{enumerate}[(i)]
\item $\mathbb{P}((\Phi_k^o, \Phi_k^e)\in\tilde M)=1$.
\item The map $G_2: M_p\big([0,1]\times (-\infty,0]\big)\times M_p\big([-1,2]\times (-\infty,0]\big)\to\mathbb{R}$ is continuous in all points of $\tilde M$.
\item The map $G_1$, trivially extended to a map $M_p\big([0,1]\times (-\infty,0]\big)\times M_p\big([-1,2]\times (-\infty,0]\big) \to \mathbb{R}$, is continuous in all points of $\tilde M$. 
\item $\lim_{n\to\infty}  S_1^1(h,k)=0.$ 
\end{enumerate}
\end{lemma}

\prf  (i) Clearly, the third, forth and fifth property from the definition of $\tilde M$ are fulfilled with probability $1$. Once these properties are fulfilled given $\Phi_k^e$ there is a set $A\subseteq (0,1]\times (-\infty, 0]$ of Lebesgue measure zero such that
\[ \big\{ (\Phi_k^o, \Phi_k^e)\in\tilde M \big\} = \big\{ \mbox{there is no point $(\mathfrak{X},\mathfrak{Y}) \in \Phi_k^o$ with $(\mathfrak{X},\mathfrak{Y})\in A$} \big\}. \]
Since the event on the right-hand side has a.s.\ conditional probability $1$, these events must also have unconditional probability $1$.

(ii) Let $(\phi^o,\phi^e)\in\tilde M$ and let $(\phi_n^o)_{n\in\mathbb{N}}$ and $(\phi_n^e)_{n\in\mathbb{N}}$ be sequences in $M_p([0,1]\times(-\infty,0])$ and $M_p([-1,2]\times(-\infty,0])$ converging to $\phi^o$ and $\phi^e$ respectively. 

Let $\tilde g_\infty$ denote the estimator $\tilde g$ calculated based on the data set $\phi^e$ and let $\tilde g_n$ be the estimator $\tilde g$ based on the data set $\phi_n^e$. At first we show that there is a constant $L$ independent of $n$ such that
\begin{equation} |\tilde g_\infty(x_2) - \tilde g_\infty(x_1)| \le L \cdot |x_2-x_1| \quad \mbox{and} \quad |\tilde g_n(x_2) - \tilde g_n(x_1)| \le L \cdot |x_2-x_1|,\ n\in\mathbb{N}, \label{eq:tildegLip}\end{equation}
whenever 
\begin{enumerate}[(A)]
\item there is no $(\mathfrak{X},\mathfrak{Y}) \in \phi^e$ or $(\mathfrak{X},\mathfrak{Y}) \in \phi_n^e$ respectively such that $\mathfrak{X}+1\in[x_1\wedge x_2,x_1\vee x_2]$ or $\mathfrak{X}-1\in[x_1\wedge x_2, x_1\vee x_2]$.\label{i:propA}
\end{enumerate}

By the definition of $\tilde M$ we have
\[ S:=\min\big\{ \max\big\{ \mathfrak{Y} \mid (\mathfrak{X},\mathfrak{Y})\in\phi^e, \ \mathfrak{X}\in(-1+\tfrac{k-1}{8}, -1+\tfrac{k}{8}] \big\} \mid k\in\{1,\dots, 24\} \big\} >-\infty.\]
Moreover put
\[ \Delta:= \min\big\{ |\mathfrak{X}-\mathfrak{X}'| \mid (\mathfrak{X},\mathfrak{Y}), (\mathfrak{X}',\mathfrak{Y}') \in \phi^e, (\mathfrak{X},\mathfrak{Y}) \ne (\mathfrak{X}',\mathfrak{Y}'), \mathfrak{Y}>2S-1, \mathfrak{Y}'>2S-1 \big\}. \]
We define
\[ R_n:= \inf\{ r>0 \mid \phi_n^e(B) \le \phi^e(B^r), \ \phi^e(B) \le \phi^e_n(B^r), \quad B\in\mathcal{B}([-1,2] \times (2S-1,0]) \}, \]
where $B^r = \{v\in [-1,2]\times (2S-2, 0] \mid \|b-v\|< r \mbox{ for one $b\in B$} \}$.  
Then we have 
\begin{equation} R_n\to 0. \label{e:Rntozero} \end{equation}
While the proof of \eqref{e:Rntozero} essentially relies on ideas of the proof of the fact that weak convergence of measures implies convergence in the Prohorov metric (see e.g.\ \cite[p.\ 72]{Bil99}), the new parts are so large that we decided to give the complete proof here. 

Let $\epsilon\in(0,\frac{1}{3})$. Let $\{A_i\mid i=1, \dots, k\}$ be a partition of $[-1,2]\times (2S-1,0]$ in finitely many sets such that
\[ \diam A_i := \sup \{\|x-y\| \mid x,y\in A_i\}<\epsilon.\]
Consider the system of open sets 
\[\mathcal{G} := \{(A_{i_1}\cup \dots \cup A_{i_m})^\epsilon \mid \{i_1, \dots, i_m\} \subseteq \{1, \dots, k\}\}. \]
Similar to the portmentau theorem we have $\liminf_{n\to\infty} \phi_n^e(G)\ge \phi^e(G)$ for all $G\in \mathcal{G}$; see e.g.\ \cite[Proposition 3.12]{Res87}. Hence there is $n_0\in\mathbb{N}$ with $\phi_n^e(G) \ge \phi^e(G)$ for all $n\ge n_0$ and $G\in\mathcal{G}$. For $B\in\mathcal{B}([-1,2]\times (2S-1,0])$ we put $B_0:= \bigcup \{ A_k \mid A_k \cap B \ne \emptyset \}$. Now we have
\[ \phi^e(B) \le \phi^e(B_0^\epsilon) \le \phi_n^e(B_0^\epsilon) \le \phi_n^e(B^{2\epsilon}) \]
for all $n\ge n_0$.

On the other hand, consider the system of compact sets
\[\mathcal{C}:= \{ \cl ((A_{i_1}\cup \dots \cup A_{i_m})^\epsilon) \mid \{i_1,\dots, i_m\} \subseteq \{1, \dots, k\}\},\]
where $\cl A$ denotes the closure of $A$. 
Using again \cite[Proposition 3.12]{Res87} we have $\limsup_{n\to\infty} \phi_n^e(K) \le \phi^e(K)$ for all $K\in\mathcal{C}$. Hence there is $n_1\in\mathbb{N}$ with $\phi_n^e(K)\le \phi^e(K)$ for all $n\ge n_1$ and $K\in\mathcal{C}$. So
\[ \phi_n^e(B) \le \phi_n^e(\cl B_0^\epsilon) \le \phi^e(\cl B_0^\epsilon) \le \phi^e(B^{3\epsilon}) \]
for all $n\ge n_1$. So $R_n\le 3\epsilon$. Hence \eqref{e:Rntozero} is proven. 

Let $n$ be so large that $R_n<\min\{\Delta/3, 1/8\}$ and let $x_1,x_2\in[0,1]$ with the property (\ref{i:propA}). There are $(\mathfrak{X},\mathfrak{Y})\in\phi^e_n$ and $(\mathfrak{X}', \mathfrak{Y}')\in\phi^e_n$ such that $(x_1,\tilde g_n(x_1))$, $(\mathfrak{X},\mathfrak{Y})$ and $(\mathfrak{X}',\mathfrak{Y}')$ are on one line. From property (\ref{i:propA}) we get that $|\mathfrak{X}-x_2|<1$ and $|\mathfrak{X}'-x_2|<1$. Hence $(x_2,\tilde g_n(x_2))$ cannot be below that line and therefore we get
\[ \frac{\tilde g_n(x_2)- \tilde g_n(x_1)}{|x_2-x_1|} \ge \frac{\mathfrak{Y}-\mathfrak{Y}'}{|\mathfrak{X}-\mathfrak{X}'|}, \]
where we assumed w.l.o.g.\ that $(x_2-x_1)/(\mathfrak{X}-\mathfrak{X}')>0$. Since the roles of $x_1$ and $x_2$ can be interchanged if we also interchange the roles of $(\mathfrak{X},\mathfrak{Y})$ and $(\mathfrak{X}',\mathfrak{Y}')$, we get the reverse inequality and hence equality holds. 

Now for each $k\in\{1,\dots, 24\}$ there is a point $(\mathfrak{X}_k, \mathfrak{Y}_k)\in\phi_n^e$ with $\mathfrak{X}_k\in (-1+\frac{k-2}{8}, -1+\frac{k+1}{8}]$ and $\mathfrak{Y}_k>S-1/8$. There is $k_1$ with $(-1+\frac{k_1-2}{8}, -1+\frac{k_1+1}{8}] \subseteq [x_1-1,x_1-\frac{1}{2}]$ and $k_2$ with $(-1+\frac{k_2-2}{8}, -1+\frac{k_2+1}{8}] \subseteq [x_1+\frac{1}{2},x_1+1]$. Since neither $(\mathfrak{X}_{k_1},\mathfrak{Y}_{k_1})$ nor $(\mathfrak{X}_{k_2},\mathfrak{Y}_{k_2})$ can lie above the line on which the four points mentioned above lie, we conclude $\mathfrak{Y}>2S-1$ and $\mathfrak{Y}'>2S-1$. Because moreover $|\mathfrak{X}-\mathfrak{X'}| \ge \Delta-2R_n > \Delta/3$, we get 
\[ \frac{|\tilde g_n(x_2)-\tilde g_n(x_1)|}{|x_2-x_1|} \le \frac{6|S|+3}{\Delta}. \]
Clearly, this relation holds also for $\tilde g_\infty$ instead of $\tilde g_n$. For small values of $n$ we may have higher bounds, but since we have a finite bound for each fixed $n$, we have a finite bound independent of $n$. So \eqref{eq:tildegLip} is proven.

Let $n$ be large enough that $R_n<\min\{\Delta/2, 1\}$. Then for each point $(\mathfrak{X}_i^o,\mathfrak{Y}_i^o)\in \phi^o$ with $\mathfrak{Y}_i^o>S-1$ there is a uniquely determined point $(\mathfrak{X}_{i,n}^o,\mathfrak{Y}_{i,n}^o)\in \phi_n^o$ with $\|(\mathfrak{X}_i^o,\mathfrak{Y}_i^o)-(\mathfrak{X}_{i,n}^o,\mathfrak{Y}_{i,n}^o)\| \le R_n$. Clearly every point $(\mathfrak{X}^o,\mathfrak{Y}^o)\in \phi_n^o$ with $\mathfrak{Y}^o>S$ is of this form. We have
\begin{equation} \mathbf{1}_{\{\mathfrak{Y}_{i,n}^o \ge \tilde g_n(\mathfrak{X}_{i,n}^o) \}}  =  \mathbf{1}_{\{\mathfrak{Y}_{i}^o \ge \tilde g_\infty(\mathfrak{X}_{i}^o) \}} \label{e:indi=} \end{equation}
for all sufficiently large $n$. In order to show this, we treat the cases $\mathfrak{Y}_{i}^o \ge \tilde g_\infty(\mathfrak{X}_{i}^o)$  and $\mathfrak{Y}_{i}^o < \tilde g_\infty(\mathfrak{X}_{i}^o)$ jointly. In the first case we have even $\mathfrak{Y}_{i}^o > \tilde g_\infty(\mathfrak{X}_{i}^o)$, since $(\phi^o,\phi^e)\in\tilde M$. Put $\delta:= |\tilde g_\infty(\mathfrak{X}_i^o)-\mathfrak{Y}_i^o|$. Let $n$ be large enough that $R_n<\min\{ \frac{\delta}{4}, \frac{\delta}{4(1+L)} \}$, where $L$ is the constant from \eqref{eq:tildegLip}. Then $|\mathfrak{Y}_i^o-\mathfrak{Y}_{i,n}^o|<\delta/4$ and $|\mathfrak{X}_i^o-\mathfrak{X}_{i,n}^o|<\frac{\delta}{4L}$. Recall that there is no $(\mathfrak{X}^e, \mathfrak{Y}^e) \in \phi^e$ with $\mathfrak{X}^e-\mathfrak{X}_i^o\in\{\pm 1\}$. Hence $\mathfrak{X}^e-\mathfrak{X}_{i,n}^o$ cannot convergence to $\pm 1$ -- not even along a subsequence. Thus for sufficiently large $n$ Assumption \eqref{i:propA} is fulfilled with $x_1=\mathfrak{X}_{i,n}^o$ and $x_2=\mathfrak{X}_i^o$ and hence \eqref{eq:tildegLip} holds. So $|\tilde g_n(\mathfrak{X}_{i,n}^o)-\tilde g_n(\mathfrak{X}_i^o)|<\delta/4$. Moreover $|\tilde g_\infty(\mathfrak{X}_{i}^o)-\tilde g_n(\mathfrak{X}_{i}^o)|<(1+L)R_n$ implies $|\tilde g_\infty(\mathfrak{X}_{i}^o)-\tilde g_n(\mathfrak{X}_{i}^o)|<\delta/4$ and hence
\[ | \tilde g_n(\mathfrak{X}_{i,n}^o)-\tilde g_\infty(\mathfrak{X}_i^o)|  \le  |\tilde g_n(\mathfrak{X}_{i,n}^o)-\tilde g_n(\mathfrak{X}_i^o)|+ |\tilde g_n(\mathfrak{X}_{i}^o)-\tilde g_\infty(\mathfrak{X}_i^o)| \le \frac{\delta}{4}+\frac{\delta}{4}=\frac{\delta}{2}. \]
Thus \eqref{e:indi=} holds true. 

Now we get
\begin{align*}
\sum_{(\mathfrak{X}_{i,n}^o, \mathfrak{Y}_{i,n}^o)\in \phi^o_n} \mathfrak{Y}_{i,n}^o \mathbf{1}_{\{\mathfrak{Y}_{i,n}^o \ge \tilde g_n(\mathfrak{X}_{i,n}^o)\}} &= \sum_{(\mathfrak{X}_{i,n}^o, \mathfrak{Y}_{i,n}^o)\in \phi^o_n\atop  \mathfrak{Y}_{i,n}^o\ge S-1} \mathfrak{Y}_{i,n}^o \mathbf{1}_{\{\mathfrak{Y}_{i,n}^o \ge \tilde g_n(\mathfrak{X}_{i,n}^o)\}} \\
&\le \sum_{(\mathfrak{X}_{i}^o, \mathfrak{Y}_{i}^o)\in \phi^o\atop \mathfrak{Y}_{i}^o\ge S} (\mathfrak{Y}_{i}^o+R_n) \mathbf{1}_{\{\mathfrak{Y}_{i}^o \ge \tilde g_\infty(\mathfrak{X}_{i}^o)\}} \\
 &\le \sum_{(\mathfrak{X}_{i}^o, \mathfrak{Y}_{i}^o)\in \phi^o} \mathfrak{Y}_{i}^o \mathbf{1}_{\{\mathfrak{Y}_{i}^o \ge \tilde g_\infty(\mathfrak{X}_{i}^o)\}} +R_n \cdot \nu,
\end{align*}
where $\nu:= \phi^o([0,1]\times(S-1,0])$, and 
\begin{align*}
\sum_{(\mathfrak{X}_{i,n}^o, \mathfrak{Y}_{i,n}^o)\in \phi^o_n} \mathfrak{Y}_{i,n}^o \mathbf{1}_{\{\mathfrak{Y}_{i,n}^o \ge \tilde g_n(\mathfrak{X}_{i,n}^o)\}} &= \sum_{(\mathfrak{X}_{i,n}^o, \mathfrak{Y}_{i,n}^o)\in \phi^o_n\atop \mathfrak{Y}_{i,n}^o\ge S-\frac{1}{2}} \mathfrak{Y}_{i,n}^o \mathbf{1}_{\{\mathfrak{Y}_{i,n}^o \ge \tilde g_n(\mathfrak{X}_{i,n}^o)\}} \\
&\ge \sum_{(\mathfrak{X}_{i}^o, \mathfrak{Y}_{i}^o)\in \phi^o\atop \mathfrak{Y}_{i}^o\ge S-1} (\mathfrak{Y}_{i}^o-R_n) \mathbf{1}_{\{\mathfrak{Y}_{i}^o \ge \tilde g_\infty(\mathfrak{X}_{i}^o)\}} \\
 &\ge \sum_{(\mathfrak{X}_{i}^o, \mathfrak{Y}_{i}^o)\in \phi^o} \mathfrak{Y}_{i}^o \mathbf{1}_{\{\mathfrak{Y}_{i}^o \ge \tilde g_\infty(\mathfrak{X}_{i}^o)\}} -R_n \cdot \nu.
\end{align*}
Hence 
\[ \lim_{n\to\infty} G_2(\phi_{n}^o,\phi_{n}^e) = G_2(\phi^o, \phi^e). \]
(iii)
We have
\[ | \tilde g_\infty(x) - \tilde g_n(x) | \le \begin{cases} 
(1+L)R_n & \mbox{if there is no $(\mathfrak{X},\mathfrak{Y})\in \phi^e$ with $\mathfrak{X} \in [x-R_n-1, x+R_n-1) \cup [x-R_n+1, x+R_n+1)$} \\
& \mbox{and $\mathfrak{Y} \ge 2S-1$}\\
|S|+R_n & \mbox{anyway}
\end{cases}\] 
and hence
\[ \Big|\int_0^1 \tilde g_\infty(x)^2 \, dx - \int_0^1 \tilde g_n(x)^2 \, dx\Big| \le 2(1+L)\cdot R_n \cdot (|S|+R_n) + 2\phi^e([-1,2]\times(2S-1,0]) \cdot R_n \cdot (|S|+R_n)^2 \longrightarrow 0. \] 
(iv) Observe
\[S_1^1(h,k) = \frac{(nh)^2}{nh} \sum_{i:x_{2i-1}\in ((k-1)h,kh]} (\hat g(x_{2i-1}) - g(x_{2i-1}))^2 - \frac{(nh)^2}{2h} \int_{(k-1)h}^{kh} (\hat g(x) -g(x))^2 \, dx = \frac{1}{2} \big( \frac{2}{nh}\sum_{i=1}^\mathfrak{n} \tilde g_n(\mathfrak{X}_i) -\int_0^1 \tilde g_n(x)\, dx\big), \]
where $\{\frac{x_{2i-1}}{h}-(k-1)\mid x_{2i-1}\in[(k-1)h, kh]\} = \{ \mathfrak{X}_i \mid i=1, \dots, \mathfrak{n}\}.$
Let $\nu:\mathbb{N}\to (0,\infty)$ be a function with $\lim_{n\to\infty} \nu(n)=0$ and
\[ |\frac{\Phi^o_{k,n}\big([x,x+\tfrac{H}{h}) \times (-\infty,0]\big)}{ \frac{nH}{2}} -1| <\nu(n)  \mbox{ for all $x\in[0,1-\frac{H}{h}]$ and $k\in\{1, \dots, \lfloor \frac{1}{h} \rfloor\}$ }.  \]
Then we have
\begin{align} \big|\frac{2}{nh} \sum_{i=1}^{\mathfrak{n}} \tilde g_n(\mathfrak{X_i})^2 - \int_0^1 \tilde g_n(x)^2 \, dx\big| &\le (1+\nu(n)) (1+\tfrac{H}{h})\tfrac{H}{h}\cdot L \cdot (|S|+R_n) \notag \\
&+ \Phi^e_{k}([-1,2]\times(2S-1,0]) \cdot \big(\tfrac{H}{h}+2R_n\big)(1+\nu(n)) \cdot (|S|+R_n)^2\notag\\ &+ \big(1+\tfrac{H}{h}\big)\cdot\nu(n) \cdot (|S|+R_n)^2 .  \label{e:diff_Rie_int}\end{align}
Indeed, partition the interval [0,1] in $\lfloor \frac{h}{H}\rfloor$ intervals of the form $[a,a+H/h)$ and one shorter interval. For each of these intervals we have
\begin{align*}
 \big|\frac{2}{nh} &\sum_{i: \mathfrak{X_i}\in[a,a+H/h)} \tilde g_n(\mathfrak{X_i})^2 - \int_a^{a+H/h} \tilde g_n(x)^2 \, dx\big| \\
&\le  \big|\frac{2}{nh} \sum_{i: \mathfrak{X_i}\in[a,a+H/h)} \tilde g_n(\mathfrak{X_i})^2 - \frac{2\#\{i\mid \mathfrak{X_i}\in[a,a+H/h)\}}{nH}  \int_a^{a+H/h} \tilde g_n(x)^2 \, dx\big|  +  \nu(n)\frac{H}{h}(|S|+R_n)^2 
\end{align*}
and 
\begin{align*} \big|\frac{2}{nh} &\sum_{i: \mathfrak{X_i}\in[a,a+H/h)} \tilde g_n(\mathfrak{X_i})^2 - \frac{2\#\{i\mid \mathfrak{X_i}\in[a,a+H/h)\}}{nh} \frac{1}{H/h} \int_a^{a+H/h} \tilde g_n(x)^2 \, dx\big| \\
& = \big|\frac{2}{nh} \sum_{i: \mathfrak{X_i}\in[a,a+H/h)} \Big(\tilde g_n(\mathfrak{X_i})^2 - \frac{1}{H/h} \int_a^{a+H/h} \tilde g_n(x)^2 \, dx\Big)\big|\\
&\le \begin{cases} 
\frac{2\#\{i\mid \mathfrak{X_i}\in[a,a+H/h)\}}{nh}L \frac{H}{h}(|S|+R_n)& \mbox{if there is no $(\mathfrak{X},\mathfrak{Y})\in \phi_k^e$ with}
\\ &\mathfrak{X} \in [a-R_n-1, a+R_n+\frac{H}{h}-1) \cup [a-R_n+1, a+R_n+\frac{H}{h}+1)\\
\frac{2\#\{i\mid \mathfrak{X_i}\in[a,a+H/h)\}}{nh}(|S|+R_n)^2 & \mbox{anyway}.
\end{cases}
\end{align*} 
Since this calculation also holds for the shorter interval, we arrive at \eqref{e:diff_Rie_int}. 
Clearly the right-hand side and thus the left-hand side of \eqref{e:diff_Rie_int} converges to $0$ as $n\to\infty$. 
\qed\medskip

We derive the following consequence of Proposition \ref{p:Gcon} and Lemma \ref{l:Gparts_continuous}(iv).

\begin{cor}
Assume that \eqref{e:lim=gamma} and  (A\ref{i:ass_scatter_strong}) are fulfilled, let $\hat \gamma$ be a weakly consistent estimator for $\gamma$ and let $(k_n)_{n\in\mathbb{N}}$ be a sequence with $k_n\in\{1, \dots, \lfloor \frac{1}{h}\rfloor \}$ for all $n\in\mathbb{N}$. Then we have
\[ \big(S_1(h,k+k_n)\big)_{k=1, \dots, K} \to \big(G(\Phi^o_{k}, \Phi^e_{k}, \gamma)\big)_{k=1, \dots, K}, \qquad n\to\infty, \]
for any $K\in\mathbb{N}$ in distribution.
\end{cor}

\begin{prop}\label{p:var_con}
Assume that \eqref{e:lim=gamma}, \eqref{e:ass_distr_tail} and (A\ref{i:ass_scatter_strong}) are fulfilled. Let $\hat \gamma$ be a weakly consistent estimator for $\gamma$ fulfilling (G\ref{e:gamma_large-dev}), let $(k_n)_{n\in\mathbb{N}}$ be a sequence with $k_n\in\{1, \dots, \lfloor \frac{1}{h}\rfloor \}$ for all $n\in\mathbb{N}$ and let $k,l\in\mathbb{N}$.  Then we have
\[ \lim_{n\to\infty} \Cov\big( S_1(h,k+k_n), S_1(h,l+k_n) \big) = \Cov\big( G(\Phi_k^o,\Phi_k^e, \gamma), G(\Phi_{l}^o, \Phi_{l}^e, \gamma ) \big). \]
\end{prop}
\prf All we have to show is the uniform integrability
\begin{equation} \int_0^\infty \sup_n \mathbb{P}(|S_1(h,k)|>s)s\, ds<\infty.\label{e:uni_int} \end{equation}
We have
\begin{align*}
\mathbb{P}(|S_1(h,k)|>s) \le &\mathbb{P}\big(|nh \sum_{i:x_{2i-1}\in((k-1)h,kh]} (\hat g(x_{2i-1}) - g(x_{2i-1}))^2 |>\tfrac{s}{2}\big) + \mathbb{P}\big(\hat \gamma < \tfrac{2}{\sqrt{s}}\big)\\ & + \mathbb{P}\big(|nh\sum_{i:x_{2i-1}\in ((k-1)h,kh]} (Y_{2i-1}-g(x_{2i-1}))\mathbf{1}_{\{Y_{2i-1} \ge \hat g(x_{2i-1})\}}|>\tfrac{\sqrt{s}}{2}\big). 
\end{align*}
So let us treat these three summands one by one. We have 
\begin{align*}
\mathbb{P}\big(|nh\sum_{i:x_{2i-1}\in ((k-1)h,kh]}& (Y_{2i-1}-g(x_{2i-1}))\mathbf{1}_{\{Y_{2i-1} \ge \hat g(x_{2i-1})\}}|>\tfrac{\sqrt{s}}{2}\big)\\
& \le  \mathbb{P}\Big( \Big\{ hn\sum_{i: \frac{x_{2i-1}}{h}\in(k-1,k]} \epsilon_{2i-1}\cdot \mathbf{1}_{\{ \epsilon_{2i-1}\ge \hat g(x_{2i-1}) - g(x_{2i-1}) \}} < -\frac{\sqrt{s}}{2} \Big\} \cap \Big\{ m_\epsilon< \frac{1}{nh} \sqrt{\tfrac{\sqrt{s}}{8\bar\gamma}} \Big\} \Big)\\
+& \mathbb{P}\big(m_\epsilon \ge \frac{1}{nh} \sqrt{\tfrac{\sqrt{s}}{8\bar\gamma}} \big), 
\end{align*}
 where 
\[ m_\epsilon:= \max\big\{ \big| \max\{\epsilon_{2i}\mid x_{2i}\in ((j-1)h/2,jh/2] \} \big| \mid j=2k-2, \dots, 2k+1 \big\} \]
and $\bar\gamma$ is a constant such that
\[ 1-F(t) \le \bar\gamma\cdot |t|, \quad t\le 0.\]
Now Proposition \ref{p:JMR3.1} yields
\begin{align*}
\mathbb{P}\Big( \Big \{hn \sum_{i: \frac{x_{2i-1}}{h}\in(k-1,k]} &\epsilon_{2i-1} \mathbf{1}_{\{ \epsilon_{2i-1}\ge \hat g(x_{2i-1}) - g(x_{2i-1}) \}} < -\frac{\sqrt{s}}{2} \Big\} \cap \Big\{ m_\epsilon< \frac{1}{nh} \sqrt{\tfrac{\sqrt{s}}{8\bar\gamma}} \Big\} \Big)\\
& \le \mathbb{P}\Big( hn \sum_{i: \frac{2i-1}{h}\in(k-1,k]} -\frac{1}{nh} \sqrt{\tfrac{\sqrt{s}}{8\bar\gamma}} \mathbf{1}_{\{ \epsilon_{2i-1}\ge -\frac{1}{nh} \sqrt{\frac{\sqrt{s}}{8\bar\gamma}} \}} < -\frac{\sqrt{s}}{2}  \Big)\\
& = \mathbb{P}\Big( \sum_{i: \frac{2i-1}{h}\in(k-1,k]}  \mathbf{1}_{\{ \epsilon_{2i-1}\ge -\frac{1}{nh} \sqrt{\frac{\sqrt{s}}{8\bar\gamma}} \}} > \sqrt{2\sqrt{s}\bar\gamma}  \Big).
\end{align*}
The indicators $\mathbf{1}_{\{ \epsilon_{2i-1}\ge -\frac{1}{nh} \sqrt{\frac{\sqrt{s}}{8\bar\gamma}} \}}$ are independent with success probability 
\[ 1-F\big(-\frac{1}{nh} \sqrt{\tfrac{\sqrt{s}}{8\bar\gamma}} \big) \le \frac{\sqrt{\sqrt{s}\bar\gamma}}{\sqrt{8} nh}  . \]
Hence Lemma \ref{l:binom}, \ref{l:poisson} and \ref{l:stirling} give 
\begin{align*}
\mathbb{P}\Big( \sum_{i: \frac{x_{2i-1}}{h}\in(k-1,k]} & \mathbf{1}_{\{ \epsilon_{2i-1}\ge -\frac{1}{nh} \sqrt{\frac{\sqrt{s}}{8\bar\gamma}} \}} > \sqrt{2\sqrt{s}\bar\gamma}  \Big)
 \le 2\binom{m}{\mathfrak{k}} p^\mathfrak{k}(1-p)^\mathfrak{k} \le 2\frac{(mp)^\mathfrak{k}}{\mathfrak{k}!} \le 2\frac{(mp)^\mathfrak{k}}{(\frac{\mathfrak{k}}{e})^\mathfrak{k}\sqrt{\frac{\pi \mathfrak{k}}{2}}} \le \sqrt{\frac{8}{\pi \mathfrak{k}}} \big( \frac{e}{4} \big)^\mathfrak{k},
\end{align*}
where
\[ \mathfrak{k}= \big\lceil \sqrt{2\sqrt{s}\bar\gamma}\rceil , \qquad p = \frac{\sqrt{\sqrt{s}\bar\gamma}}{\sqrt{8} nh},\qquad m = \nu^+(n), \]
provided that $n$ is large enough that $\nu^+(n)<nh$ and hence $\mathfrak{k}>4mp$ -- recall $\nu^-(n):=\inf\{ \#\{i\in\{1, \dots, n\}\mid x_i \in [x,x+h/2)\} \mid x\in[-h/2,1] \}$ and $\nu^+(n):=\sup\{ \#\{i\in\{1, \dots, n\}\mid x_i \in [x,x+h/2)\} \mid x\in[-h/2,1] \}$.
Therefore
\begin{align}
 \int_0^\infty \sup_n \mathbb{P}\Big( \big \{\sum_{i: \frac{x_{2i-1}}{h}\in(k-1,k]}& \epsilon_{2i-1}\cdot hn \mathbf{1}_{\{ \epsilon_{2i-1}\ge \hat g(x_{2i-1}) - g(x_{2i-1}) \}} < -\frac{\sqrt{s}}{2} \big\} \cap \big\{ m_\epsilon< \frac{1}{nh} \sqrt{\tfrac{\sqrt{s}}{8\bar\gamma}} \big\} \Big)s\, ds \notag\\
&\le \int_0^\infty \sqrt{\tfrac{8}{\pi \lceil \sqrt{2\sqrt{s}\bar\gamma}\rceil}} \big( \frac{e}{4} \big)^{\lceil \sqrt{2\sqrt{s}\bar\gamma}\rceil} s \, ds 
<\infty.\label{e:uni_int_1}\end{align}
Turning to the other summand we have
\begin{align*}
\mathbb{P}\big(m_\epsilon \ge \frac{1}{nh} \sqrt{\tfrac{\sqrt{s}}{8\bar\gamma}} \big) &\le \sum_{j=2k-2}^{2k+1} \mathbb{P}\big( \max\{\epsilon_{2i}\mid x_{2i}\in ((j-1)h/2,jh/2] \} \le -\frac{1}{nh} \sqrt{\tfrac{\sqrt{s}}{8\bar\gamma}}\big)
 \le 4 F\big( -\frac{1}{nh} \sqrt{\tfrac{\sqrt{s}}{8\bar\gamma}} \big)^{\lfloor \nu^-(n)/2\rfloor}.
\end{align*}
Recall that by \eqref{e:ass_distr_tail} there are $\Gamma>0$ and $\tilde C_F\ge 1$ with $F(t)\le \tilde C_F \cdot |t|^{-\Gamma}$ for all $t\in (-\infty,-1)$. We consider at first the case that $\frac{\tilde C_F^{-1/\Gamma}}{nh}\sqrt{\tfrac{\sqrt{s}}{8\bar\gamma}}>2$. We get
\begin{align*}
F\big( -\frac{1}{nh} \sqrt{\tfrac{\sqrt{s}}{8\bar\gamma}} \big)^{\lfloor \nu^-(n)/2\rfloor} &\le \big( \tilde C_F |\frac{1}{nh}\sqrt{\tfrac{\sqrt{s}}{8\bar\gamma}}|^{-\Gamma} \big)^{\lfloor \nu^-(n)/2\rfloor} \\
& = \Big( \big|\frac{\tilde C_F^{-1/\Gamma}}{nh}\sqrt{\tfrac{\sqrt{s}}{8\bar\gamma}} \big|^{\Gamma\lfloor \nu^-(n)/2\rfloor/9} \Big)^{-9} \\
&\le \Big( \big|\frac{\tilde C_F^{-1/\Gamma}}{nh}\sqrt{\tfrac{\sqrt{s}}{8\bar\gamma}}-1 \big|\cdot{\Gamma\lfloor \nu^-(n)/2\rfloor/9} +1\Big)^{-9} \\
&\le \Big( \big|\frac{\tilde C_F^{-1/\Gamma}}{nh}\sqrt{\tfrac{\sqrt{s}}{8\bar\gamma}} \big|\cdot{\Gamma \nu^-(n)/72} \Big)^{-9} \\
&= (144\cdot\sqrt{2})^9 \cdot (\tilde C_x)^9\tilde C_F^{9/\Gamma} \big( \frac{1}{\Gamma}\sqrt{\tfrac{\bar\gamma}{\sqrt{s}}} \big)^{9},
\end{align*}
where $\tilde C_x$ is a constant with $\tilde C_x nh\le \nu^-(n)$ for all $n\in\mathbb{N}$ and we assumed $\nu^-(n) \ge 4$. 
Now turn to the case $\frac{\tilde C_F^{-1/\Gamma}}{nh}\sqrt{\tfrac{\sqrt{s}}{8\bar\gamma}} \le 2$.  There is $\tilde\gamma>0$ with 
\[ F(t) \le \exp\{\tilde\gamma t\}, \qquad t\in [-2\tilde C_F^{1/\Gamma}, 0). \]
Hence
\[ 
F\big( -\frac{1}{nh} \sqrt{\tfrac{\sqrt{s}}{8\bar\gamma}} \big)^{\lfloor \nu^-(n)/2\rfloor}  \le    \exp\Big\{ -\frac{\tilde\gamma}{nh} \sqrt{\tfrac{\sqrt{s}}{8\bar\gamma}} {\lfloor \nu^-(n)/2\rfloor}\Big\}  
\le    \exp\Big\{ -\frac{\tilde\gamma}{nh} \sqrt{\tfrac{\sqrt{s}}{8\bar\gamma}} \nu^-(n)/4\Big\} \le   \exp\Big\{ -\frac{\tilde\gamma}{8} \tilde C_x \sqrt{\tfrac{\sqrt{s}}{2\bar\gamma}} \Big\} 
\]
Altogether we get
\begin{equation} \int_0^\infty \sup_n \mathbb{P}\big(m_\epsilon \ge \frac{1}{nh} \sqrt{\tfrac{\sqrt{s}}{8\bar\gamma}} \big) s\, ds \le \int_0^\infty 4\Big( \min\Big\{ (144\cdot\sqrt{2})^9\cdot \tilde C_x^9 \tilde C_F^{9/\Gamma} \big( \frac{1}{\Gamma}\sqrt{\tfrac{\bar\gamma}{\sqrt{s}}} \big)^{9}, 1 \Big\} + \exp\big\{ -\frac{\tilde\gamma}{8} \tilde C_x \sqrt{\tfrac{\sqrt{s}}{2\bar\gamma}} \big\}\Big) s\, ds,\label{e:uni_int2}\end{equation}
where the supremum is taken over all $n$ with $\nu^-(n) \ge 4$. 
Moreover,
\begin{equation} \int_0^\infty \sup_n \mathbb{P}\big(\hat \gamma <\tfrac{2}{\sqrt{s}}\big) s\, ds <\infty \label{e:uni_int3} \end{equation}
follows from (G\ref{e:gamma_large-dev}). 
 
Since
\[ \big( \hat g(x)- g(x) \big)^2 \le (m_\epsilon)^2, \qquad x\in[0,1], \]
we get
\[ \mathbb{P}\big(|nh \sum_{i:x_{2i-1}\in((k-1)h,kh]} (\hat g(x_{2i-1}) - g(x_{2i-1}))^2|>\frac{s}{2}\big) \le \mathbb{P}\Big(nh \cdot \nu^+(n) (m_\epsilon)^2>\frac{s}{2}\Big). \]
Now one gets 
\[ \int_0^\infty \sup_n \mathbb{P}(|nh \sum_{i:x_{2i-1}\in((k-1)h,kh]} (\hat g(x_{2i-1}) - g(x_{2i-1}))^2|>\tfrac{s}{2}) s\, ds <\infty \]
the same way as \eqref{e:uni_int2}. 
Together with \eqref{e:uni_int_1}, \eqref{e:uni_int2} and \eqref{e:uni_int3} this implies \eqref{e:uni_int} and hence the assertion of the proposition. \qed

\begin{remark}
Probably, one can show 
\[ \lim_{n\to\infty} \Ex\big[ S_1(h,k) \big] = \Ex\big[ G(\Phi_k^o,\Phi_k^e, \gamma) \big], \]
similar to Proposition \ref{p:var_con}, but we have not checked completely. Anyway, this relation is not too interesting, since both sides of the equality are zero. Showing a non-degenerate limit relation involving the same terms is an open problem that will require different methods. 
\end{remark}

\begin{theorem}\label{T:T_con_P} 
Assume that the errors satisfy \eqref{e:ass_distr} and \eqref{e:ass_distr_tail}, that $\hat\gamma$ is an estimator for $\gamma$ with (G\ref{i:G1}), (G\ref{i:G4}) and (G\ref{e:gamma_large-dev}) and that $\lim_{n\to\infty} nh^2=\infty$ and (A\ref{i:ass_scatter_strong}) hold. Then we have
\[ T \sqrt{n^2h^3} \to \mathcal{N}\Big(0,\sum_{l=1}^5\Cov\big( G(\Phi_3^o,\Phi_3^e, \gamma), G(\Phi_l^o, \Phi_l^e, \gamma) \big) \Big), \qquad n\to\infty,\]
in distribution.
\end{theorem}
\prf From Corollary \ref{c:CLT_T} we have
\[\frac{ T}{\sqrt{\Var(S_1)}} \to \mathcal{N}(0,1), \qquad n\to\infty,\]
in distribution. We can treat the variance of $S_1$ by Proposition \ref{p:var_con} except for a remainder part $\mathfrak{R}$, where we recall
\begin{align*}
\mathfrak{R} 
& = hn \Big(\sum_{i:\frac{2i-1}{h} \in (\lfloor 1/h \rfloor, 1/h]} (\hat g(x_{2i-1})-g(x_{2i-1}))^2 + \frac{2}{\hat\gamma} \big(Y_{2i-1}-g(x_{2i-1})\big) \mathbf{1}_{\{Y_{2i-1}\ge \hat g(x_{2i-1})\}}\Big).
\end{align*}
Observe that $\hat g(x_{2i-1})$ is positively correlated for different $i$, while the $Y_{2i-1}$ are independent. Hence the summands are positively correlated and so the variance of the sum increases with increasing number of summands. In particular $\Var(\mathfrak{R}) \le \Var( S_1(h,k))$. Thus we get 
\begin{equation} \lim_{n\to\infty} n^2h^3 \Var(S_1) = \lim_{n\to\infty} h \sum_{k=1}^{\lfloor 1/h \rfloor} \sum_{l=1}^{\lfloor 1/h \rfloor} \Cov\big( S_1(h,k), S_1(h,l) \big) = \sum_{l=1}^5\Cov\big( G(\Phi_3^o,\Phi_3^e, \gamma), G(\Phi_{l}^o, \Phi_{l}^e, \gamma) \big). \label{e:S1} \end{equation}
Indeed, if $k_n$ denotes for any $n\in\mathbb{N}$ the number with
\[ \sum_{l=-2}^2  \Cov\big( S_1(h,k_n), S_1(h,k_n+l) \big) = \max\big\{ \sum_{l=-2}^2  \Cov\big( S_1(h,k), S_1(h,k+l) \big) \mid k\in\{1, \dots, \lfloor \frac{1}{h} \rfloor \} \big\},\]
then 
\[ \lim_{n\to\infty} h \sum_{k=1}^{\lfloor 1/h \rfloor} \sum_{l=1}^{\lfloor 1/h \rfloor} \Cov\big( S_1(h,k), S_1(h,l) \big)  \le \sum_{l=-2}^2  \Cov\big( S_1(h,k_n), S_1(h,k_n+l) \big) =  \sum_{l=1}^5\Cov\big( G(\Phi_3^o,\Phi_3^e, \gamma), G(\Phi_{l}^o, \Phi_{l}^e, \gamma) \big).\]
On the other hand, if we let $\bar k_n$ denote the number with 
\[ \sum_{l=-2}^2  \Cov\big( S_1(h,\bar k_n), S_1(h,\bar k_n+l) \big) = \min\big\{ \sum_{l=-2}^2  \Cov\big( S_1(h,k), S_1(h,k+l) \big) \mid k\in\{1, \dots, \lfloor \frac{1}{h} \rfloor \} \big\},\]
then 
\[ \lim_{n\to\infty} h \sum_{k=1}^{\lfloor 1/h \rfloor} \sum_{l=1}^{\lfloor 1/h \rfloor} \Cov\big( S_1(h,k), S_1(h,l) \big)  \ge \sum_{l=-2}^2  \Cov\big( S_1(h,\bar k_n), S_1(h,\bar k_n+l) \big) =  \sum_{l=1}^5\Cov\big( G(\Phi_3^o,\Phi_3^e, \gamma), G(\Phi_{l}^o, \Phi_{l}^e, \gamma) \big).\]
Since both sides of \eqref{e:S1} are positive and finite -- see e.g.\ Lemma \ref{l:variances} -- we get
\[ \Var(S_1) \sim n^{-2}h^{-3} \cdot \sum_{l=1}^5\Cov\big( G(\Phi_3^o,\Phi_3^e, \gamma), G(\Phi_{l}^o, \Phi_{l}^e, \gamma) \big). \]
Hence
\[ \frac{T}{\sqrt{n^{-2}h^{-3} \cdot \sum_{l=1}^5\Cov\big( G(\Phi_3^o,\Phi_3^e, \gamma), G(\Phi_{l}^o, \Phi_{l}^e, \gamma) \big)}} \to \mathcal{N}(0,1), \quad n\to\infty, \]
which is equivalent to the assertion.  \qed\medskip

Recall $A_\gamma:= \sum_{l=1}^5\Cov\big( G(\Phi_3^o,\Phi_3^e, \gamma), G(\Phi_{l}^o, \Phi_{l}^e, \gamma) \big)$, where the intensity of the Poisson processes is $\gamma$.

\begin{lemma}\label{l:Agamma} We have $A_\gamma=A_1/\gamma^4$. \end{lemma}
\prf Let $\Phi^{o,\gamma}$ and $\Phi^{e,\gamma}$ be Poisson processes on $[0,1]\times (-\infty,0]$ or $[-1,2]\times (-\infty,0]$ respectively of intensity $\gamma$. Define maps
\[L_\gamma: [0,1]\times (-\infty,0] \to [0,1]\times (-\infty, 0],\ (x,y)\mapsto (x, \gamma y), \qquad L'_\gamma: [-1,2]\times (-\infty,0] \to [-1,2]\times (-\infty, 0],\ (x,y)\mapsto (x, \gamma y)\]
and denote the induced maps $M_p([0,1]\times (-\infty,0]) \to M_p([0,1] \times (-\infty,0])$ respectively $M_p([-1,2]\times (-\infty,0]) \to M_p([-1,2]\times (-\infty,0])$ by the same symbols. Then
\[ G(\Phi^{o,1}, \Phi^{e,1}, 1) \stackrel{d}{=} G(L_\gamma \Phi^{o,\gamma}, L'_\gamma \Phi^{e,\gamma}, 1) = \gamma^2 G(\Phi^{o,\gamma}, \Phi^{e,\gamma}, \gamma).\]
Hence $A_\gamma=\gamma^{-4}A_1$.\qed


\begin{remark} From the proof of Lemma \ref{l:variances} we see $\Var(S_1) \sim \Ex[ \Var(S_1''\mid \hat g)]$. So one can also construct a test based on the Poisson approximation of $\Ex[ \Var(S_1''\mid \hat g)]$ -- however, this idea seems less natural and thus we did not carry it out. 
\end{remark}

\subsection{Proof of Theorem \ref{c:test} and Theorem \ref{c:tests}}\label{ss:proofs}

In this subsection we prove Theorem \ref{c:test} and Theorem \ref{c:tests}. 

\prf[ of Theorem \ref{c:test}]
From Corollary \ref{c:VarS1} and Corollary \ref{c:CLT_T} we get
\begin{align*}
\lim_{n\to \infty} \mathbb{E}[\varphi_1]	&= \lim_{n\to\infty} \mathbb{P}\Big(T \ge   z_{1-\lambda} \cdot \sqrt{\frac{8}{(C_x)^3} n^{-2}h^{-3}\hat\gamma^{-4}} \Big) \\
& \le \lim_{n\to\infty} \mathbb{P}\big(T \ge z_{1-\lambda} \cdot \sqrt{\Var(S_1)} \big)\\
& = \lambda. \qquad \qed 
\end{align*}

\prf[ of Theorem \ref{c:tests}]
 By Slutzky's theorem, Lemma \ref{l:Agamma} and Theorem \ref{T:T_con_P} we have
\[ \lim_{n\to\infty}\Ex[\varphi_2] =\lim_{n\to\infty}\mathbb{P}\big(T \ge z_{1-\lambda} \sqrt{n^{-2}h^{-3}\hat\gamma^{-4}A_1}\big) =\lim_{n\to\infty}\mathbb{P}\big(T \ge z_{1-\lambda} \sqrt{n^{-2}h^{-3}A_\gamma}\big) =\lambda. \qed \]

\subsection{Consistency}\label{ss:consistency}

In this section we will show that the tests $\varphi_1$ and $\varphi_2$ are consistent. We start by showing the counterpart of Proposition \ref{p:asy_ex} under the alternative. 

\begin{lemma}\label{l:asy_ex_alt} Let $g$ be uniformly continuous, let the errors satisfy \eqref{e:ass_distr_tail} and let the sample points fulfill \eqref{e:ass_scatter_alt}. Then we have
\[\lim_{n\to\infty} \sup\big\{\Ex\big[ |\hat g(x)-g(x)|^k \big]\mid x\in[0,1] \big\}= 0. \]
\end{lemma}
\prf Put $Z_1(x):=\max\{\epsilon_{2i} \mid x_{2i} \in (x-h,x) \}$ and $Z_2(x):=\max\{\epsilon_{2i} \mid x_{2i} \in (x,x+h) \}$ and let $Z_n$ be the maximum of $\nu^-(n)$ independent random variables with distribution function $F$, where $\nu^-(n):=\min\{ \# \{i\in\{1, \dots, n\} \mid x_i\in (x, x+h/2)\} \mid x\in(-h/2, 1)\}$. Then $Z_1(x)\le Z_n$ a.s.\ and $Z_2(x)\le Z_n$ a.s.\ for appropriate couplings.  Moreover, $Z_n\to 0$ in probability as $n\to \infty$, which implies $\Ex[|Z_n|^k] \to 0$ as $n\to\infty$, since
\[ \int_0^\infty \sup_n \mathbb{P}(|Z_n|>t) t^{k-1} \, dt =  \int_0^\infty \sup_n F(-t)^{\nu^-(n)} t^{k-1} \, dt \le \int_0^\infty \min\{ \tilde C_F^{(k+1)/\Gamma} t^{-(k+1)}, 1 \} t^{k-1}\, dt < \infty, \]
where the supremum is taken over all $n$ such that $\nu^-(n)\Gamma \ge k+1$. 
Let
\[\omega_g(\delta) := \sup\{ |g(x_1)-g(x_2)| \mid |x_1-x_2| \le \delta \} \]
be the modulus of continuity of $g$. Then
\[ |\hat g(x)-g(x)| \le \omega_g(h)+\max\{|Z_1(x)|, |Z_2(x)|\} \]
and hence
\[ \Ex\big[ |\hat g(x) -g(x)|^k \big] \le 2^k\Ex[|Z_1(x)|^k] + 2^k\Ex[|Z_2(x)|^k] +2^k\omega_g(h)^k \longrightarrow 0, \] 
as $n\to\infty$ uniformly in $x$. \qed\medskip

Put
\[
a= \min_{m,c}\int_0^1 \Big( g(x)-mx-c\Big)^2\, dx
\]
so that $a$ is the squared $\mathcal{L}^2$-distance between the restriction of $g$ to $[0,1]$ and the set of restrictions of affine functions to $[0,1]$.

\begin{lemma} Let $g$ be uniformly continuous, let the errors satisfy \eqref{e:lim=gamma} and \eqref{e:ass_distr_tail}, let $\hat\gamma$ be an estimator satisfying (G\ref{i:G_alternative}) and let the design points satisfy \eqref{e:ass_scatter_alt} and (A\ref{i:ass_scatter_converge}). Then
\[\frac{2T}{n} \to a \]
as $n\to\infty$ in distribution. 
\end{lemma}
\prf One can show
\[a = \int_0^1 g(x)^2\, \mu(dx) - \Big( \int_0^1 g(x) \, \mu(dx) \Big)^2 -  \mathcal{C}  \Big(\int_0^1  g(x) (x-\mathcal{X}_0) \, \mu(dx)\Big)^2,\]
where
\[ \mathcal{X}_0= \int_0^1 x\, \mu(dx)\quad \mbox{and} \quad \mathcal{C} = \frac{1}{ \int_0^1 x^2\, \mu(dx) - (\int_0^1 x\, \mu(dx))^2}, \]
the same way as Proposition \ref{p:T1_explicit} is shown. Hence 
\begin{align*}
\mathbb{P}\bigg(&\big|\frac{2T}{n}-a\big| > \epsilon \bigg) \\
&\le \begin{aligned}[t] \mathbb{P}\bigg(&\Big|\frac{2}{n}\sum_{i=1}^{n/2} \hat g(x_{2i-1})^2 + \frac{4}{\hat\gamma n}\sum_{i=1}^{n/2} Y_{2i-1}\mathbf{1}_{\{Y_{2i-1}\ge \hat g(x_{2i-1})\}} - \int_0^1 g(x)^2\, \mu(dx) \Big| > \frac{\epsilon}{3} \bigg)\\
&+ \mathbb{P}\bigg(\Big|\frac{4}{n^2} \Big( \sum_{i=1}^{n/2} \hat g(x_{2i-1})+ \frac{1}{\hat\gamma}\sum_{i=1}^{n/2} \mathbf{1}_{\{Y_{2i-1}\ge \hat g(x_{2i-1})\}}\Big)^2 - \Big(\int_0^1 g(x)\, \mu(dx) \Big)^2\Big| >\frac{\epsilon}{3} \bigg)\\
&+ \mathbb{P}\bigg( \Big|\frac{2}{n} \frac{ n\Big(\sum_{i=1}^{n/2} \big(\hat g(x_{2i-1}) + \frac{1}{\hat\gamma} \mathbf{1}_{\{Y_{2i-1} \ge \hat g(x_{2i-1})\}} \big) \cdot (x_{2i-1}-\tfrac{S}{n/2}) \Big)^2 }{ Rn-2S^2} - \mathcal{C} \Big(\int_0^1 g(x)\big(x-\mathcal{X}_0\big)\, \mu(dx) \Big)^2\Big| >\frac{\epsilon}{3} \bigg)
\end{aligned}
\\
& =: A_I + A_{II} + A_{III}.
\end{align*}
Let $n$ be large enough that $|\frac{2}{n}\sum_{i=1}^{n/2} g(x_{2i-1})^2 - \int_0^1 g(x)^2\, \mu(dx)| < \epsilon/6$. Then
\begin{align*}
A_I &\le   \mathbb{P}\Big(\big|\frac{2}{n}\sum_{i=1}^{n/2} \big(\hat g(x_{2i-1})^2 - g(x_{2i-1})^2\big) + \frac{4}{\hat\gamma n}\sum_{i=1}^{n/2} Y_{2i-1}\mathbf{1}_{\{Y_{2i-1}\ge \hat g(x_{2i-1})\}}  \big| > \frac{\epsilon}{6} \Big)\\
& \le \frac{ \Ex\Big[ \big|\frac{2}{n}\sum_{i=1}^{n/2} \big(\hat g(x_{2i-1})^2 - g(x_{2i-1})^2\big) + \frac{4}{\hat\gamma n}\sum_{i=1}^{n/2} Y_{2i-1}\mathbf{1}_{\{Y_{2i-1}\ge \hat g(x_{2i-1})\}}  \big| \Big] }{ \frac{\epsilon}{6} } \\
& \le \frac{ \frac{2}{n}\sum_{i=1}^{n/2}  \Ex\big[\big|\hat g(x_{2i-1}) - g(x_{2i-1})\big| \cdot  2|g(x_{2i-1})| + \big(\hat g(x_{2i-1}) - g(x_{2i-1})\big)^2\big] + \Ex\big[\frac{4}{\hat\gamma n}\sum_{i=1}^{n/2}  |Y_{2i-1}| \mathbf{1}_{\{Y_{2i-1}\ge \hat g(x_{2i-1})\}} \big] }{ \frac{\epsilon}{6} } \\
& \longrightarrow 0 
\end{align*}
by Lemma \ref{l:asy_ex_alt} and the Cauchy-Schwarz inequality, since 
\[\Ex\big[ |Y_{2i-1}|^2 \mathbf{1}_{\{Y_{2i-1}\ge \hat g(x_{2i-1})\}} \mid \hat g \big]\le \max\{|g(x_{2i-1})|^2, |\hat g(x_{2i-1})|^2\} \big( \bar\gamma |\hat g(x_{2i-1}) - g(x_{2i-1})| \big), \]
where $\bar\gamma$ is a constant with $1-F(t)\le \bar\gamma|t|$, $t\le 0$. 
Now let $n$ be large enough that $|\frac{2}{n}\sum_{i=1}^{n/2} g(x_{2i-1}) - \int_0^1 g(x)\, \mu(dx)| < \sqrt{\epsilon/6}$. Then 
\begin{align*}
A_{II} &\le \begin{aligned}[t] \mathbb{P}\bigg(&\frac{4}{n^2} \Big( \sum_{i=1}^{n/2} \hat g(x_{2i-1})+ \frac{1}{\hat\gamma}\sum_{i=1}^{n/2} \mathbf{1}_{\{Y_{2i-1}\ge \hat g(x_{2i-1})\}}- \sum_{i=1}^{n/2} g(x_{2i-1})\Big)^2\\
& + \frac{8}{n^2} \Big| \sum_{i=1}^{n/2} \hat g(x_{2i-1})+ \frac{1}{\hat\gamma}\sum_{i=1}^{n/2} \mathbf{1}_{\{Y_{2i-1}\ge \hat g(x_{2i-1})\}}- \sum_{i=1}^{n/2} g(x_{2i-1})\Big|\cdot\Big|\sum_{i=1}^{n/2} g(x_{2i-1})\Big| >\frac{\epsilon}{6} \bigg)\end{aligned} \\
& \le \mathbb{P}\bigg(\frac{2}{n} \bigg( \sum_{i=1}^{n/2} \hat g(x_{2i-1})+ \frac{1}{\hat\gamma}\sum_{i=1}^{n/2} \mathbf{1}_{\{Y_{2i-1}\ge \hat g(x_{2i-1})\}}- \sum_{i=1}^{n/2} g(x_{2i-1})\Big) > \min\Big\{ \sqrt{\frac{\epsilon}{12}}, \frac{\epsilon}{48/n\cdot |\sum_{i=1}^{n/2} g(x_{2i-1})|} \Big\} \bigg)\\
& \le \frac{\frac{2}{n} \sum_{i=1}^{n/2} \Ex\big[| \hat g(x_{2i-1})+ \frac{1}{\hat\gamma}\mathbf{1}_{\{Y_{2i-1}\ge \hat g(x_{2i-1})\}} - g(x_{2i-1})|\big]}{\min\big\{ \sqrt{\frac{\epsilon}{12}}, \frac{\epsilon}{48/n\cdot |\sum_{i=1}^{n/2} g(x_{2i-1})|} \big\}}\\
& \le  \frac{ \sup\big\{\Ex\big[ |\hat g(x)- g(x)| \big] + \sqrt{\Ex[\frac{1}{\hat\gamma^2}]} \sqrt{ \Ex\big[ \bar\gamma |\hat g(x) -g(x)| \big]} \mid x\in [0,1] \big\}}{\min\big\{ \sqrt{\frac{\epsilon}{12}}, \frac{\epsilon}{48/n\cdot |\sum_{i=1}^{n/2} g(x_{2i-1})|} \big\}}\\
& \longrightarrow 0. 
\end{align*}
In order to treat $A_{III}$ assume that $|\frac{2}{n}\sum_{i=1}^{n/2} g(x_{2i-1})(x_{2i-1}-\frac{S}{n/2}) - \int_0^1 g(x)(x-\mathcal{X}_0)\, \mu(dx)| < \epsilon\cdot (Rn-2S^2)/9$ and 
\[\Big(\mathcal{C} - \frac{n^2/2}{Rn-2S^2} \Big) \cdot \int_0^1 g(x)\big(x-\mathcal{X}_0\big) \, \mu(dx) < \frac{\epsilon}{9}. \]
 Then
\begin{align*}
 A_{III} & \le \begin{aligned}[t] \mathbb{P}\bigg(&\frac{4}{n^2} \Big( \sum_{i=1}^{n/2} \big(\hat g(x_{2i-1})+ \frac{1}{\hat\gamma} \mathbf{1}_{\{Y_{2i-1}\ge \hat g(x_{2i-1})\}}-  g(x_{2i-1}) \big) \cdot \big( x_{2i-1} - \mathcal{X}_0\big) \Big)^2\\
& + \frac{8}{n^2} \Big| \sum_{i=1}^{n/2} \big( \hat g(x_{2i-1})+ \frac{1}{\hat\gamma} \mathbf{1}_{\{Y_{2i-1}\ge \hat g(x_{2i-1})\}}-  g(x_{2i-1})\big) \cdot \big(x_{2i-1}-\mathcal{X}_0\big)\Big|\cdot\Big|\sum_{i=1}^{n/2} g(x_{2i-1}) \cdot (x_{2i-1}-\mathcal{X}_0 \big)\Big| >\frac{\epsilon}{9} \bigg)\end{aligned} 
\\
&\le  \frac{ \frac{n^2}{Rn-2S^2} \cdot \frac{1}{2} \sup\big\{\Ex\big[ \hat g(x)- g(x) \big] + \sqrt{\Ex[\frac{1}{\hat\gamma^2}]} \sqrt{\Ex\big[ \bar\gamma |\hat g(x) -g(x)|  \big]} \mid x\in [0,1] \big\}}{\min\big\{ \sqrt{\frac{\epsilon}{18}}, \frac{\epsilon}{72/n\cdot |\sum_{i=1}^{n/2} g(x_{2i-1})|} \big\}} \longrightarrow 0\end{align*}
can be shown the same way as the corresponding estimate for $A_{II}$.  \qed\medskip

\prf[ of Theorem \ref{t:consistent}]  The critical value $c_{\lambda,n}$ of either test fulfills $c_{\lambda,n}<(na)/4$ with a probability tending to $1$. So
 \[ \lim_{n\to\infty} \mathbb{P}(\varphi=1) \ge \lim_{n\to\infty} \mathbb{P}\big( \tfrac{T}{n} >\tfrac{a}{4}\big) =1. \qquad \qed \]

\subsection{Properties of the scale parameter estimator}\label{s:negHill}

In this subsection we investigate the properties of the estimator $\hat\gamma$ from \eqref{e:negHill}.  

\begin{lemma}\label{l:4moment}
If \eqref{e:lim=gamma}, \eqref{e:ass_distr_tail}, \eqref{e:ass_scatter}, \eqref{e:ass1}, \eqref{e:ass2} and \eqref{e:ass2a} hold, then (G\ref{i:G4}) holds. 
\end{lemma}

The proof of this lemma is based on several further lemmata. 

\begin{lemma}\label{l:order_st}
Let $U_{1:n}, \dots, U_{n:n}$ be the ascendingly sorted order statistics of $n$ independent random variables that are uniformly distributed on $[0,1]$. Then the density of $U_{k:n}$ is bounded by
\[  \frac{20}{\pi} n \sqrt{\tfrac{n}{(k-1)(n-k)}} \exp\Big\{ -\Big(\frac{(n-1)^2}{6(k-1)} + \frac{(n-1)^2}{6(n-k)} \Big) \big(x-\frac{k-1}{n-1}\big)^2 \Big\} \]
if $|x-\frac{k-1}{n-1}| \le (1- \sqrt[3]{\frac{1}{2}}) \min\{ \frac{k-1}{n-1}, \frac{n-k}{n-1} \}$ and 
\[  \frac{20}{\pi} n \sqrt{\tfrac{n}{(k-1)(n-k)}} \exp\Big\{ -\Big(\frac{(n-1)^2}{6(k-1)} + \frac{(n-1)^2}{6(n-k)} \Big) c \Big\} \]
if $|x-\frac{k-1}{n-1}| \ge c$ for some $c\le (1- \sqrt[3]{\frac{1}{2}}) \min\{ \frac{k-1}{n-1}, \frac{n-k}{n-1} \}$.
\end{lemma}

\begin{lemma}\label{l:stirling_up}
We have
\[ k! \le 10 \big( \frac{k}{e} \big)^k \sqrt{k}. \]
\end{lemma}
\prf The same way as in the proof of Lemma \ref{l:stirling} we get
\[ k! = k^{k+1} \int_0^\infty e^{k(\log y - y)} \, dy.\]
Now a Taylor series expansion of $f(y):= \log y -y$ yields
\[ f(y)=f(1)+f'(1)\cdot (y-1)+ \frac{1}{2} f''(\tilde y) \cdot (y-1)^2 \le -1- \frac{1}{2} \frac{(y-1)^2}{4} \]
for all $y\in[0,2]$ and some appropriate $\tilde y$ between $1$ and $y$. For $y>2$ we have $f(y)\le -y/2$. Hence
\begin{align*}
\int_0^\infty e^{k(\log y -y)} \, dy &\le \int_0^2 e^{k(-1-(y-1)^2/8)}\, dy + \int_2^\infty e^{-ky/2}\, dy\\
&\le e^{-k} \cdot \int_{-\infty}^\infty e^{-k(y-1)^2/8}\, dy + \int_2^\infty e^{-ky/2}\, dy\\
&= e^{-k} \sqrt{\frac{2\pi \cdot 8}{k}} + \tfrac{2}{k}e^{-k} \\
&\le e^{-k} \frac{8}{\sqrt{k}} +\frac{2}{\sqrt{k}}e^{-k} \\
& = \frac{10}{\sqrt{k}} e^{-k}.
\end{align*}
Hence the assertion follows. \qed \medskip

\prf[ of Lemma \ref{l:order_st}]
It is well-known that the density of $U_{k:n}$ is
\[ \frac{n!}{(k-1)!(n-k)!} x^{k-1}(1-x)^{n-k}. \]
This can be bounded from above by
\begin{align*}
&\frac{ 10 \big(\frac{n}{e}\big)^n \sqrt{n} }{ \sqrt{\frac{\pi(k-1)}{2}}\big(\frac{k-1}{e}\big)^{k-1} \sqrt{\frac{\pi(n-k)}{2}} \big(\frac{n-k}{e}\big)^{n-k}} x^{k-1} (1-x)^{n-k} \\
&= \frac{20}{\pi} \frac{n}{e} \sqrt{\frac{n}{(k-1)(n-k)}} \big( \tfrac{n}{n-1}\big)^{n-1} \Big(\frac{x}{\frac{k-1}{n-1}}\Big)^{k-1} \Big( \frac{1-x}{\frac{n-k}{n-1}} \Big)^{n-k}\\
&= \frac{20}{\pi} \frac{n}{e} \sqrt{\frac{n}{(k-1)(n-k)}} \big( \tfrac{n}{n-1}\big)^{n-1} \exp\Big\{ (k-1) \cdot \big(\log(x)-\log\big(\tfrac{k-1}{n-1}\big)\big) + (n-k) \cdot \big(\log(1-x)-\log\big(\tfrac{n-k}{n-1}\big) \big) \Big\}.
\end{align*}
Using a Taylor expansion of $\log$ about $\frac{k-1}{n-1}$ we get
\[ \log x = \log\big(\tfrac{k-1}{n-1}\big) + \frac{x-\frac{k-1}{n-1}}{\frac{k-1}{n-1}} - \tfrac{1}{2} \frac{\big(x-\frac{k-1}{n-1}\big)^2}{\big(\frac{k-1}{n-1}\big)^2} + \tfrac{1}{6} \frac{2\big(x-\frac{k-1}{n-1}\big)^3}{\tilde x^3} \]
for some $\tilde x$ between $\frac{k-1}{n-1}$ and $x$.  Using a Taylor expansion about $\frac{n-k}{n-1}$ we get
\[ \log(1-x) = \log\big(\tfrac{n-k}{n-1}\big) + \frac{1-x-\frac{n-k}{n-1}}{\frac{n-k}{n-1}} - \tfrac{1}{2} \frac{\big((1-x)-\frac{n-k}{n-1}\big)^2}{\big(\frac{n-k}{n-1}\big)^2} + \tfrac{1}{6} \frac{2\big(1-x-\frac{n-k}{n-1}\big)^3}{\hat x^3} \]
for some $\hat x$ between $\frac{n-k}{n-1}$ and $1-x$. For $x$ with
\begin{equation}  \big| x - \frac{k-1}{n-1} \big| \le \big(1-\sqrt[3]{\tfrac{1}{2}}\big) \min\Big\{ \frac{k-1}{n-1}, \frac{n-k}{n-1} \Big\} \label{e:x_near} \end{equation}
we thus get 
\begin{align*}
 (k-1) \cdot \big(&\log(x)-\log\big(\tfrac{k-1}{n-1}\big)\big) + (n-k) \cdot \big(\log(1-x)-\log\big(\frac{n-k}{n-1}\big) \big) \\
&= \begin{aligned}[t] (k-1) \cdot \Big(& \frac{x-\frac{k-1}{n-1}}{\frac{k-1}{n-1}} - \tfrac{1}{2} \frac{\big(x-\frac{k-1}{n-1}\big)^2}{\big(\frac{k-1}{n-1}\big)^2} + \tfrac{1}{6} \frac{2\big(x-\frac{k-1}{n-1}\big)^3}{\tilde x^3} \Big)\\ & + (n-k) \cdot  \Big(\frac{1-x-\frac{n-k}{n-1}}{\frac{n-k}{n-1}} - \tfrac{1}{2} \frac{\big((1-x)-\frac{n-k}{n-1}\big)^2}{\big(\frac{n-k}{n-1}\big)^2} + \tfrac{1}{6} \frac{2\big(1-x-\frac{n-k}{n-1}\big)^3}{\hat x^3}\Big)\end{aligned}\\
&\le \begin{aligned}[t] (k-1) \cdot \frac{x-\frac{k-1}{n-1}}{\frac{k-1}{n-1}} + (n-k) \cdot \frac{1-x-\frac{n-k}{n-1}}{\frac{n-k}{n-1}} &+(k-1)\cdot \Big( -\tfrac{1}{2} \frac{\big(x-\frac{k-1}{n-1}\big)^2}{\big(\frac{k-1}{n-1}\big)^2} +\tfrac{1}{6} \frac{2\big(x-\frac{k-1}{n-1}\big)^2\frac{1}{2} \frac{k-1}{n-1}}{\frac{1}{2}\big(\frac{k-1}{n-1}\big)^3} \Big) \\ &+ (n-k) \cdot \Big( -\tfrac{1}{2} \frac{\big((1-x)-\frac{n-k}{n-1}\big)^2}{\big(\frac{n-k}{n-1}\big)^2} + \tfrac{1}{6} \frac{2\big((1-x)-\frac{n-k}{n-1}\big)^2\frac{1}{2} \frac{n-k}{n-1}}{\frac{1}{2}\big(\frac{n-k}{n-1}\big)^3} \Big)\end{aligned} \\
&= -\frac{k-1}{6}  \frac{\big(x-\frac{k-1}{n-1}\big)^2}{\big(\frac{k-1}{n-1}\big)^2} -\frac{n-k}{6} \frac{\big((1-x)-\frac{n-k}{n-1}\big)^2}{\big(\frac{n-k}{n-1}\big)^2}\\
&= -\Big( \frac{(n-1)^2}{6(k-1)} + \frac{(n-1)^2}{6(n-k)} \Big) \cdot \big(x- \frac{k-1}{n-1}\big)^2.
\end{align*}
Hence for $x$ satisfying \eqref{e:x_near} the density of $U_{k:n}$ can be bounded from above by
\[ \frac{20}{\pi} n\sqrt{\tfrac{n}{(n-1)(n-k)}} \exp\Big\{ - \Big( \frac{(n-1)^2}{6(k-1)} + \frac{(n-1)^2}{6(n-k)} \Big) \big(x - \frac{k-1}{n-1}\big)^2 \Big\}. \]
Furthermore, observing
\[ \frac{d}{dx}\big( (k-1) \cdot \log x + (n-k)\cdot \log(1-x) \big) = \frac{k-1}{x} - \frac{n-k}{1-x} \]
we see that $(k-1) \cdot \log x + (n-k)\cdot \log(1-x)$ is monotonically increasing on $[0, \frac{k-1}{n-1}]$ and monotonically decreasing on $[\frac{k-1}{n-1}, 1]$. Thus the density of $U_{k:n}$ is bounded from above by 
\[ \frac{20}{\pi} n\sqrt{\tfrac{n}{(n-1)(n-k)}} \exp\Big\{ - \Big( \frac{(n-1)^2}{6(k-1)} + \frac{(n-1)^2}{6(n-k)} \Big) c \Big\} \]
for all $x\in[0,1]$ with $|x-\frac{k-1}{n-1}|\ge c$ if $c\le \big(1-\sqrt[3]{\tfrac{1}{2}}\big) \min\Big\{ \frac{k-1}{n-1}, \frac{n-k}{n-1} \Big\}$. 
Hence the assertion follows. \qed

\begin{lemma}\label{l:order_to_zero}
Assume that \eqref{e:lim=gamma} holds, let $(k_n)_{n\in\mathbb{N}}$ be a sequence with $\lim_{n\to\infty} k=\infty$ and $\lim_{n\to\infty} k/n=0$ and let $(a_n)_{n\in\mathbb{N}}$ be a sequence fulfilling $\lim_{n\to\infty} a_n/\sqrt{k}=0$ and $\lim_{n\to\infty} a_nk/n=0$. 
Then
\[ a_n \cdot \big(\frac{n}{2k} \epsilon_{\frac{n}{2}-k:\frac{n}{2}} + \frac{1}{\gamma } \big) \longrightarrow 0 \]
in probability as $n\to\infty$. 
\end{lemma}
A central limit theorem for $\epsilon_{\frac{n}{2}-k:\frac{n}{2}}$ is given in \cite[Theorem 2.2.1]{dHFe06} under the additional assumption that the distribution function $F$ is twice differentiable. Showing a central limit theorem for $\epsilon_{\frac{n}{2}-k:\frac{n}{2}}$ without differentiability assumptions on $F$ would be an interesting project, but since this lemma is not a central result of this paper, we will not do it here. While the proof of this lemma uses some ideas of the proof of \cite[Theorem 2.2.1]{dHFe06}, it is essentially new.

\prf For a monotonically increasing function $f$ denote by $f^\leftarrow(y):=\inf\{s \mid f(s)>y\}$ its right-continuous inverse. Then $\epsilon_{2i}=F^\leftarrow(1-U_i)$ in distribution for $i=1, \dots, n/2$, where the random variables $U_i$ are distributed uniformly on $[0,1]$. We assume that $\epsilon_{2i}$ and $U_i$ are defined on the same probability space and that $\epsilon_{2i}=F^\leftarrow(1-U_i)$ holds almost surely. Then we have $\epsilon_{\frac{n}{2}-k:\frac{n}{2}} = F^\leftarrow(1-U_{k+1:\frac{n}{2}})$. 

Now for each $\epsilon>0$ there is $\delta>0$ such that
\begin{equation} \frac{y-1}{\gamma} - \frac{C_F+\epsilon}{\gamma} \big( \frac{y-1}{\gamma}\big)^2  \le  F^\leftarrow(y)  \le  \frac{y-1}{\gamma} + \frac{C_F}{\gamma} \big( \frac{y-1}{\gamma}\big)^2 \label{e:F_inv_ungl} \end{equation}
for all $y\in(1-\delta,1)$. 

Indeed, let $\epsilon>0$. Choose $\eta>0$ such that $\frac{C_F+\epsilon}{\gamma^3}> \frac{C_F}{\gamma} (\frac{1}{\gamma}+\eta)^2$. Further choose $\delta>0$ such that
\[ \frac{C_F}{\gamma} \Big( \big(\frac{1}{\gamma}+\eta\big)(y-1) \Big)^2< \eta \cdot |y-1| \qquad \mbox{and} \qquad C_F\Big( \big(\frac{1}{\gamma}+\eta\big)\big(y-1\big) \Big)^2 < |y-1| \]
for all $y\in(1-\delta,1)$ and let $y\in(1-\delta,1)$. 

In order to prove the first inequality of \eqref{e:F_inv_ungl}, let $s\le \frac{y-1}{\gamma} - \frac{C_F+\epsilon}{\gamma}\big(\frac{y-1}{\gamma}\big)^2$. Then, in particular, $s \le \frac{y-1}{\gamma} - \frac{C_F}{\gamma} \big(\frac{1}{\gamma}+\eta)^2(y-1)^2$ and hence
\begin{align*}
F(s) & \le F\big( \frac{y-1}{\gamma} - \frac{C_F}{\gamma}(\frac{1}{\gamma}+\eta)^2(y-1)^2 \big) \\
&\le 1  +  \gamma \cdot \big( \frac{y-1}{\gamma} - \frac{C_F}{\gamma}(\frac{1}{\gamma}+\eta)^2(y-1)^2 \big) + C_F \cdot \big( \frac{y-1}{\gamma} - \frac{C_F}{\gamma}(\frac{1}{\gamma}+\eta)^2(y-1)^2 \big)^2 \\
& \le 1 +  {y-1}- C_F(\frac{1}{\gamma}+\eta)^2(y-1)^2  + C_F \cdot \big( \frac{y-1}{\gamma} - \eta \cdot |y-1| \big)^2\\
&= y,
\end{align*}
since $x\mapsto (\frac{y-1}{\gamma}-x)^2$ is monotonically increasing on $(\frac{y-1}{\gamma},\infty)$.

Therefore we get
\[ \frac{y-1}{\gamma} - \frac{C_F+\epsilon}{\gamma} \big(\frac{y-1}{\gamma}\big)^2 \le \inf\{s \mid F(s)>y\} = F^\leftarrow(y). \]

In order to show the second inequality of \eqref{e:F_inv_ungl}, observe
\begin{align*}
F\big( \frac{y-1}{\gamma} + \frac{C_F}{\gamma} \big( \frac{y-1}{\gamma}\big)^2 \big)  & \ge   1 + \gamma\cdot\big( \frac{y-1}{\gamma} + \frac{C_F}{\gamma} \big( \frac{y-1}{\gamma}\big)^2 \big)  -  C_F\cdot \big( \frac{y-1}{\gamma} + \frac{C_F}{\gamma} \big( \frac{y-1}{\gamma}\big)^2 \big)^2 \\
& \ge 1 + y-1 + C_F \big( \frac{y-1}{\gamma}\big)^2	  -  C_F\cdot \big( \frac{y-1}{\gamma} \big)^2 \\
&=y, 
\end{align*} 
since $x \mapsto C_F \cdot (\frac{y-1}{\gamma} +x )^2$ is monotonically decreasing $(-\infty, \frac{|y-1|}{\gamma})$.
Hence 
\[\frac{y-1}{\gamma} + \frac{C_F}{\gamma} (\frac{y-1}{\gamma})^2 \ge \inf\{s \mid F(s)>y\} = F^\leftarrow(y). \]
So \eqref{e:F_inv_ungl} is proven. 

Now let $\epsilon>0$ and choose $\delta>0$ such that \eqref{e:F_inv_ungl} holds. Then we have
\[ -\frac{U_{k+1:\frac{n}{2}}}{\gamma} - \frac{C_F+\epsilon}{\gamma} \big( \frac{U_{k+1:\frac{n}{2}}}{\gamma}\big)^2  \le  \epsilon_{\frac{n}{2}-k:\frac{n}{2}}  \le  -\frac{U_{k+1:\frac{n}{2}}}{\gamma} + \frac{C_F}{\gamma} \big( \frac{U_{k+1:\frac{n}{2}}}{\gamma}\big)^2 \]
provided that $U_{k+1:n/2}<\delta$ and in particular with a probability tending to $1$. 
Since 
\[ \sqrt{k} \big( \frac{n}{2k} U_{k+1:\frac{n}{2}} -1 \big) \stackrel{d}{\longrightarrow} \mathcal{N}(0,1) \]
in distribution as $n\to\infty$ as a consequence of Smirnov's lemma (see e.g. \cite[Lemma 2.2.3]{dHFe06}), we get
\[ a_n \big( \frac{n}{2k} U_{k+1:\frac{n}{2}} -1 \big) \stackrel{P}{\longrightarrow} 0. \] 
Since, moreover, 
\[ a_n \big( \frac{n}{2k} U_{k+1:\frac{n}{2}}^2  \big) \stackrel{P}{\longrightarrow} 0 \]
we get 
\[ a_n \big( \frac{n}{2k} \big(-\frac{U_{k+1:\frac{n}{2}}}{\gamma} - \frac{C_F+\epsilon}{\gamma} \big( \frac{U_{k+1:\frac{n}{2}}}{\gamma}\big)^2 \big) -1 \big) \stackrel{P}{\longrightarrow} 0 \] 
and
\[ a_n \big( \frac{n}{2k} \big(-\frac{U_{k+1:\frac{n}{2}}}{\gamma} + \frac{C_F}{\gamma} \big( \frac{U_{k+1:\frac{n}{2}}}{\gamma}\big)^2 \big)-1 \big) \stackrel{P}{\longrightarrow} 0. \]  
Hence the assertion follows. \qed

\prf[ of Lemma \ref{l:4moment}]
We have
\begin{align*} 
\mathbb{E}\big[ \big(\tfrac{1}{\hat\gamma} - \tfrac{1}{\gamma} \big)^4 \big] &= \mathbb{E}\big[ \big(\tfrac{n}{2k}(\hat\epsilon_{\frac{n}{2}:\frac{n}{2}}-\hat\epsilon_{\frac{n}{2}-k:\frac{n}{2}}) - \tfrac{1}{\gamma} \big)^4 \big] \\
&\le 81\mathbb{E}\big[ \big(\tfrac{n}{2k}(\epsilon_{\frac{n}{2}:\frac{n}{2}}-\epsilon_{\frac{n}{2}-k:\frac{n}{2}}) - \tfrac{1}{\gamma} \big)^4 \big] + 81\mathbb{E}\big[ \big(\tfrac{n}{2k}(\hat\epsilon_{\frac{n}{2}-k:\frac{n}{2}}-\epsilon_{\frac{n}{2}-k:\frac{n}{2}})  \big)^4 \big] + 81\mathbb{E}\big[ \big(\tfrac{n}{2k}(\hat\epsilon_{\frac{n}{2}:\frac{n}{2}}-\epsilon_{\frac{n}{2}:\frac{n}{2}}) \big)^4 \big].
\end{align*}
Now Proposition \ref{p:JMR3.1} and Proposition \ref{p:asy_ex} yield
\begin{align*}
\mathbb{E}\big[ & \big(\tfrac{n}{2k}(\hat\epsilon_{\frac{n}{2}-k:\frac{n}{2}}-\epsilon_{\frac{n}{2}-k:\frac{n}{2}})  \big)^4 \big] \\
& \le \big( \tfrac{n}{2k} \big)^4 \mathbb{E}\big[ \sup\big\{ |\hat g(x) - g(x)| \mid x\in[0,1] \big\}^4 \big] \\
& \le \big( \tfrac{n}{2k} \big)^4 \mathbb{E}\big[ \max\big\{ \big|\max\{\epsilon_{2i} \mid x_{2i}\in ((j-1)\tfrac{h_1}{2}, j\tfrac{h_1}{2} ] \}\big| \mid j=0, \dots, \lceil \tfrac{2}{h_1} \rceil+1 \big\}^4 \big] \\ 
& \le \big( \tfrac{n}{2k} \big)^4 \sum_{j=0}^{\lceil \frac{2}{h_1} \rceil+1} \mathbb{E}\big[ \max\{\epsilon_{2i} \mid x_{2i}\in ((j-1)\tfrac{h_1}{2}, j\tfrac{h_1}{2} ] \}^4 \big]\\
& \in O\big( \big( \tfrac{n}{k} \big)^4  \tfrac{1}{h_1} (nh_1)^{-4}\big) \\
&= O\big( k^{-4}h_1^{-5} \big)\\
&\subseteq O\big(h^{-2}n^{-2}\big)
\end{align*}
by \eqref{e:ass1}.
The same way one obtains
\[ \mathbb{E}\big[ \big(\tfrac{n}{2k}(\hat\epsilon_{\frac{n}{2}:\frac{n}{2}}-\epsilon_{\frac{n}{2}:\frac{n}{2}}) \big)^4 \big] \in O\big(h^{-2}n^{-2}\big). \]
So let us turn to 
\[ \mathbb{E}\big[ \big(\tfrac{n}{2k}(\epsilon_{\frac{n}{2}:\frac{n}{2}}-\epsilon_{\frac{n}{2}-k:\frac{n}{2}}) - \tfrac{1}{\gamma} \big)^4 \big]. \]
Due to \eqref{e:ass2} and \eqref{e:ass2a} we can apply Lemma \ref{l:order_to_zero} with $a_n=\sqrt{nh}$. Hence
\[ \sqrt{nh} \big( \tfrac{n}{2k}\epsilon_{\frac{n}{2}-k:\frac{n}{2}} + \tfrac{1}{\gamma} \big) \stackrel{P}{\longrightarrow} 0  \]
and together with
\[ \sqrt{nh} \big( \tfrac{n}{2k}\epsilon_{\frac{n}{2}:\frac{n}{2}} \big) \stackrel{P}{\longrightarrow} 0, \]
which follows from \eqref{e:ass2} and the convergence of $n\epsilon_{\frac{n}{2}:\frac{n}{2}}$ to the Weibull distribution, 
we get 
\[ \sqrt{nh} \big( \tfrac{n}{2k}(\epsilon_{\frac{n}{2}:\frac{n}{2}}-\epsilon_{\frac{n}{2}-k:\frac{n}{2}}) -\tfrac{1}{\gamma}\big) \stackrel{P}{\longrightarrow} 0. \]
So it suffices to show
\[ \int_0^\infty \sup_n \mathbb{P}\big( \sqrt{nh} \big| \tfrac{n}{2k}(\epsilon_{\frac{n}{2}:\frac{n}{2}}-\epsilon_{\frac{n}{2}-k:\frac{n}{2}}) -\frac{1}{\gamma}\big|  \ge S \big) S^3 dS <\infty.\]
We have
\begin{align*}
\mathbb{P}\big( \sqrt{nh} \big| \frac{n}{2k}(\epsilon_{\frac{n}{2}:\frac{n}{2}}-&\epsilon_{\frac{n}{2}-k:\frac{n}{2}}) -\frac{1}{\gamma}\big|  \ge S \big) = \mathbb{P} \big( \big| \frac{n}{2k} \big( \epsilon_{\frac{n}{2}:\frac{n}{2}} - \epsilon_{\frac{n}{2}-k:\frac{n}{2}} - \frac{2k}{n\gamma}\big) \big| \ge \frac{S}{\sqrt{nh}} \big) \\
& \le \mathbb{P} \big(  \epsilon_{\frac{n}{2}-k:\frac{n}{2}}  \le -\frac{2k}{n\gamma} - \frac{Sk}{n\sqrt{nh}} \big) + \mathbb{P} \big(  \epsilon_{\frac{n}{2}-k:\frac{n}{2}}  \ge -\frac{2k}{n\gamma}+ \frac{Sk}{n\sqrt{nh}} \big) + \mathbb{P} \big( \frac{n}{2k} \epsilon_{\frac{n}{2}:\frac{n}{2}} \le \frac{S}{2\sqrt{nh}} \big). \end{align*}
Consider at first the case 
\[ \frac{\gamma Sk}{2n\sqrt{nh}} \le \big(1- \sqrt[3]{\tfrac{1}{2}}\big) \min\Big\{ \frac{k-1}{n/2-1}, \frac{n/2-k}{n/2-1} \Big\}. \] 
From \eqref{e:ass2} and \eqref{e:ass2a} we get 
\[ \frac{k}{\sqrt{nh}}>1, \quad  C_F \cdot \frac{8k^2}{n^2\gamma^2} \le \frac{k}{n\sqrt{nh}}\quad \mbox{and} \quad C_F\frac{k-1}{n/2-1} \le \frac{\gamma^2}{8} \]
 for all sufficiently large $n$, which will be assumed in the sequel. Let $S>4/\gamma$. Then 
\begin{align*}
 \mathbb{P} \big(  \epsilon_{\frac{n}{2}-k:\frac{n}{2}}  \le -\frac{2k}{n\gamma} - \frac{Sk}{n\sqrt{nh}} \big) &\le  \mathbb{P} \big(  F(\epsilon_{\frac{n}{2}-k:\frac{n}{2}})  \le F(-\frac{2k}{n\gamma} - \frac{Sk}{n\sqrt{nh}}) \big) \\
&\le  \mathbb{P} \big(  F(\epsilon_{\frac{n}{2}-k:\frac{n}{2}})  \le 1-\gamma\cdot(\frac{2k}{n\gamma} + \frac{Sk}{n\sqrt{nh}}) + C_F \cdot(\frac{2k}{n\gamma} + \frac{Sk}{n\sqrt{nh}})^2 \big) \\
&\le  \mathbb{P} \big(  F(\epsilon_{\frac{n}{2}-k:\frac{n}{2}})  \le 1-\frac{k-1}{n/2-1}  - \frac{\gamma Sk}{n\sqrt{nh}} + \frac{k}{n\sqrt{nh}} + \frac{\gamma Sk}{4n\sqrt{nh}} \big) \\
&\le  \mathbb{P} \big(  F(\epsilon_{\frac{n}{2}-k:\frac{n}{2}})  \le 1-\frac{k-1}{n/2-1}  - \frac{\gamma Sk}{2n\sqrt{nh}} \big)
\end{align*}
and, if
\[ \frac{Sk}{n \sqrt{nh}} \le \frac{2k}{n\gamma}, \]
then 
\begin{align*}
 \mathbb{P} \big(  \epsilon_{\frac{n}{2}-k:\frac{n}{2}}  \ge -\frac{2k}{n\gamma} + \frac{Sk}{n\sqrt{nh}} \big) &\le  \mathbb{P} \big(  F(\epsilon_{\frac{n}{2}-k:\frac{n}{2}})  \ge F(-\frac{2k}{n\gamma} + \frac{Sk}{n\sqrt{nh}}) \big) \\
&\le  \mathbb{P} \big(  F(\epsilon_{\frac{n}{2}-k:\frac{n}{2}})  \ge 1-\gamma\cdot(\frac{2k}{n\gamma} - \frac{Sk}{n\sqrt{nh}}) - C_F \cdot(\frac{2k}{n\gamma} - \frac{Sk}{n\sqrt{nh}})^2 \big) \\
&\le  \mathbb{P} \big(  F(\epsilon_{\frac{n}{2}-k:\frac{n}{2}})  \ge 1-\frac{k-1}{n/2-1} -\frac{1}{n/2-1} + \frac{\gamma Sk}{n\sqrt{nh}} - \frac{k}{n\sqrt{nh}}  \big) \\
&\le  \mathbb{P} \big(  F(\epsilon_{\frac{n}{2}-k:\frac{n}{2}})  \ge 1-\frac{k-1}{n/2-1}  + \frac{\gamma Sk}{2n\sqrt{nh}} \big).
\end{align*}
Notice that 
 \[  \mathbb{P} \big(  \epsilon_{\frac{n}{2}-k:\frac{n}{2}}  \ge -\frac{2k}{n\gamma} + \frac{Sk}{n\sqrt{nh}} \big) \le \mathbb{P} \big(  F(\epsilon_{\frac{n}{2}-k:\frac{n}{2}})  \ge 1-\frac{k-1}{n/2-1}  + \frac{\gamma Sk}{2n\sqrt{nh}} \big) \]
holds trivially if
\[ \frac{Sk}{n \sqrt{nh}} \ge \frac{2k}{n\gamma}. \]
Assume $2\le k<n/4$ and $\frac{(n/2-1)^2}{(n/2)^2} > \frac{24}{25}$. Then Lemma \ref{l:order_st} yields
\begin{align*}
 \mathbb{P} \Big( \big|  \epsilon_{\frac{n}{2}-k:\frac{n}{2}} +\frac{2k}{n\gamma} \big| \ge &\frac{ Sk}{n\sqrt{nh}} \Big) \\
&\le \frac{20}{\pi} \frac{n}{2} \sqrt{\tfrac{n/2}{(k-1)(n/2-k)}} \exp\Big\{ -\big(\frac{(n/2-1)^2}{6(k-1)} + \frac{(n/2-1)^2}{6(n/2-k)} \big) \big(\frac{\gamma Sk}{2n\sqrt{nh}}\big)^2 \Big\} \Big(1-\frac{\gamma Sk}{2n\sqrt{nh}}\Big)^+\\
& \le \frac{10}{\pi} n \sqrt{\tfrac{4}{k}} \exp\Big\{ -\frac{n^2}{25k}  \big(\frac{4k}{2n\sqrt{nh}}\big)^2 \Big\} \big(\frac{2n\sqrt{nh}}{\gamma Sk}\big)^5\\
& = \frac{640}{\pi\gamma^5} \exp\Big\{ \tfrac{17}{2}\log n + \tfrac{5}{2} \log h -\tfrac{11}{2} \log k - \frac{4k}{25nh} \Big\} S^{-5}.
\end{align*}
Notice that by \eqref{e:ass2} we have
\[ B_1:= \sup_n \exp\Big\{ \tfrac{17}{2}\log n + \tfrac{5}{2} \log h -\tfrac{11}{2} \log k - \frac{4k}{25nh} \Big\} < \infty. \]
Now assume
\[ \frac{\gamma Sk}{2n\sqrt{nh}} \ge \big(1- \sqrt[3]{\tfrac{1}{2}}\big) \min\Big\{ \frac{k-1}{n/2-1}, \frac{n/2-k}{n/2-1} \Big\}. \] 
Put $N:=\lceil 5/\Gamma \rceil$ and $X_j:=\max\{\epsilon_{(j-1)N+1}, \dots, \epsilon_{jN}\}$ for $j=1,\dots, l:= \lfloor \frac{n}{2N} \rfloor$ and let $G$ denote the distribution function of the $X_j$.
Then, if $l>k$, 
\begin{align*}
 \mathbb{P} \Big(  \epsilon_{\frac{n}{2}-k:\frac{n}{2}}  \le -\frac{2k}{n\gamma} - \frac{Sk}{n\sqrt{nh}} \Big) \le   \mathbb{P} \Big(  X_{l-k:l}  \le -\frac{2k}{n\gamma} - \frac{Sk}{n\sqrt{nh}} \Big)
&\le \mathbb{P} \Big( G( X_{l-k:l})  \le G( -\frac{2k}{n\gamma} - \frac{Sk}{n\sqrt{nh}}) \Big)
\end{align*}
and
\[ G\big(-\frac{2k}{n\gamma} - \frac{Sk}{n\sqrt{nh}}\big) \le F\big(-\frac{Sk}{n\sqrt{nh}} \big)^{5/\Gamma} \le \tilde C_F^{5/\Gamma} \cdot \big( \frac{n\sqrt{nh}}{Sk}\big)^5. \]
There is $\delta>0$ such that $F(t) \le 1+\frac{23}{24}\gamma t$ and $\exp\{t\}\le 1+ \frac{22}{23}t$ for $t\in(-\delta,0)$. Let $n$ be large enough that
\[ \frac{6(k-1)}{5(n/2-1)\gamma} \le \delta, \quad \frac{23}{20}N\frac{k-1}{n/2-1} \le \delta, \quad \mbox{and} \quad \frac{(n/2-1)/N}{l-1} \le \frac{22}{21}. \]
Then
\[G\Big(-\frac{2k}{n\gamma} - \frac{Sk}{n\sqrt{nh}}\Big) \le F\Big(-\frac{k-1}{(n/2-1)\gamma} - \frac{1}{5} \frac{k-1}{(n/2-1)\gamma} \Big)^N \le \Big(1-\frac{23}{20}\frac{k-1}{n/2-1} \Big)^N \le \exp\Big\{ -\frac{23}{20}N\frac{k-1}{n/2-1} \Big\} \le 1-\frac{21}{20} \frac{k-1}{l-1}. \]
Lemma \ref{l:order_st} yields
\begin{align*}
 \mathbb{P} \Big(&  \epsilon_{\frac{n}{2}-k:\frac{n}{2}} +\frac{2k}{n\gamma} \le -\frac{Sk}{n\sqrt{nh}} \Big) \\
&\le \mathbb{P} \Big( G( X_{l-k:l})  \le \min\big\{1-\frac{21}{20} \frac{k-1}{l-1}, \tilde C_F^{5/\Gamma} \cdot \big( \frac{n\sqrt{nh}}{Sk}\big)^5 \big\}\Big)\\
&=\frac{20}{\pi} l \sqrt{\tfrac{l}{(k-1)(l-k)}} \exp\Big\{ -\big(\frac{(l-1)^2}{6(k-1)} + \frac{(l-1)^2}{6(l-k)} \big) \big(\tfrac{1}{20}\big)^2\big( \frac{k-1}{l-1} \big)^2 \Big\} \tilde C_F^{5/\Gamma}\big(\frac{n\sqrt{nh}}{Sk}\big)^5\\
& \le \frac{20}{\pi} l \sqrt{\frac{4}{k}} \exp\Big\{ -\frac{k-1}{6N} \big(\tfrac{1}{20}\big)^2 \Big\} \tilde C_F^{5/\Gamma}\big(\frac{n\sqrt{nh}}{Sk}\big)^5\\
& = \frac{40 \tilde C_F^{5/\Gamma}}{\pi} \exp\Big\{ \tfrac{15}{2}\log n + \log l + \tfrac{5}{2} \log h -\tfrac{11}{2} \log k - \frac{k-1}{2400} \Big\} S^{-5}.
\end{align*}
Notice that \eqref{e:ass2} in particular implies $\limsup n^\alpha/k=0$ and hence
\[ B_2:= \sup_n  \exp\Big\{ \tfrac{15}{2}\log n + \log l + \tfrac{5}{2} \log h -\tfrac{11}{2} \log k - \frac{k-1}{2400} \Big\} < \infty. \]
Moreover, the same way as in the first case we get
\[
 \mathbb{P} \big(  \epsilon_{\frac{n}{2}-k:\frac{n}{2}}  \ge -\frac{2k}{n\gamma} + \frac{Sk}{n\sqrt{nh}} \big) \le \mathbb{P} \big(  F(\epsilon_{\frac{n}{2}-k:\frac{n}{2}})  \ge 1-\frac{k-1}{n/2-1}  - \frac{\gamma Sk}{2n\sqrt{nh}} \big)
\]
and thus
\begin{align*}
 \mathbb{P} \Big(&  \epsilon_{\frac{n}{2}-k:\frac{n}{2}} +\frac{2k}{n\gamma} \ge \frac{Sk}{n\sqrt{nh}} \Big) \\
&\le \frac{20}{\pi} \frac{n}{2} \sqrt{\tfrac{n/2}{(k-1)(n/2-k)}} \exp\Big\{ -\big(\frac{(n/2-1)^2}{6(k-1)} + \frac{(n/2-1)^2}{6(n/2-k)} \big) \big(1-\sqrt[3]{\frac{1}{2}}\big)^2 \min\big\{\frac{k-1}{n/2-1}, \frac{n/2-k}{n/2-1}\big\}^2 \Big\} \Big(1-\frac{\gamma Sk}{2n\sqrt{nh}}\Big)^+\\
& \le \frac{10}{\pi} n \sqrt{\tfrac{4}{k}} \exp\Big\{ -\frac{(n/2-1)^2}{6(k-1)} (\tfrac{1}{5})^2 \big(\frac{k-1}{n/2-1}\big)^2 \Big\} \big(\frac{2n\sqrt{nh}}{\gamma Sk}\big)^5\\
& \le \frac{10}{\pi} n \sqrt{\tfrac{4}{k}} \exp\Big\{ -\frac{k-1}{6} (\tfrac{1}{5})^2  \Big\} \big(\frac{2n\sqrt{nh}}{\gamma Sk}\big)^5\\
& = \frac{640}{\pi\gamma^5} \exp\Big\{ \tfrac{17}{2}\log n + \tfrac{5}{2} \log h -\tfrac{11}{2} \log k - \frac{k-1}{150} \Big\} S^{-5}.
\end{align*}
Now
\[ B_3 := \sup_n \exp\Big\{ \tfrac{17}{2}\log n + \tfrac{5}{2} \log h -\tfrac{11}{2} \log k - \frac{k-1}{150} \Big\} < \infty. \] 
Observe that \eqref{e:lim=gamma} implies that there is some $s<0$ with $F(t) \le \exp\{\gamma t/2\}$ for all $t\in[s,0]$. For $n$ that are large enough that $k>\sqrt{nh}$ we get
\begin{align*}
\mathbb{P} \big( \big| \tfrac{n}{2k} \epsilon_{\frac{n}{2}:\frac{n}{2}} \big| \ge \frac{S}{2\sqrt{nh}} \big) &= F\big(-\frac{Sk}{n\sqrt{nh}}\big)^{n/2} \\
&\le  \exp\big\{ - \frac{\gamma Sk}{2n\sqrt{nh}} \big\}^{n/2} +  F\big(\frac{Sk}{n\sqrt{nh}}\big)^{5/\Gamma}\cdot F(s)^{n/2-5/\Gamma} \\
& \le \exp\Big\{ - \gamma\frac{kS}{4\sqrt{nh}} \Big\} + \tilde C_F^{5/\Gamma} \big(\frac{n\sqrt{nh}}{Sk}\big)^5\cdot F(s)^{n/2-5/\Gamma} \\
& \le \exp\{-\gamma S/4\} + \tilde C_F^{5/\Gamma} \big(\frac{n\sqrt{nh}}{k}\big)^5\cdot F(s)^{n/2-5/\Gamma} \cdot S^{-5}.
\end{align*}
Observe
\[ B_4:= \sup_n \big(\frac{n\sqrt{nh}}{k}\big)^5\cdot F(s)^{n/2-5/\Gamma} < \infty.\]
Altogether, there is some $n_0$ such that
\begin{align*}
\int_0^\infty \sup_{n>n_0}& \mathbb{P}\big( \sqrt{nh} \big| \tfrac{n}{2k}(\epsilon_{\frac{n}{2}:\frac{n}{2}}-\epsilon_{\frac{n}{2}-k:\frac{n}{2}}) -\frac{1}{\gamma}\big|  \ge S \big) S^3 \, dS\\
& \le \int_0^{4/\gamma} 1 \cdot S^3\, dS + \int_{4/\gamma}^\infty \Big( \frac{640}{\pi \gamma^5} B_1 S^{-5} + \frac{40 \tilde C_F^{5/\Gamma}}{\pi} B_2 S^{-5} + \frac{640}{\pi \gamma^5} B_3 S^{-5} + \exp\{-\gamma S/4\} + \tilde C_F^{5/\Gamma} B_4S^{-5} \Big) S^3 \, dS < \infty. \qquad \qed
\end{align*}

\begin{lemma} Assume that \eqref{e:ass_distr_tail}, \eqref{e:ass3} and \eqref{e:ass4} hold. Then (G\ref{e:gamma_large-dev}) holds as well. 
\end{lemma}
\prf Put $N:=\lceil \frac{9}{2\Gamma} \rceil$ and $X_j:=\max\{\epsilon_{(j-1)N+1}, \dots, \epsilon_{jN}\}$ for $j=1, \dots, l:=\lfloor \frac{n}{2N} \rfloor$. We have
\[
\mathbb{P}(\hat \gamma<t) = \mathbb{P}\big( \tfrac{n}{2k} (\hat\epsilon_{\frac{n}{2}:\frac{n}{2}} - \hat\epsilon_{\frac{n}{2}-k:\frac{n}{2}}) >\tfrac{1}{t} \big)\le \mathbb{P}\big( \epsilon_{\frac{n}{2}-k:\frac{n}{2}} <-\tfrac{2k}{nt} \big)\le \mathbb{P}\big( X_{l-k:l}< - \tfrac{2k}{nt} \big).
\]
Let $G$ denote the distribution function of the $X_j$ and let $U_{1:l}, \dots, U_{l:l}$ be the ascendingly sorted order statistics of independent random variables $U_1, \dots, U_l$ that are distributed uniformly on $[0,1]$. Then  
\[ \mathbb{P}\big( X_{l-k:l} < -\frac{2k}{nt} \big) \le \mathbb{P}\big( G(X_{l-k:l}) \le G(-\tfrac{2k}{nt}) \big) \le \mathbb{P}\big( U_{l-k:l} \le G(-\tfrac{2k}{nt}) \big). \] 
Now we will show
\begin{equation} \Big|G\big(-\tfrac{2k}{nt}\big) - \frac{n-k}{n-1} \Big|\ge (1-\sqrt[3]{\tfrac{1}{2}}) \cdot \frac{k-1}{n-1} \label{e:G_near} \end{equation}
for all $t\le \bar\gamma$ and all sufficiently large $n$. Here $\bar\gamma$ is a number such that $F(x) \le 1 + \bar\gamma x$ for all $x\in(-\delta, 0)$ for some $\delta>0$. Assume $\frac{2k-2}{n-1} \le \bar\gamma\delta$. In case $\frac{2k}{nt}<\delta$ we have
\[ G\big(-\tfrac{2k}{nt}\big)\le F\big(-\tfrac{2k}{nt}\big) \le 1-\bar\gamma \frac{2k}{nt} \le 1 - \frac{2k}{n} = \frac{n-2k}{n} \le \frac{n-2k+1}{n-1}. \]
In case $\frac{2k}{nt} \ge \delta$ we have
\[ G\big(-\tfrac{2k}{nt}\big)\le F(-\delta) \le 1-\bar\gamma \delta\le \frac{n-2k+1}{n-1}.\]
In both cases we obtain \eqref{e:G_near}. So we get
\begin{align*}
\mathbb{P}(\hat \gamma<t) & \le \mathbb{P}\big( U_{l-k:l} \le G(-\tfrac{2k}{nt}) \big)\\
&\le \frac{20}{\pi} l \sqrt{\tfrac{l}{(k-1)(l-k)}} \exp\Big\{ -\Big(\frac{(l-1)^2}{6(k-1)} + \frac{(l-1)^2}{6(l-k)} \Big) \big(1- \sqrt[3]{\tfrac{1}{2}}\big)^2 \min\Big\{ \frac{k-1}{l-1}, \frac{l-k}{l-1} \Big\}^2 \Big\} G\Big(-\frac{2k}{nt}\Big)\\
&\le \frac{20}{\pi} l \sqrt{\tfrac{4}{k}} \exp\Big\{ -\frac{(l-1)^2}{6(k-1)}  (\tfrac{1}{5})^2 \big( \frac{k-1}{l-1} \big)^2 \Big\} F\Big(-\frac{2k}{nt}\Big)^N\\
&\le \frac{20}{\pi} l \sqrt{\tfrac{4}{k}} \exp\Big\{ -\frac{k-1}{150}   \Big\} \tilde C_F^{9/(2\Gamma)} \cdot \Big(\frac{2k}{nt}\Big)^{-9/2}\\
&\le \frac{5\tilde C_F^{9/(2\Gamma)}}{2^{3/2}\pi} \exp\Big\{ \tfrac{11}{2} \log n -\log(2 N) -5 \log k -\frac{k-1}{150}\Big\} t^{9/2}.
\end{align*}
Since
\[ \lim_{n\to\infty} \tfrac{11}{2} \log n -5 \log k -\frac{k-1}{150} = -\infty \]
in view of \eqref{e:ass4}, we get
\[ B_5 := \sup_n \exp\Big\{ \tfrac{11}{2} \log n -\log (2N) -5 \log k -\frac{k-1}{150}\Big\}<\infty. \]
Hence we have
\[ \mathbb{P}(\hat\gamma<t) \le \max\Big\{ \frac{5\tilde C_F^{9/(2\Gamma)}}{2^{3/2}\pi} B_5, \bar\gamma^{-9/2} \Big\} t^{9/2} \]
for all $t>0$ and all $n$ large enough that $k<l/2$. \qed

\begin{lemma} Let the function $g:[0,1]\to\mathbb{R}$ be uniformly continuous, assume that the design points are equidistant, $x_i=i/n, i=1, \dots, n$, and assume \eqref{e:ass4}, \eqref{e:ass5}, \eqref{e:ass6} and \eqref{e:ass8}. Then (G\ref{i:G_alternative}) holds. \end{lemma}
\prf We have 
\begin{align*} \Ex\big[ (\frac{1}{\hat \gamma})^2\big] &= \Ex\big[ (\frac{n}{2k})^2 (\hat\epsilon_{\frac{n}{2}:\frac{n}{2}} - \hat\epsilon_{\frac{n}{2}-k:\frac{n}{2}})^2 \big] \\
&\le  9\Ex\big[ (\frac{n}{2k})^2 (\hat\epsilon_{\frac{n}{2}:\frac{n}{2}} - \epsilon_{\frac{n}{2}:\frac{n}{2}})^2 \big] + 9\Ex\big[ (\frac{n}{2k})^2 (\epsilon_{\frac{n}{2}:\frac{n}{2}} - \epsilon_{\frac{n}{2}-k:\frac{n}{2}})^2 \big] + 9\Ex\big[ (\frac{n}{2k})^2 (\epsilon_{\frac{n}{2}-k:\frac{n}{2}} - \hat\epsilon_{\frac{n}{2}-k:\frac{n}{2}})^2 \big] \\
&\le 18 \Ex\big[(\frac{n}{2k})^2 \sup\big\{ |\hat g(x)-g(x)| \mid x\in[0,1]\big\}^2\big] + 9\Ex\big[ (\frac{n}{2k})^2 ( \epsilon_{\frac{n}{2}-k:\frac{n}{2}})^2 \big]. 
\end{align*}
Put $Z_j:=\max\{\epsilon_{2i} \mid x_{2i} \in ((j-1)\tfrac{h_1}{2},j\tfrac{h_1}{2}) \}$ for $j=0,\dots,\lceil \frac{2}{h_1}\rceil+1$ and let
\[\omega_g(\delta) := \sup\{ |g(x_1)-g(x_2)| \mid |x_1-x_2| \le \delta \} \]
be the modulus of continuity of $g$. Then
\[ |\hat g(x)-g(x)| \le \omega_g(h_1)+\max\{|Z_j|\mid j=0,\dots,\lceil \frac{2}{h_1}\rceil+1\} \]
and hence
\[ \Ex\big[ \sup\{|\hat g(x) -g(x)| \mid x\in[0,1] \}^2\big] \le \sum_{j=0}^{\lceil \frac{2}{h_1} \rceil+1} 2\Ex[|Z_j|^2] + 2\omega_g(h_1)^2. \] 
Now $\frac{n}{2k} \sum Z_j\to 0$ in probability as $n\to \infty$ by \eqref{e:ass6}, since $nh_1Z_0$ converges to a Weibull distribution. This implies $\Ex[(\frac{n}{2k})^2 \sum |Z_j|^2] \to 0$ as $n\to\infty$, since
\[ \int_0^\infty \sup_n \mathbb{P}\big((\lceil \tfrac{2}{h_1} \rceil +2)\frac{n}{2k}|Z_0|>t\big) t \, dt  < \infty, \]
where the supremum is taken over all $n$ such that $\lfloor nh_1/2\rfloor\Gamma \ge 3$. Indeed, let $s<0$ be a number with $F(t) \le \exp\{ \gamma t/2\}$ for $t\in[s,0]$. Then
\begin{align*}
\int_0^\infty \sup_n \mathbb{P}((\lceil \tfrac{2}{h_1} \rceil +2)\frac{n}{2k}|Z_0|>t) t \, dt  &\le \int_0^\infty \sup_n \mathbb{P}(Z_0<-\frac{2kh_1}{3n}t) t \, dt\\
&\le  \int_0^\infty \sup_n F\big(-\frac{2kh_1}{3n}t\big)^{\lfloor nh_1/2 \rfloor} t \, dt\\
& \le \int_0^\infty \sup_n \big\{ \exp\{-\gamma \frac{kh_1}{3n} t \lfloor nh_1/2 \rfloor \} + \tilde C_F^{3/\Gamma} (\frac{3n}{2kh_1})^3 t^{-3} F(s)^{\lfloor nh_1/2 \rfloor-3/\Gamma} \big\} t\, dt < \infty \end{align*}
in view of \eqref{e:ass6} and \eqref{e:ass8}. Now \eqref{e:ass5} implies 
\[ \Ex\big[(\frac{n}{2k})^2 \sup\big\{ |\hat g(x)-g(x)| \mid x\in[0,1]\big\}^2\big] \to 0.\]
Now Lemma \ref{l:order_to_zero}, applied with $a_n=1$, yields
\[ \frac{n}{2k} \epsilon_{\frac{n}{2}-k:\frac{n}{2}} \stackrel{P}{\longrightarrow} \frac{1}{\gamma}, \qquad n\to\infty. \]
So it remains to prove
\[ \int_0^\infty \sup_n \mathbb{P}\big( \frac{n}{2k} \epsilon_{\frac{n}{2}-k:\frac{n}{2}} \le -t\big) t \, dt <\infty. \]
Recall that by \eqref{e:ass_distr_tail} there are $\Gamma>0$ and $\tilde C_F>0$ with $F(t)\le \tilde C_F\cdot |t|^{-\Gamma}$ for $t<0$. Put $N:=\lceil 3/\Gamma \rceil$ and $X_j:=\max\{\epsilon_{j(N-1)+1}, \dots, \epsilon_{jN}\}$ for $j=1, \dots, l:=\lfloor \frac{n}{2N} \rfloor$ and let $G$ denote the distribution function of the $X_j$. There is $\delta>0$ such that $F(t)\le 1+ \frac{\gamma}{2}t$ for all $t\in(-\delta,0)$. Assume that $n$ is large enough that $8k/(n\gamma)<\delta$, $k<l$ and 
\[\exp\{-2\frac{k-1}{l-1}\}\le 1- \frac{6}{5} \frac{k-1}{l-1}.\]
Then for all $t\ge 4/\gamma$ we have
\begin{align*}
\mathbb{P}\Big( \frac{n}{2k} \epsilon_{\frac{n}{2}-k:\frac{n}{2}} \le -t\big) &\le \mathbb{P}\Big(  \epsilon_{l-k:l} \le -\frac{2kt}{n}\Big) \\
&\le \mathbb{P}\Big(  G(\epsilon_{l-k:l}) \le F(-\frac{2kt}{n})^N\Big) \\
&\le \mathbb{P}\Big( G( \epsilon_{l-k:l}) \le \min\big\{(1-\frac{\gamma}{2}\frac{8k}{n\gamma})^N, \tilde C_F^{3/\Gamma}(\frac{2kt}{n})^{-3}\big\}\Big) \\ 
&\le \mathbb{P}\Big( G( \epsilon_{l-k:l}) \le \min\big\{\exp\{-2\frac{k-1}{l-1}\}, \tilde C_F^{3/\Gamma}(\frac{2kt}{n})^{-3}\big\}\Big)\\
&\le \frac{20}{\pi} l \sqrt{\frac{l}{(k-1)(l-k)}} \exp\Big\{ -\Big( \frac{(l-1)^2}{6(k-1)} + \frac{(l-1)^2}{6(l-k)} \Big) \big(\tfrac{1}{5}\big)^2 \big( \frac{k-1}{l-1} \big)^2 \Big\} \tilde C_F^{3/\Gamma}(\frac{2kt}{n})^{-3}\\
&\le \frac{5}{2\pi} l \sqrt{\frac{4}{k}} \exp\Big\{ - \frac{(l-1)^2}{6(k-1)} \big(\tfrac{1}{5}\big)^2 \big( \frac{k-1}{l-1}\big)^2 \Big\} \tilde C_F^{3/\Gamma}(\frac{kt}{n})^{-3}\\
&= \frac{5}{\pi} l \sqrt{\frac{1}{k}} \exp\Big\{ - \frac{k-1}{150}  \Big\} \tilde C_F^{3/\Gamma}(\frac{kt}{n})^{-3}\\
&= \frac{5\tilde C_F^{3/\Gamma}}{\pi} \exp\Big\{ \log l + 3 \log n - \frac{7}{2} \log k - \frac{k-1}{150}  \Big\} t^{-3}.
\end{align*}
In view of \eqref{e:ass4} we get
\[ B_6:= \sup_n \exp\big\{ \log l + 3 \log n - \frac{7}{2} \log k - \frac{k-1}{150} \big\} < \infty. \]
Hence
\[ \int_0^\infty \sup_n \mathbb{P}\big( \frac{n}{2k} \epsilon_{\frac{n}{2}-k:\frac{n}{2}} \le -t\big) t \, dt \le \int_0^{4/\gamma} t\, dt + \int_{4/\gamma}^\infty \frac{5\tilde C_F^{3/\Gamma}}{\pi} B_6 t^{-3} t \, dt <\infty. \qquad \qed \]
\end{appendix}


\begin{thebibliography}{99}
\bibitem{Bil99} P.~Billingsley: \emph{Convergence of Probability Measures}. Wiley (1999). 
\bibitem{CFR07} J.~Cuesta-Albertos, R.~Fraiman and T.~Ransford: \emph{A sharp form of the Cr\'amer-Wold theorem}. Journal of Theoretical Probability {\bf 20} (2007), 201--209.
\bibitem{DLF19} C.~Dong, G.~Li and X.~Feng: \emph{Lack-of-fit tests for quantile regression models}. Journal of the Royal Statistical Society Series B {\bf 81} (2019), 629--648.  
\bibitem{DNS17} H.~Drees, N.~Neumeyer and L.~Selk: \emph{Estimation and hypotheses testing in boundary regression problems}. Bernoulli {\bf 25} (2019), 424--463.
\bibitem{Fa95} M.~Falk: \emph{Some best parameter estimates for distributions with finite endpoint}. Statistics: A Journal of Theoretical and Applied Statistics {\bf 27} (1995), 115--125.  
\bibitem{GC13} W.~Gonzalez-Manteiga and R.M.~Crujeiras: \emph{An updated review of goodness-of-fit tests for regression models}. TEST {\bf 22} (2013), 361-–411.
\bibitem{JMR14} M.~Jirak, A.~Meister and M.~Rei\ss: \emph{Adaptive function estimation in nonparametric regression with one-sided errors}. The Annals of Statistics {\bf 42} (2014), 1970--2002. 
\bibitem{dHFe06} L.~d.~Haan and A.~Ferreira: \emph{Extreme Value Theory -- An Introduction}. Springer (2006). 
\bibitem{HvK09} P.~Hall and I.~v.~Keilegom: \emph{Nonparametric ``regression'' when errors are positioned at end-points}. Bernoulli {\bf 15} (2009), 614--633.
\bibitem{HuKu15} S.~Huet and E.~Kuhn: Goodness-of-fit test for Gaussian regression with block correlated errors. Statistics {\bf 49} (2015), 239--266.
\bibitem{HMD} Human Mortality Database (www.mortality.org). University of California at Berkeley and the Max Planck Institute
for Demographic Research.
\bibitem{MvKY19} E.~Mammen, I.~v.~Keilegom and K.~Yu: \emph{Expansion for moments of regression quantiles with applications to nonparametric testing}. Bernoulli {\bf 25} (2019), 793--827. 
\bibitem{OV02} J.~Oeppen and J.~Vaupel: \emph{Broken limits to life expectancy}. Science {\bf 296} (2002), 1029--1031.
\bibitem{Ph69} W.~Philipp: \emph{The central limit theorem for mixing sequences of random variables}. Zeitschrift f\"ur Wahrscheinlichkeitstheorie und verwandte Gebiete {\bf 12} (1969), 155--171.
\bibitem{ReWa17} M.~Rei{\ss} and M.~Wahl: \emph{Functional estimation and hypotheses testing in nonparametric boundary models}. Bernoulli {\bf 25} (2019), 2597--2619.   
\bibitem{Res87} S.~Resnick: \emph{Extreme Values, Regular Variation, and Point Processes}. Springer (1987).
\bibitem{Zh98} J.~Zheng: \emph{A consistent nonparametric test of parametric regression models under conditional quantile restrictions}. Economic Theory {\bf 14} (1998), 123--138.
\end{thebibliography}
\end{document}